\documentclass[12pt, a4paper, twoside]{amsart}
\usepackage{amsmath}
\usepackage{amssymb}
\usepackage{amsfonts}
\usepackage{a4}
\usepackage[all]{xy}

\usepackage[latin1]{inputenc}
\usepackage[latin1]{inputenc}

\numberwithin{equation}{section}

\newtheorem{theorem}{Theorem}[section]
\newtheorem{lemma}[theorem]{Lemma}

\theoremstyle{definition} 

\theoremstyle{plain} 

\newtheorem{proposition}[theorem]{Proposition}

\newtheorem{conclusion}[theorem]{Conclusion}
\newtheorem{convention}[theorem]{Convention}
\newtheorem{claim}[theorem]{Claim}

\newtheorem*{maintheorem*}{Main Theorem} 
\newtheorem*{conjecture*}{Conjecture} 
\newtheorem{definition}[theorem]{Definition}

\newtheorem{hypothesis}[theorem]{Hypothesis}
\newtheorem{observation}[theorem]{Observation}

\theoremstyle{remark}  
\newtheorem{remark}[theorem]{Remark}
\newtheorem*{remarks*}{Remarks}
\newtheorem*{remark*}{Remark}
\newtheorem*{claim*}{Claim}

\newcommand{\nc}{\newcommand}
\nc{\nothing}[1]{}

\nothing{   
\nc{\dom}{{\rm dom}}
\nc{\card}{{\rm card}}
\nc{\lh}{{\rm lh}}
\nc{\lgg}{{\rm lg}}
\nc{\rge}{\mbox{\rm range}}
\nc{\cf}{{\rm cf}}
\nc{\nex}{\mbox{\rm next}}
\nc{\uhr}{\restriction}
\nc{\supt}{{\rm supt}}
\nc{\supp}{{\rm supp}}
\nc{\Lim}{{\rm Lim}}
\nc{\Leb}{{\rm Leb}}
\nc{\modd}{{\rm mod}}
\nc{\RO}{{\rm RO}}
\nc{\prob}{{\rm Prob}}
}

\nc{\On}{{\rm On}}
\nc{\Ord}{{\rm On}}
\nc{\nco}{\DeclareMathOperator}
\nco{\rk}{rk}
\nco{\order}{o} 
\nco{\ppower}{pp} 
\nco{\pcf}{pcf} 
\nco{\tcf}{tcf} 
\nco{\tlim}{tlim} 
\nco{\limtext}{lim} 
\nco{\prodt}{{\textstyle \prod}} 
\nco{\symdiff}{\triangle}
\nco{\dom}{dom}
\nco{\card}{card}
\nco{\lh}{lh}
\nco{\lgg}{lg}
\nco{\hgt}{ht}
\nco{\rge}{range}
\nco{\otp}{otp}
\nco{\trunk}{tr}
\nco{\cf}{cf}
\nco{\nex}{next}
\nc{\uhr}{\restriction}
\nco{\reduction}{red}
\nco{\supt}{supt}
\nco{\supp}{supp}
\nco{\Lim}{Lim}
\nco{\Leb}{Leb}
\nco{\modd}{mod}
\nco{\invariant}{inv}
\nco{\id}{id}
\nco{\RO}{RO}
\nco{\poss}{pos}
\nco{\Inc}{Inc} 
\nco{\Ge}{Ge}
\nco{\hdrop}{\hat{drop}}

\nc{\potom}{\ensuremath{{\cal P}(\omega)}}
\nc{\potinf}{\ensuremath{[\omega]^\omega}}
\nc{\pfin}{\ensuremath{{\cal P}(\omega)/{\rm fin}}}

\nc{\potfin}{\ensuremath{[\omega]^{<\omega}}}
\nc{\inn}{\ensuremath{{\omega^{\uparrow \omega}}}}
\nc{\hoch}{^{<\omega}}
\nc{\hocho}{^{\omega}}
\nc{\tree}[1]{{[} #1 {]}_0}
\nc{\tre}[2]{ {#1}_{#2}}
   
\nc{\prooff}[1]{{\bf Proof} of #1:}
\nc{\proofend}{\makebox{} \hfill ${\bf \square}$ \\}
\nc{\proofendof}[1]{\makebox{} \hfill $\boldmath{\square}_{\rm #1}$ \\}
\nc{\beq}{\begin{eqnarray*}}
\nc{\eeq}{\end{eqnarray*}}
\nc{\bde}{\begin{list}}
\nc{\ede}{\end{list}}

\newenvironment{myrules}
{\begin{list}{}
{
 \setlength{\leftmargin}{0.5in}
 \setlength{\labelwidth}{1cm}
 \setlength{\labelsep}{0.2in}
 \setlength{\parsep}{0.5ex plus 0.2ex minus 0.1 ex}
 \setlength{\itemsep}{0.3ex plus 0.2 ex minus 0ex}
}}{\end{list}}

{\end{list}}

\newcounter{subalph}
{\end{list}}

\newcommand{\greek}[1]{\ifthenelse{\value{#1}=1}{\mbox{$\alpha$}}%
  {\ifthenelse{\value{#1}=2}{\mbox{$\beta$}}{%
   \ifthenelse{\value{#1}=3}{\mbox{$\gamma$}}{%
   \ifthenelse{\value{#1}=4}{\mbox{$\delta$}}{%
   \ifthenelse{\value{#1}=5}{\mbox{$\varepsilon$}}{%
   \ifthenelse{\value{#1}=6}{\mbox{$\zeta$}}{%
   \ifthenelse{\value{#1}=7}{\mbox{$\eta$}}{%
   \ifthenelse{\value{#1}=8}{\mbox{$\theta$}}{%
   \ifthenelse{\value{#1}=9}{\mbox{$\iota$}}{%
   \ifthenelse{\value{#1}=10}{\mbox{$\kappa$}}{%
   \ifthenelse{\value{#1}=11}{\mbox{$\lambda$}}{%
   \ifthenelse{\value{#1}=12}{\mbox{$\mu$}}{%
   \ifthenelse{\value{#1}=13}{\mbox{$\nu$}}{%
   \ifthenelse{\value{#1}=14}{\mbox{$\xi$}}{%
   \ifthenelse{\value{#1}=15}{\mbox{$\rm o$}}{%
   \ifthenelse{\value{#1}=16}{\mbox{$\pi$}}{%
   \ifthenelse{\value{#1}=17}{\mbox{$\varrho$}}{%
   \ifthenelse{\value{#1}=18}{\mbox{$\sigma$}}{%
   \ifthenelse{\value{#1}=19}{\mbox{$\tau$}}{%
   \ifthenelse{\value{#1}=20}{\mbox{$\upsilon$}}{%
   \ifthenelse{\value{#1}=21}{\mbox{$\varphi$}}{%
   \ifthenelse{\value{#1}=22}{\mbox{$\chi$}}{%
   \ifthenelse{\value{#1}=23}{\mbox{$\psi$}}{\mbox{$\omega$}%
  }}}}}}}}}}}}}}}}}}}}}}}}

\newcounter{subgreek}
{\end{list}}

\newcounter{subarabic}
{\end{list}}

\newcounter{subroman}
{\end{list}}

\newcount\skewfactor

\def\mathunderaccent#1#2 {\let\theaccent#1\skewfactor#2
\mathpalette\putaccentunder}
\def\putaccentunder#1#2{\oalign{$#1#2$\crcr\hidewidth
\vbox to.2ex{\hbox{$#1\skew\skewfactor\theaccent{}$}\vss}\hidewidth}}
\def\name{\mathunderaccent\tilde-3 }

\nc{\nname}{\name}

\nc{\even}{\ensuremath{\rm Even}}
\nc{\odd}{\ensuremath{\rm Odd}}

\nc{\al}{$\alpha$\  }
\nc{\om}{\omega}
\nc{\omm}{\ensuremath{\omega_1}}
\nc{\ep}{\varepsilon}
\nc{\tk}{\tilde{K}}
\nc{\concat}{\,\hat{} \,}   
\nc{\force}{\Vdash}
\nc{\fb}{f_{\overline{M}}}
\nc{\such}{\, : \,}   
\newcommand{\la}{\langle}
\renewcommand{\r}{\rangle}

\nc{\meager}{\ensuremath{{\cal M}}}
\nc{\lebesgue}{\ensuremath{{\cal N}}}
\nc{\nulll}{\ensuremath{{\cal N}}}
\nc{\ksigma}{\ensuremath{{\bf K}_\sigma}}
\nc{\ideal}{\ensuremath{{\cal I}}}
\nc{\ga}{\ensuremath{\frak a}}
\nc{\AAA}{{\cal A}}   
\nc{\gc}{\ensuremath{\frak c}}
\nc{\gs}{\ensuremath{\frak s}}
\nc{\gh}{\ensuremath{\frak h}}
\nc{\gd}{\ensuremath{\frak d}}
\nc{\gb}{\ensuremath{\frak b}}
\nc{\gro}{\ensuremath{\frak g}}
\nc{\gu}{\ensuremath{\frak u}} 
\nc{\gr}{\ensuremath{\frak r}} 
\nc{\gt}{\ensuremath{\frak t}}
\nc{\fff}{\ensuremath{\frak f}}
\nc{\gm}{\ensuremath{\mathfrak{mcf}}}
\nc{\gge}{\ensuremath{\mathfrak e}}
\nc{\cfupro}{\ensuremath{\cf(\upro)}}
\nc{\cfvpro}{\ensuremath{\cf(\vpro)}}
\nc{\gp}{\ensuremath{\frak p}}
\nc{\gk}{\ensuremath{\frak k}}

\nc{\add}[1]{\mbox{\ensuremath{{\rm add}(#1)}}}
\nc{\cov}[1]{\mbox{\ensuremath{{\rm cov}(#1)}}}
\nc{\unif}[1]{\mbox{\ensuremath{{\rm unif}(#1)}}}
\nc{\cof}[1]{{\mbox{\ensuremath{\rm cof}(#1)}}}

\nc{\addd}[2]{\mbox{\ensuremath{{\rm add}^{#1}(#2)}}}   
\nc{\covv}[2]{\mbox{\ensuremath{{\rm cov}^{#1}(#2)}}}   
\nc{\uniff}[2]{\mbox{\ensuremath{{\rm unif}^{#1}(#2)}}} 
\nc{\coff}[2]{{\mbox{\ensuremath{\rm cof}^{#1}(#2)}}}

\nc{\cd}{Cicho\'n's Diagram}

\nc{\MA}{\mbox{\rm MA}}
\nc{\PFA}{\mbox{\rm PFA}}
\nc{\OCA}{\mbox{\rm OCA}}
\nc{\GCH}{\mbox{\rm GCH}}
\nc{\CH}{\mbox{\rm CH}}
\nc{\zfc}{\mbox{\rm ZFC}}
\nc{\sch}{\mbox{\rm SCH}} 
\nc{\ZF}{\mbox{\rm ZF}}
\nc{\NCF}{\mbox{\rm NCF}} 
\nc{\FD}{\mbox{\rm FD}}   

\nc{\fourG}{\mbox{\rm 4G}} 
\nc{\fourI}{\mbox{\rm 4I}}   

\nc{\Borelhood}{Borel measurability} 
\nc{\Pieinseins}{\mbox{${\bf \Pi}^1_1$}}
\nc{\seinseins}{\mbox{${\bf\Sigma}^1_1$}}
\nc{\seinszwei}{\mbox{${\bf\Sigma}^1_2$}}
\nc{\seinsdrei}{\mbox{${\bf\Sigma}^1_3$}}
\nc{\Deleinszwei}{\mbox{${\bf\Delta}^1_2$}}

\nc{\up}{\ensuremath{{\cal U}\mbox{\ensuremath{\rm -prod}}\,\omega}}
\nc{\upp}{\ensuremath{{\cal U}'\mbox{\ensuremath{\rm -prod}}\,\omega}}
\nc{\upro}{\ensuremath{{\cal U}\mbox{\ensuremath{\rm -prod}}\,\om}}
\nc{\fupro}{\ensuremath{f({\cal U})\mbox{\ensuremath{\rm -prod}}\,\om}}
\nc{\vpro}{\ensuremath{{\cal V}\mbox{\ensuremath{\rm -prod}}\,\om}}
\nc{\fpro}{\ensuremath{{\cal F}\mbox{\ensuremath{\rm -prod}}\,\om}}

\nc{\cff}[1]{{\text{cf}\,(#1)}}           
\nc{\cu}{\ensuremath{\cal U}}             
\nc{\ai}{\ensuremath{\forall^\infty}}     
\nc{\ei}{\ensuremath{\exists^\infty}}     
\nc{\ww}{\ensuremath{\omega^\omega}}      

\nc{\gw}{groupwise dense}

\nc{\kk}{car\-dinal cha\-rac\-teris\-tic}
\nc{\joker}{\ast}
\nc{\gtc}{Galois-Tukey connection} 

\nc{\av}[1]{{\rm Av}_{#1}}
\nc{\eps}{\varepsilon}
\nc{\n}{{\bf n}}                 
\nc{\m}{{\bf m}}

\nc{\marginparr}[1]{}
\nc{\footnoteee}{} 
\nc{\footnotee}{}  

\newcommand{\cal}{\mathcal}

\nc{\divs}{{c_0 \setminus \ell^1}}
\nc{\divser}{(\divs, \leq^*)/\thickapproy}
\nc{\bfin}{\RO(\pfin \setminus\{0\},\subseteq^*)}
\nc{\bdivser}{\RO(\divser)}
\nc{\inc}{{\rm INC}}
\nc{\com}{{\rm COM}}
\nc{\thickapproy}{\makebox{}\!\!\thickapprox}
\nc{\approy}{\makebox{}\!\!\approx}
\nc{\lessi}{\leqslant}
\nc{\gessi}{\geqslant}
\nc{\interior}[1]{{\rm int}(#1)}
\nc{\closure}[1]{{\rm cl}(#1)}
\nc{\Vo}{Vojt\'a\v{s}}

\nc{\precedeseq}{\leq^*} 
\nc{\precedes}{\prec}
\nc{\stronger}{\leqslant_{\bf P}}
\nc{\underlline}[1]{\hat{#1}}
\nc{\PO}{{\bf P}}
\nc{\charak}{\text{ch}}
\nc{\symom}{{\rm{Sym}(\omega)}}

\nc{\needed}{needed\ }
\nc{\neededc}{needed}
\nc{\Needed}{Needed\ }
\nc{\wneeded}{weakly needed\ }
\nc{\Wneeded}{Weakly needed\ }
\nc{\wneededc}{weakly needed}

\nc{\mup}{m_{\rm up}}
\nc{\mdn}{m_{\rm dn}}
\nco{\may}{may}
\nco{\aver}{av} 
\nco{\norm}{nor} 
\nco{\val}{val} 
\nco{\dis}{dis} 
\nco{\basis}{basis}
\nco{\pos}{pos}
\nco{\spec}{spec}
\nc{\err}{\mbox{err}}
\nc{\eee}{\mbox{e}}
\nco{\Expect}{Exp}
\nco{\rt}{rt}
\nco{\pr}{pr}
\nco{\suc}{suc}
\nco{\splitt}{sp}
\nco{\halv}{h}

\nc{\bbforcing}{\mathbb A}
\nc{\itername}{\mathfrak q}
\nc{\iterp}{\mathfrak p}
\nc{\iterq}{\mathfrak q}
\nc{\invcm}{\rm inv_{cm}}
\nc{\invcf}{\rm inv_{cf}}
\nc{\invgm}{\rm inv_{gm}}

\begin{document}


\title{Long low iterations}

\author{Heike Mildenberger and Saharon Shelah}

\thanks{
The first author was supported by the Austrian 
``Fonds zur wissenschaftlichen F\"orderung'', grants no.\ 13983 and 16334
and by travel support by the Landau Center.\\
The second author's research is supported by the United States-Israel Binational Science
Foundation (Grant no.\ 2002323). This is his Publication 843.}

\address{Heike Mildenberger, Universit\"at Wien,
Institut f\"ur formale Logik, W\"ahringer Str.\ 25,
A-1090 Vienna, Austria
}
\address{Saharon Shelah, Institute of Mathematics, The Hebrew University of Jerusalem,
 Givat Ram, 91904 Jerusalem, Israel, and 
Mathematics Department, Rutgers University, New Brunswick,
NJ, USA}

\email{heike@logic.univie.ac.at}
\email{shelah@math.huji.ac.il}

\begin{abstract}
We try to control many cardinal characteristics by working with a notion of
orthogonality between two families of forcings.
We show that $\gb^+ < \gro$ is consistent.

\medskip

{\em This work is dedicated to James Baumgartner on the occasion of his 60th birthday.}

\end{abstract}

\date{March 30, 2004}
\maketitle

\setcounter{section}{-1}
\nothing{\section*{Annotated content}
0. Introduction
1. Pseudo creature forcing
2. Forcing bigness notions
3. Long low c.c.c.\ iterations with ${\bf q} \perp \Gamma$
4. Examples
5. Towards many cardinal characteristics
6. On the consistency of $\gb^+ < \gro$}

\tableofcontents
\section{Introduction}

In this section we define cardinal characteristics that
are the candidates for our arrangements. 
Our aim is to arrange that one  characteristic be small, say be $\kappa$,  and 
to increase the other
characteristic to $\lambda > \kappa^+$. 
The technical means is some
strategically $(<\lambda)$-complete pattern of  c.c.c.\ forcings. 
If the  relations underlying the two characteristics
 are sufficiently orthogonal, this will be possible.

\smallskip

For more than two characteristics at the same time, only very
little is known.

\smallskip

The kinds of invariants that will be forced to be big
are connected to notions of forcings $({\bf Q}, \leq_{\bf q})$
that add special reals
and thus increase some cardinal characteristic.
Which characteristics? We give a pragmatic definition:

\begin{definition}\label{0.1}
Let ${\bf q}=({\bf Q},\leq_{\bf q})$ be a notion of forcing, ${\bf Q} \subseteq {}^\omega \omega$.
\begin{myrules}
\item[(1)]
$\invcm({\bf q}) = \min\{\kappa \such $ there is a $\leq_{\bf q}$-increasing
sequence of length $\kappa$ 
with no $\leq_{\bf q}$ upper bound.$\}$.
\item[(2)]
$\invcm({\bf q},R) = \min\{\kappa \such $ there is a $\leq_{\bf q}$-increasing
sequence $\la \bar{\gc}_\zeta \such \zeta < \kappa\r$ such that for all
 $\eta \in {}^\omega \omega$ there is some
$\zeta < \kappa$ such that $\neg \bar{\gc}_\zeta R\eta\}$.
\item[(3)]
$\invcf({\bf q}) = \min\{|{\mathcal C}|
\such {\mathcal C}\subseteq {\bf Q} \wedge (\forall \bar{\gd} \in {\bf Q})
(\exists \bar{\gc} \in {\mathcal C}) (\bar{\gd} \leq_{\bf q} \bar{\gc})\}$.
\item[(4)]
$\invgm({\bf q}) $ is the minimal $\kappa$ such that in the following game
 $\Game_\kappa({\bf q})$ the empty player
has a winning strategy: 
A play lasts $\kappa +1$ moves and in the $\alpha$-th move the 
non-empty player chooses $\bar{\gc}_\alpha \in {\bf Q}$
satisfying $\beta < \alpha \Rightarrow \bar{\gd}_\beta \leq_{\bf q} \bar{\gc}_\alpha$ and then
the empty player chooses $\bar{\gd}_\alpha \in {\bf Q}$ such 
that $\bar{\gc}_\alpha \leq_{\bf q} \bar{\gd}_\alpha$.
Then the non-empty player wins iff he always has a legal move.
\end{myrules}
\end{definition}

\begin{definition}\label{0.2}
For ${\bf q}$ as in Definition~\ref{0.1} 
we define ${\rm Spec}(\bf q)$ as the set of regular (uncountable $\leq 2^\omega$)
cardinals $\kappa$ such that in the game $\Game^*_\kappa(\bf q)$ the empty player 
has no winning strategy. The game  $\Game^*_\kappa(\bf q)$ is defined just like 
$\Game_\kappa(\bf q)$ except that in stage $\kappa$ non-empty wins if
there were $\bar{\gc}_\alpha$ for all $\alpha <\kappa$ but  there is no 
$\bar{\gc}_\kappa$.
\end{definition}

Now we introduce the orthogonal relations, that shall have
 small characteristics: Some relations $R$, 
different from the ones in the previous definition,  shall have a
small $R$-unbounded set $\{\name{\eta}_i \such i < \omega_1\}$
that needs to be preserved by the forcing ${\mathfrak K} *
\name{{\mathbb P}}$, though we build 
a $\mathfrak K$-generic forcing $\name{{\mathbb P}}$
increasing some of the ``creature''-invariants above. We work with forcing 
bigness notions $\Gamma$ such that $R$ fits $\Gamma$ (see Definition~\ref{0.3}).
The $\name{{\mathbb P}}$ will be described by better and better approximations
$({\mathbb P}, {\mathbb P}^+, \name{\eta})\in \Gamma$ on which the same $\name{\eta}$
keeps its r\^ole as a member of an $R$-unbounded family.
Preliminarily the reader may think of a forcing bigness notion $\Gamma$ as a way of
finding suitable extensions in the $\lessdot$ order of the forcings.
The technical definition of $\Gamma$ will be given in  Definition~\ref{2.1}.

\begin{definition}\label{0.3}
\begin{myrules}
\item[(1)] We say a binary Borel relation $R$ on ${}^\omega
({\mathcal H}(\omega))
$ fits a bigness notion $\Gamma$ if
$({\mathbb P}, {\mathbb P}^+, \name{\eta}) \in \Gamma$ implies that 
$$\Vdash_{{\mathbb P}^+}
\forall \nu \in (^\omega \omega)^{{\bf V}[\mathbb P]} (\neg \name{\eta} R \nu).$$

\item[(2)] Let ${\rm inv}(R) = \min\{ |Y| \such Y \subseteq {}^\omega \omega \wedge
(\forall \nu \in {}^\omega \omega)(\exists \eta \in Y)(\neg \eta R\nu)\}$.

\item[(3)] In $\bf V$ or in any generic extension ${\bf V}[\mathbb P]$ of ${\bf V}$, 
we say that $\eta \in {}^\omega \omega$ is
$(\Gamma,R)$-big over $A\subseteq {}^\omega \omega$ if $R$ is a binary Borel relation on 
${}^\omega({\mathcal H}(\omega))$ which fits $\Gamma$ and $\nu \in A$ implies
 $\neg \eta R \nu$.

\item[(4)] Let $\Gamma$ be a forcing bigness notion in ${\bf V}_1$, ${\bf V}_2$ extend ${\bf V}_1$, 
$A \in {\bf V}_2$, $A \subseteq {\bf V}_1$, $\eta \in (^\omega \omega)^{{\bf V}_1}$.
We say that $\eta$ is $\Gamma$-big over $A$ in $({\bf V}_1,{\bf V}_2)$ if for some 
$({\mathbb P}, {\mathbb P}^+,\name{\eta}) \in \Gamma^{{\bf V}_1}$ and 
$G \in {\bf V}_2$, $G$ is ${\mathbb P}^+$-generic over ${\bf V}_1$, 
$A \in {\bf V}_1[G \cap {\mathbb P}]$, and
$\eta = \name{\eta}[G]$.
\end{myrules}
\end{definition}

\section{Pseudo creature forcing}\label{S1}

Now we introduce some forcings ${\bf q}= ({\bf Q}, \leq_{\bf q})$, whose invariants
according to Definition~\ref{0.1} shall be 
increased by some compound forcings later
in Section~3.

\smallskip

We think of forcings with linear creatures as in Shelah \cite{Sh:207}, Blass and Shelah
 \cite{BlassShelah},
Proper and Improper Forcing  \cite[VI,\S 6]{Sh:f} and
Ros\l anowski and  Shelah \cite{RoSh:470}. \nothing{In this section, though, 
we do not require that the norms
diverge, because this requirement does not play a r\^ole here.}

\begin{definition}\label{1.1}
We say ${\bf q} = (K_{\bf q}, \Sigma_{\bf q}, \norm_{\bf q},\val_{\bf q}) = (K,\Sigma,\norm,\val)$
is a {\em pseudo creature forcing} (abbreviated by pcrf) iff
\begin{myrules}

\item[(a)] $K \subseteq {\mathcal H}(\aleph_0)$ is a set of creatures of the form
$\gc \colon [\mdn(\gc),\mup(\gc)) \to {\mathcal H}(\aleph_0)$,
and for each $\gc \in K$ we have $\mdn(\gc)  < \mup(\gc)  < \omega$.

\item[(b)] The subcomposition function $\Sigma$ is a partial function from ${}^{\omega >} K$ 
to ${\mathcal P}(K) \setminus\{\emptyset\}$ such that
\begin{myrules}

\item[($\alpha$)] $\dom(\Sigma) \subseteq \{ (\gc_n,\dots,\gc_{m-1}) \such
n<m<\omega \wedge
 \forall \ell \in [n,m)(\gc_\ell \in K \wedge \mup(\gc_\ell) \leq \mdn(\gc_{\ell+1}))\}$.

\item[($\beta$)]
$\Sigma(\gc_n,\dots,\gc_{m-1}) \subseteq \{ \gc \in K \such 
(\exists k,\ell)(n\leq k <\ell \leq m
\wedge \mdn(\gc) = \mdn(\gc_k), \mup(\gc) = \mup(\gc_{\ell-1}))\}$.

\item[($\gamma$)] $\Sigma$ is transitive: If $\gc_\ell^1 \in \Sigma(\gc_{m_\ell}^0,\dots,
\gc_{m_{\ell+1}-1}^0)$ for $\ell <k$ and $\gc^2 \in \Sigma(\gc_0^1,\dots,\gc_{k-1}^1)$, then
$\gc^2 \in \Sigma(\gc_{m_0}^0,\dots,\gc_{{m_1}-1}^0,\gc^1_{m_1},\dots ,\gc^1_{{m_2}-1},\dots
, \gc^{k-1}_{m_k -1})$.
That also means, that the domain of $\Sigma$ is closed under the mentioned
concatenations. 

\item[($\delta$)] $(\forall \gc \in K) (\gc \in \Sigma(\gc))$
and $\Sigma $ is monotonous in the following sense:
If $i_0 < \cdots < i_{m-1} <n$ then
$$\Sigma(\gc_{i_0}, \dots, \gc_{i_{n-1}}) \subseteq \Sigma(\gc_{0}, \dots, \gc_{n-1}).$$

\end{myrules}

\item[(c)] $\norm \colon K \to {\mathbb R}^{\geq 0 }$ (usually $\omega$).

\item[(d)] $\val \colon K \to {\omega}^{<\aleph_0}$,
$\emptyset \neq \val({\gc}) \subseteq [\mdn(\gc),\mup(\gc))$.

\end{myrules}
\end{definition}

{\bf Remark}: Note that Definition~1.1. (b)($\beta$) is a stronger
requirement on the subcomposition functions than the usual
$\mdn(\gc) \geq \mdn(\gc_n)$, $\mup(\gc) \leq \mup(\gc_{m-1})$.

\begin{definition}\label{1.2}
1. The set of pure conditions of ${\bf q}$ is
\begin{equation}
\begin{split}
{\bf Q}^{pr}_{\bf q}= \{ \bar{\gc}\such \bar{\gc} = &\langle 
\gc_n \such n<\omega\rangle, \gc_n \in K, \mup(\gc_n) = \mdn(\gc_{n+1}),
\mdn(\gc_0) = 0\\
& \wedge \lim_{n \to \omega}\norm(\gc_n) = \infty\}.
\end{split}
\end{equation}

We usually write $\bar{\gc} \in K$.
 Let ${\bf Q}^{p}_{\bf q}$ be defined similarly, allowing
$\mdn(\gc_0) >0$.

\medskip

2. The set of not necessarily pure conditions ${\bf Q}^{np}_{\bf q}$ of ${\bf q}$ is
\begin{equation}
\begin{split}
{\bf Q}^{np}_{\bf q} = \{ \bar{\gc}\such \bar{\gc} = & \langle 
\gc_n \such n<\omega\rangle, \gc_n \in K, \mup(\gc_n) \leq \mdn(\gc_{n+1}),\\
& \wedge \lim_{n \to \omega}\norm(\gc_n) = \infty\}.
\end{split}
\end{equation}

\medskip

3. The true creature forcing ${\bf Q}^{tr}_{\bf q}$ is a subset of 
\begin{equation*}
\begin{split}
\{\langle w,\gc_0,\gc_1,\dots \rangle \such & (\exists n_0 <n_1 <\dots) 
(\forall i < \omega)(\exists u_i \subseteq [n_i,n_{i+1}))(\exists \bar{\gc}^* \in 
{\bf Q}^{np}_{\bf q})\\
&(u_i \neq\emptyset \wedge
w \in {\rm pos}(\langle\rangle,\gc_0^*,\dots \gc_{n_0-1}^*) \wedge
\gc_i \in \Sigma(\gc_m^*)_{m\in u_i})\}, 
\end{split}
\end{equation*}

This is not quite the general case of creature forcing. In the terminology of
\cite[Definition 2.1.1.]{RoSh:470} all our forcings are omittory.
The pos (possibility) operation is defined as follows:
\begin{equation}
\begin{split}
{\rm pos}(w,\gc_0,\dots \gc_{n-1}) = &
\{w \cup  \bigcup_{\ell <k} w_\ell \such \exists \la u_\ell, \gd_\ell \such \ell < k \r,
u_\ell \neq \emptyset, u_\ell \subseteq n,\\
&
 \max(u_\ell)
<\min(u_{\ell +1}), \gd_\ell \in \Sigma((\gc_i)_{i \in u_\ell }),
w_\ell \in \val(\gd_\ell)\}.
\end{split}
\end{equation}

\end{definition}

\begin{remark}\label{1.3}
${\bf Q}^{pr}_{\bf q}$ and ${\bf Q}^{p}_{\bf q}$ are the pure elements in \cite{RoSh:470}.
Pure means: $\mup(\bar{\gc}_n) = \mdn(\bar{\gc}_{n+1})$ for $n<\omega$, so the
union of the domains is $\omega$ in the case of ${\bf Q}^{pr}_{\bf q}$ and
$\omega \setminus n$ for some $n$  in the case of ${\bf Q}^{p}_{\bf q}$.
Of course, one could also introduce a pure true creature forcing.
But for our intended applications, pureness is often  too 
strong a requirement.
\end{remark}

Now we equip ${\bf Q}^{pr}_{\bf q}$ and ${\bf Q}^{p}_{\bf q}$ 
with two partial orderings $\leq_{\bf q}$ 
and the stricter  $\leq_{\rm full}$. 
\nothing{For $({\bf Q_q},\leq_{\bf q})$ we shall just write
${\bf Q_q}$ if there is no ambiguity, because this is the most frequently
used partial order.}

\begin{definition}\label{1.4}
1.
For $\bar{\gc}^1,\bar{\gc}^2 \in {\bf Q}^{pr}_{\bf q}$ or in ${\bf Q}^{p}_{\bf q}$ we write
$\bar{\gc}^1 \leq_{\bf q}^{pr} \bar{\gc}^2 $ iff there is some $i(*) \in \omega$ and there is
some sequence $n_{i(*)} < n_{i(*)+1} < \dots$ such that for every
$i \geq i(*)$ $\gc^2_i \in \Sigma(\gc_{n_i}^1,\dots,\gc_{n_{i+1}-1}^1)$
and $\mdn(\gc^2_i) = \mdn(\gc^1_{n_i})$ and
$\mup(\gc^2_i) = \mup(\gc^1_{n_{i+1}-1})$
(so, and in Definition~\ref{1.1}(b)($\beta$) there is equality in the pure case).

We say that $(i(*),n_{i(*)})$ witnesses $\bar{\gc}^1 \leq_{\bf q}^{pr} \bar{\gc}^2$.

\smallskip

If we have in addition $i(*)=0$ (then, in case $\bar{\gc}^1,\bar{\gc}^2 
\in {\bf Q}^{\pr}_{\bf q}$, also
$n_{i(*)} =0$),
then we write  $\bar{\gc}^1 \leq_{\rm full}^{pr} \bar{\gc}^2 $.
The reversed relations are $\geq_{\bf q}^{pr}$ and $\geq_{\rm full}^{pr}$ respectively.

\medskip

2.
For $\bar{\gc}^1,\bar{\gc}^2 \in {\bf Q}^{np}_{\bf q}$ we write
$\bar{\gc}^1 \leq^{np}_{\bf q} \bar{\gc}^2$ iff there is some $i(*) \in \omega$ and there is
some sequence $n_{i(*)} < n_{i(*)+1} < \dots$ and some
sequence $\emptyset \neq u_i \subseteq [n_i,n_{i+1})$ such that for every
$i \geq i(*)$, $\gc^2_i \in \Sigma((\gc_{z}^1)_{z \in u_i})$.
So 
$\mdn(\gc_i^2) = \mdn(\gc_{\min(u_i)}^1)$, $\mup(\gc_i^2) = \mdn(\gc_{\max(u_i)}^1)$.

\medskip

3. For $(w^1,\bar{\gc}^1),(w^2,\bar{\gc}^2) \in {\bf Q}^{tr}_{\bf q}$ we write
$(w^1,\bar{\gc}^1) \leq^{tr}_{\bf q} (w^2,\bar{\gc}^2) $ iff there is some 
some sequence $0\leq n_{0} < n_{1} < \dots$ and some
sequence $u_i \subseteq [n_i,n_{i+1})$ such that for every
$i \in \omega$, $\gc^2_i \in \Sigma((\gc_{z}^1)_{z \in u_i})$,
and such that $w^1 \trianglelefteq w^2$ and 
$w^2 \setminus w^1 \in \pos(\bar{\gc}^1 \restriction n_0)$.

\end{definition}

\nothing{
\begin{definition}\label{1.b}
For $\bar{\gc}^* \in{\bf Q_q}$ we define the restricted partial order
${\bf Q_q}(\geq^* \bar{\gc}^*)= \{ \bar{\gc} \in {\bf Q_q} \such \bar{\gc} \geq^* \bar{\gc}^*\}$.
\end{definition}
{\sf This definition is not yet used.}}

\begin{convention}\label{1.5}
1. In the pure context, we write $\langle \gc_{n_0},\dots, 
\gc_{n_1-1}\rangle \in \Sigma(\gc_{k_0}, \dots,\gc_{k_1-1})$
if there are
$k_0= m_0 < m_1 < \cdots <  m_{n_1} =k_1$ such that for all  $j \in [n_0,n_1)$, 
$\gc_j \in \Sigma(\gc_{m_j}, \dots,\gc_{m_{j+1}-1})$.

2. In the not necessarily pure context, we write 
$\langle \gc_{n_0},\dots, \gc_{n_1-1}\rangle \in \Sigma(\gc_{k_0}, \dots,\gc_{k_1-1})$
if there are
$k_0= m_0 < m_1 < \cdots <  m_{n_1} =k_1$ and $\emptyset \neq u_i \subseteq [m_i,n_{i+1})$
such that for all  $j \in [n_0,n_1)$, 
$\gc_j \in \Sigma((\gc_{z})_{z \in u_i})$.
\end{convention}

The following property can be used for all the variants of notions of forcing
$({\bf Q}^{xy}_{\bf q},\leq^{xy}_{{\bf q},{\rm full}})$ so far defined.

\begin{definition}\label{1.6}
We say that a pcrf ${\bf q}$ is forgetful if the following holds:
If $\bar{\gc}^1 \leq_{\bf q} \bar{\gc}^2$ and $n\leq k<\omega$ and 
$\bar{\gc}^1\restriction k = \bar{\gc}^2\restriction k$ 
and $\bar{\gge}^\ell \in \Sigma(\bar{\gc}^1\restriction [n,k))$ for $\ell =1,2$, 
then we can find $m_1,m_2$ such that $\mdn(\gc^1_{m_1}) = \mdn(\gc^2_{m_2})$, $m_1 >k$, $m_2 > n$
and $$
\bigcap_{\ell =1}^{2} \Sigma(\bar{\gge}^\ell \concat \bar{\gc}^1\restriction [k,m_1)) \cap
\Sigma(\bar{\gc}^2 \restriction [n,m_2)) \neq \emptyset.
$$
\end{definition}

Note that 
$ \lg(\bar{\gge}^\ell) \leq k -n$ and may be strictly less,
and that
$$(\bar{\gge}^\ell \concat \bar{\gc}^1)(i)=
\left\{ \begin{array}{ll}
\gge^\ell_i, &\mbox{ if } i < \lg(\bar{\gge}^\ell),\\
\gc^1_{k+j}, & \mbox{ if } i = \lg(\bar{\gge}^\ell) +j.
\end{array}
\right.
$$
but since we do not need to look so
close at the conditions this shift in indexing 
will not appear explicitly. 
 
\begin{convention}\label{1.7}
We let ${\bf q}$ denote a forgetful pcrf unless said otherwise.
\end{convention}

\begin{definition}\label{1.8}
We define an equivalence relation $\equiv_{\bf q}$ on ${\bf Q}^{np}_{\bf q}$
and on ${\bf Q}^{tr}_{\bf q}$:
$\bar{\gc}^1 \equiv_{\bf_q} \bar{\gc}^2$ iff there are $k_1, k_2 \in \omega$ 
such that for all $n<\omega$, $\gc^1_{k_1+n} = \gc^2_{k_2 +n}$.
\end{definition}

All equivalence classes of this relation are countable. In the language of
Kechris, Hjorth, Louveau and others this is a countable equivalence relation.

\begin{definition}\label{1.9}
Let $\bf q =({\bf Q}_{\bf q}, \leq_{\bf q})$ any of our versions of forcings.
\begin{myrules}
\item[1.]
${\mathcal C}$ is a ${\bf q}$-pre-directed family if 
\begin{myrules}
\item[(a)]
${\mathcal C} \subseteq {\bf Q}_{\bf q}$,
\item[(b)] for some proper filter $D$ on
$\omega$ containing the cofinite sets of $\omega$, any finitely many members of 
$\mathcal C$ have a common upper bound $\bar{\gc}$ w.r.t.\ $\leq_{\bf q}$
with $\supp(\bar{\gc}) = \{\mdn(\gc_n) \such n \in \omega \} \in D$.
\end{myrules}

Note that in (b) we did not require that $\gc \in {\mathcal C}$. This explains the
prefix ``pre''. Such a filter $D$ is called a witness for ${\mathcal C}$ being 
$\bf q$-pre-directed.

\smallskip

\item[2.]
 ${\mathcal C}$ is a $\bf q$-directed family if 
\begin{myrules}
\item[(a)] ${\mathcal C} \subseteq {\bf Q_q}$,

\item[(b)] $({\mathcal C}, \leq_{\bf q})$ is directed, i.e., every two members 
of $\mathcal C$ have a common upper bound in $\mathcal C$ w.r.t.\ $\leq_{\bf q}$,

\item[(c)] If $\bar{\gc}^1 \equiv_{\bf q} \bar{\gc}^2 \in {\bf Q_q}$, then
$\bar{\gc}^1 \in {\mathcal C}$ iff $\bar{\gc}^2 \in {\mathcal C}$.
\end{myrules}
\end{myrules}
\end{definition}

In the following lemma we use the strong requirement from Definition~\ref{1.1}(b)($\beta$).

\begin{lemma}\label{1.10}
Let the order be $\leq_{\bf q}^{pr}$, $\leq_{\bf q}^{np}$, 
or $\leq_{\bf q}^{tr}$.
Every $\bf q$-directed family is $\bf q$-pre-directed.
\end{lemma}

\proof Let ${\mathcal C}$ be $\bf q$-directed. 
Take $D$ be the filter generated by
$\{\supp(\bar{\gc}) \such \bar{\gc} \in {\mathcal C}\}$. 
If $\bar{\gc}^1 \leq_{\bf q}^{pr} \bar{\gc}^2$,
then $\supp(\bar{\gc}^2) \subseteq^* \supp(\bar{\gc}^1)$, and for $\leq_{\rm full}^{pr}$
there are no exceptions to the inclusion.
\proofend

\begin{definition}\label{1.11}
Note that $\bbforcing$ is not identical with $\bf Q$.
\nothing{Saharon, I changed $\mathbb Q$ into $\mathbb A$, 
following Baumgartner, Hru\v{s}\'ak and others.
This is likely to  help to prevent bad stories with the typesetters in the end.
I have bad experience with Elsevier and with Springer. The Israel Journal
seems to be more careful.}

For a $\bf q$-directed family $\mathcal C$ in $({\bf Q}_{\bf q}, \leq_{\bf q})$
 we define a partial order 
${\bbforcing}_{\bf q}({\mathcal C})=
{\bbforcing}({\mathcal C})$ as follows:
\begin{myrules}
\item[(a)] the elements are pairs $(n,\bar{\gc})$ where $n<\omega$ and $\bar{\gc} \in {\mathcal C}$,
\item[(b)] the order is $(n_1,\bar{\gc}^1) \leq_{{\bbforcing}({\mathcal C})}
(n_2,\bar{\gc}^2)$ iff
\begin{myrules}
\item[($\alpha$)] $n_1\leq n_2$,
\item[($\beta$)] $\bar{\gc}^1 \restriction n_1 = \bar{\gc}^2 \restriction n_1$,
\item[($\gamma$)] $\bar{\gc}^1 \restriction [n_1,\omega) \leq_{\rm full} \bar{\gc}^2 \restriction 
[n_1,\omega)$, which means, by the definition of $\leq_{\rm full}$: for some
$m_0=n_1 < m_1 <  m_2 \dots$ we have for all $i<\omega$, $\gc^2_{n_1 +i} \in 
\Sigma(\gc^1_{m_i},\dots,\gc_{m_{i+1}-1}^1)$.
\end{myrules}
\end{myrules}
Often we write only $\leq$ for $\leq_{{\bbforcing}({\mathcal C})}$.

\end{definition}

\nothing
{\bf Saharon, there is a change:
We do not really want to work with ${\bbforcing}(\mathcal C)$ but
with
 $({\bf Q}_{\bf q}^{tr}(\mathcal C),\leq_{\bf q}^{tr})$ 
which looks fortunately quite similar to the versatile forcing from
Definition~\ref{1.10}, only that in item (b)($\gamma$)  we
want to go to subsets $u_i \subseteq [m_i,m_{i+1})$ and
we do not want $\la m_i\such i<\omega\r$ to be readable from the generic
real, because that will add a dominating real all the time.}
\nothing{
 what about the
quantification over the $m_i$ in the non-pure case? 
Is $\la m_i \such i <\omega \r$ definable form $\bar{\gc}_1$ and $\bar{\gc}_2$ alone?
We need some absoluteness in Claims~\ref{1.16} and Conclusion~\ref{1.17}.
So I see difficulties in not increasing $\gb$ by the 
forgetful creature forcing because of the $\la m_i\such i<\omega \r$}

\begin{definition}\label{1.12}
Let $\mathcal C$ be $\bf q$-directed.
For $\gc \in {\mathcal C}$, let ${\bf Q}(\bar{\gc}) = \{ \bar{\gc}' \in {\bf Q} 
\such \bar{\gc}' \equiv_{\bf q} \bar{\gc} \}$
and let ${\bbforcing}(\bar{\gc}) = \{ (n,\bar{\gc}') \such n < \omega, \bar{\gc}' \in 
{\bf Q}(\bar{\gc})\}$.
Let ${\bf Q}_{\mathcal C}(\bar{\gc}) = \{ \bar{\gc}' \in {\mathcal C} 
\such \bar{\gc}' \equiv_{\bf q} \bar{\gc} \}$
and let ${\bbforcing}_{\mathcal C}(\bar{\gc}) = \{ (n,\bar{\gc}') \such n < \omega, \bar{\gc}' \in 
{\bf Q}_{\mathcal C}(\bar{\gc})\}$.
Each $\equiv_{\bf q}$-class is infinite and countable. Hence, if above each condition
 there are two incompatible conditions, then
the  second and the fourth forcing (with $\leq_{{\bbforcing}(\mathcal C)}$) are equivalent
to the Cohen forcing.
\end{definition}

\begin{definition}\label{1.13}
We say that $\mathcal C$ is strongly ${\bf q}$-directed if
$({\mathcal C},\leq_{\bf q})$ is ${\bf q}$-directed and  $\aleph_1$-directed, that is, 
every countable set has a common upper bound.
\end{definition}

\begin{claim}\label{1.14}
Let $\bf q$ be a forgetful pcrf.

\smallskip

1. If ${\mathcal C}$ is $\bf q$-directed, then ${\bbforcing}({\mathcal C})$ is a
$\sigma$-centered forcing notion.

\smallskip

2. If ${\mathcal C}$ is $\bf q$-pre-directed and $\bar{\gc}^n \in {\mathcal C}$ for 
$n<\omega$ (possibly with repetitions), then
\begin{equation*}
\begin{split}
\Vdash_{\rm Cohen} & \mbox{``for some } \bar{\gc} \in {\bf Q_q}^{{\bf V}[{\rm Cohen}]} 
\mbox{ we have }
 n<\omega \rightarrow \bar{\gc}^n \leq_{\bf q} \bar{\gc}\\
& \mbox{and } {\mathcal C} \cup \{ \bar{\gc}\} \mbox{ is {\bf q}-pre-directed.''}
\end{split}
\end{equation*}

\smallskip

3. Let $\mathcal C$ be a $\bf q$-pre-directed family and ${\mathbb P}$ be the forcing of adding 
$2^{\aleph_0}$ Cohen reals. Then in ${\bf V}^{\mathbb P}$ we can extend $\mathcal C$ to a strongly
$\bf q$-directed ${\mathcal C}'$.

\smallskip

4. If ${\mathcal C}$ is strongly $\bf q$-directed and ${\mathcal I}$ is a predense subset
of ${\bbforcing}({\mathcal C})$, then for some $\bar{\gc} \in
{\mathcal C}$ and some $\mathcal J$ we have
\begin{myrules}
\item[($\alpha$)] ${\mathcal J} \subseteq {\bbforcing}_{\mathcal C}(\bar{\gc})$, hence is countable,
\item[($\beta$)] $\mathcal J$ is predense in ${\bbforcing}({\mathcal C})$,
\item[($\gamma$)] $(\forall q \in {\mathcal J})(\exists p \in {\mathcal I}) 
(p \leq_{{\bbforcing}({\mathcal C})} q)$.
\end{myrules}

\smallskip

5. If ${\mathcal C}^1\subseteq {\mathcal C}^2$ are $\bf q$-directed, 
\nothing{(closure under $\equiv_{\bf q}$ is enough)} then
\begin{myrules}

\item[($\alpha$)] If $p,q \in {\bbforcing}({\mathcal C}^1)$ and
 ${\bbforcing}({\mathcal C}^1) \models p \leq q$ then
${\bbforcing}({\mathcal C}^2) \models p \leq q$,

\item[($\beta$)]  If $p,q \in {\bbforcing}({\mathcal C}^1)$ and
 ${\bbforcing}({\mathcal C}^1) \models p \not\leq q$ then
${\bbforcing}({\mathcal C}^2) \models p \not\leq q$,

\item[($\gamma$)]  If $p=(n^p,\bar{\gc}^p)$, $q=(n^q,\bar{\gc}^q) \in {\bbforcing}({\mathcal C}^1)$ 
are incompatible in  ${\bbforcing}({\mathcal C}^1)$ and $\bar{\gc}^p \leq_{\bf q} \bar{\gc}^q$ then
$p,q$ are incompatible in
${\bbforcing}({\mathcal C}^2)$.

\item[($\delta$)]   If $p=(n^p,\bar{\gc}^p)$,
$q=(n^q,\bar{\gc}^q) \in {\bbforcing}({\mathcal C}^1)$ 
are incompatible in  ${\bbforcing}({\mathcal C}^1)$, then
$p,q$ are incompatible in
${\bbforcing}({\mathcal C}^2)$.
\end{myrules}

\smallskip

6. The family of ${\bf q}$-pre-directed families ${\mathcal C}$
as well as the family of ${\bf q}$-directed families ${\mathcal C}$ are
 closed under increasing unions.
\end{claim}

\proof 
 1. Choose $\bar{\gc}^*\in {\mathcal C}$ and fix it for the rest of the
proof. (Since ${\mathcal C}$ is $\bf q$-directed, the choice does not matter.)

We set 
\begin{equation*}
Y =  \{(\bar{\gd},n_0,n_1,n_2) \such  
n_0 \leq n_1 <\omega, n_2 \in \omega, \bar{\gd} \in {}^{n_1} K,
 \nothing{\bar{\gc}^*\in {\mathcal C}} \}.
\end{equation*}

For each $y =(\bar{\gd},n_0,n_1,n_2)= (\bar{\gd}^y,n_0^y,n_1^y,n_2^y) \in Y$ we define
\begin{equation}\label{Ay0}
\begin{split}
{\bbforcing}_y^0 = \{ p \in {\mathbb Q}({\mathcal C})\such &
n^p = n_0, \bar{\gc}^p \restriction n_1 = \bar{\gd},
\bar{\gc}^* \restriction [n_2,\omega) \leq_{\rm full} \bar{\gc}^p\restriction[n_1,\omega),\\
&\mbox{so, in particular } \mdn(\gc_{n_2}^*) = \mdn(\gc^p_{n_1}) \}
\end{split}
\end{equation}
and 
\begin{equation}\label{Ay}
{\bbforcing}_y = \{ p \in {\bbforcing}({\mathcal C}) \such (\exists q \in {\bbforcing}_y^0)
(p \leq q)\}. 
 \end{equation}

Now, since $\bar{\gd} \in {}^{\omega >} K$ and $K \subseteq {\mathcal H}(\omega)$ is
countable, we have that $Y$ is countable.

\smallskip

Now we prove that ${\bbforcing}({\mathcal C}) = \bigcup \{ {\bbforcing}_y \such y \in Y\}$.
Let $p \in {\bbforcing}({\mathcal C})$. Recalling that ${\mathcal C}$ is $\leq_{\bf q}$-directed,
there is some $\bar{\gc}^1 \in {\mathcal C}$ such that $\bar{\gc}^p \leq_{\bf q} \bar{\gc}^1$
and
$\bar{\gc}^* \leq_{\bf q} \bar{\gc}^1$.
By the definition of $\leq_{\bf q}$ and considering, that once $i(*)$ is large enough,
all larger $i$ could serve as well, we find
$k_0, k_1,k_2$ such that $(k_0,k_2) $ witness $\bar{\gc}^p \leq_{\bf q} \bar{\gc}^1$
and $k_0 \leq k_1 = n_{i(*)}$ and $(k_1,k_2)$ witness  $\bar{\gc}^* \leq_{\bf q} \bar{\gc}^1$.
We choose $k_0 \geq n^p$.

Now let $y =(\bar{\gc}^p \restriction k_1, k_0,k_1,k_2)$.
Easily $y \in Y$,
and let $\bar{\gc}^2 = (\bar{\gc}^p \restriction k_1,\gc^1_{k_2}, \gc^1_{k_2+1},\dots )$ 
and $p_2 = (k_0,\bar{\gc}^2)$. Now, since $\gc \in \Sigma(\gc)$ and since
$k_0 \geq n^p$ and $k_0 \geq n_{i(*)} $ for $\bar{\gc}^p \leq_{\bf q} \bar{\gc}^1$, 
we have $\bar{\gc}^p\restriction k_0 = \bar{\gc}^2\restriction k_0$
and $\bar{\gc}^p \restriction [k_0,\infty) \leq 
\bar{\gc}^2 \restriction [k_1,\infty)$, hence
 $p \leq p_2 \in {\bbforcing}({\mathcal C})$, since
$\bar{\gc}^1 \restriction [k_2,\infty) \leq_{\rm full}
\bar{\gc}^2[k_1,\infty)$ we have
$p_2 \in {\bbforcing}_y^0$, hence $p \in {\bbforcing}_y$ as required.

\smallskip

Each ${\bbforcing}_y$ is directed in $\leq = \leq_{{\bbforcing}({\mathcal C})}$:
Let $y = (\bar{\gd}, n_0,n_1,n_2) $ and let $p_1, p_2 \in {\bbforcing}_y$. We can find
$q_1,q_2 \in {\bbforcing}_y^0 $ such that $p_\ell \leq q_\ell$ for $\ell =1,2$.
So $\bar{\gc}^{q_1}, \bar{\gc}^{q_2} \in {\mathcal C}$ and we can find
$\bar{\gc} \in {\mathcal C}$ such that $\bar{\gc}^{q_\ell} \leq_{\bf q} \bar{\gc}$
for $\ell = 1,2$. Hence there are $k_1, k_2,k \in \omega$ such that
$(k_\ell,k)$ witnesses $\bar{\gc}^{q_\ell} \leq_{\bf q} \bar{\gc}$ for $\ell =1,2$ and
$k \geq n^{q_\ell}$, $k_\ell \geq n_2$. 
Since $q_\ell \in {\bbforcing}_y^0$, for some $k^* > n_1$ we have
$$
\bar{\gc}^{q_\ell} \restriction [n_1,k_\ell) \in \Sigma(\bar{\gc}^*
\restriction [n_2,k^*)) \mbox{ for }\ell =1,2.$$
Without loss of generality $\bar{\gc} \restriction n_1 = \bar{\gc}^* \restriction n_1$.
Recall 
$\bar{\gc}^{q_\ell} \restriction n_1 = \bar{\gd}$.
Now apply the definition of forgetfulness with $n_2, n_1, 
\bar{\gc}^{q_\ell} \restriction [n_1,k_\ell)
(\ell =1,2), \bar{\gc}^*, \bar{\gc}$ standing for $n,k,\bar{\gge}_\ell (\ell =1,2),
\bar{\gc}^1, \bar{\gc}^2$ gives $m_1$ and $m_2$ and $\bar{\gc}'$ in the intersection.
Let $\bar{\gc}'' = (\bar{\gd}\concat \bar{\gc}' \concat \bar{\gc} \restriction [m_2,\omega))$.
Then $(n_0,\bar{\gc}'') \in {\bbforcing}_y^0$ is a common upper bound of
$q_1, q_2 \in {\bbforcing}^0_y$ in $\leq$. Since the upper bound is in ${\bbforcing}_y$,
 we can iterate the process
and find common upper bounds in ${\bbforcing}_y$ for finitely many $p_1, \dots,p_n$.

\medskip

2.  We shall work on the following two cases
\begin{myrules}
\item[(A)] $\{\bar{\gc}^n \such n \in \omega \}$ has just two members $\bar{\gc}^0$ 
and $\bar{\gc}^1$.
\item[(B)] For all $n$, $\bar{\gc}^n \leq_{\bf q} \bar{\gc}^{n+1}$.
\end{myrules}

In both cases we define a countable notion of forcing ${\mathbb P}$ that adds a solution.
First we show that to consider (A) and (B) is enough. We 
set $\bar{\gc}^0 = \bar{\gc}^0$ and, for $n>0$ we iteratively
add $\bar{\gd}^n \geq_{\bf q} \bar{\gc}^n, \bar{\gd}^{n-1}$. Thereafter we add $\bar{\gd}^{\omega}$
above all $\bar{\gd}^n$'s. So we add $\omega +1$ times a Cohen real and
perform this iteration with finite
supports, and this is again equivalent to Cohen forcing. Again forgetfulness makes
it work.

In case (A) we choose $m_0^\ell < m_1^\ell <\dots$ such that $\mdn(\gc^0_{m_i^0})
= \mdn(\gc^1_{m_i^1})$ for $i<\omega$ and such that $\{\mdn(\gc^0_{m_i^0}) \such i <\omega \}
\in D$. Now, by forgetfulness, the Cohen forcing, in its disguised form
$(\{\langle\gd_0,\dots \gd_{k-1}\rangle \such (\forall \ell <k) \gd_\ell \in
\Sigma(\bar{\gc}^0\restriction [m^0_\ell,m^0_{\ell+1})) \cap
\Sigma(\bar{\gc}^1\restriction [m^1_\ell,m^1_{\ell+1}))\},\trianglelefteq)$
will add a generic $\bar{\gd}$ such that $\supp(\bar{\gd}) \in D$ and 
such that $\bar{\gd} \geq_{\bf q} \bar{\gc}^0,\bar{\gc}^1$.

In case (B) we force with Cohen forcing in the guise
$(\bigcup_{n<\omega} {\bbforcing}_{\mathcal C}(\bar{\gd}^n), \leq_{{\bbforcing}({\mathcal C})})$.
By density arguments, the generic $\bar{\gd}$ 
will fulfill $\bar{\gd} \geq_{\bf q} \bar{\gd}^n$ for all $n$.
Note that we did not enlarge the filter $D$.

\medskip

3. We repeat part 2.\ $2^{\aleph_0}$ times, by bookkeeping adding common upper bounds
to all  countable subsets
of all intermediate ${\mathcal C}_1 \supseteq {\mathcal C}$ appearing along
the iteration and such that only ${\bf q}$-pre-directed families appear along the iteration.

\medskip

4. Let ${\mathcal J}^+ = \{(n,\bar{\gc})\in {\bbforcing}({\mathcal C}) \such (n,\bar{\gc}) $
is above some member of ${\mathcal I}\}$. This is dense and open, and let ${\mathcal J}'$
be a maximal set of pairwise incompatible elements of ${\mathcal J}^+$.
By part 1.\ of this claim, ${\mathcal J}'$ is countable. Hence
$\{\bar{\gc}^p \such p \in {\mathcal J}' \}$ is countable. Since ${\mathcal C}$ is
strongly $\bf q$-directed, there is a $\leq_{\bf q}$-upper bound
$\bar{\gc}^*$ of  $\{\bar{\gc}^p \such p \in {\mathcal J}' \}$.
Finally we set
$$
{\mathcal J} = \{p' \in {\bbforcing}_{\mathcal C}(\bar{\gc}^*) \such
p' \mbox{ is above some member of } {\mathcal J}'\}.
$$
Clearly ${\mathcal J}$ satisfies requirements ($\alpha$) and ($\gamma$).
As for  ($\beta$), let $r\in {\bbforcing}({\mathcal C})$ so $r$ is compatible 
with some $p \in{\mathcal J}'$. Choose $r^+ \in {\bbforcing}({\mathcal C})$ such that 
$r^+ \geq r,p$. Since $\bar{\gc}^{r^+}$ and $\bar{\gc}^*$ are compatible in 
$({\mathcal C},\leq_{\bf q})$, there is a common upper bound $\bar{\gd} \in {\mathcal C}$.
\relax From the definitions of $\bar{\gc}^p \leq^* \bar{\gc}^{r^+} \leq_{\bf q} \bar{\gd}$ and
of $\bar{\gc}^p \leq_{\bf q} \bar{\gc}^*$ we get series $\langle m_i^\ell \such i < \omega\rangle
$, $\ell = 0,1,2,3$ such that for every $i$,
\begin{equation*}
\begin{split}
&\bar{\gc}^* \restriction[m_i^1,m_{i+1}^1) \in \Sigma(\bar{\gc}^p 
\restriction[m_i^0,m_{i+1}^0)) \mbox{ and}\\
&\bar{\gd} \restriction[m_i^3,m_{i+1}^3) \in \Sigma(\bar{\gc}^* \restriction[m_i^1,m_{i+1}^1))\cap
 \Sigma(\bar{\gc}^{r^+} \restriction[m_i^2,m_{i+1}^2)).
\end{split}
\end{equation*}
W.l.o.g., $m_0^2 \geq n^{r^+}$ and $m_0^2\geq n^p$. Then $q=(m_0^2,
(\bar{\gc}^{r^+}\restriction m_0^2) \concat \bar{\gd} \restriction[m_0^3,\omega))$.
Now $q \geq p'=(n^p,\bar{\gc}^p \restriction n^p \concat \bar{\gc}^* \restriction [m^1_0,\omega))
 \geq  p$, and we have $p' \in {\mathcal J}$.
Moreover $q \geq r^+\geq r$ as required. So ${\mathcal J}$ is predense.

\medskip

5. Part ($\alpha$) follows immediately from the 
definition. For part ($\beta$) take into account
that the $m_i$, $i <\omega$, in Definition~\ref{1.1}(b)($\gamma$) 
are unique, and hence $\la m_i \such i < \omega \r$ is already in the ground model where
${\mathcal C}^1$ exists.

\nothing{\bf In the case of ${\bf Q}^{np}$ we work with absoluteness of well-foundedness 
in the search of the $m_i$.  Is this  o.k.?
We want the $m_i$ to be not retrievable from the conditions, because $\gb$ shall stay small.}

Only parts ($\gamma$) and ($\delta$) are not so short:

($\gamma$)
Suppose that $p,q$ are incompatible in ${\bbforcing}({\mathcal C}^1)$, 
and, in the first case,
 assume that 
$n^p \leq n^q$. The first subcase is $\bar{\gc}^p \restriction n^p \neq \bar{\gc}^q \restriction n^p$
or that they are equal, but $\bar{\gc}^q \restriction [n^p, n^q) \not\geq_{\rm full }
\bar{\gc}^p \restriction [n^p,\infty)$.
Then, of course, this is also a reason for incompatibility in ${\bbforcing}({\mathcal C}^2)$.
The second  subcase is $\bar{\gc}^p \restriction n^p = \bar{\gc}^q \restriction n^p$
and $\bar{\gc}^q \restriction [n^p, n^q) \geq_{\rm full }
\bar{\gc}^p \restriction [n^p,\infty)$, but
for all $\bar{\gc} \in {\mathcal C}^1$,
$\bar{\gc}^p \restriction [n^p, \omega) \not\leq_{\rm full } \bar{\gc} \restriction [n^p,\omega)$
or $\bar{\gc}^q \restriction [n^q, \omega) \not\leq_{\rm full } \bar{\gc} \restriction [n^q,\omega)$.
By the premise $\bar{\gc}^p \leq_{\bf q} \bar{\gc}^q$, there $i(*)$ and $n_{i(*)}$ such that
$\bar{\gc}^p \restriction [n_{i(*)}, \omega) \leq_{\rm full } \bar{\gc}^q 
\restriction [i(*), \omega)$,
and $i(*), n_{i(*)} \geq n^p$. By forgetfulness, there are $m_1,m_2 > n^q$ and
 $\bar{\gd}$ such that $\mdn(\gc^p_{m_1}) = \mdn(\gc^q_{m_2})$ and
$$
\bar{\gd} \in \Sigma(\bar{\gc}^p \restriction [n^q, m_1) \cap
\Sigma(\bar{\gc}^q \restriction [n^q, m_2) \neq \emptyset.$$
Then $\bar{\gc} = (\bar{\gc}^q \restriction n^q\concat 
\bar{\gd} \concat \bar{\gc}^q \restriction [m_2, \omega)
\in {\mathcal C}^1$ shows that $p,q$ are compatible in ${\bbforcing}({\mathcal C}^1)$.
Contradiction.
The case $n^q < n^p$ is similar.

\smallskip

($\delta$) 
Again, as in ($\gamma$), we may assume
that $n^p \leq n^q$ and that  we are in the second subcase.
We show that if $(n^p,\bar{\gc}^p)$ and  $(n^q,\bar{\gc}^q)$ are compatible in
${\bbforcing}({\mathcal C}^2)$, then they are compatible in
${\bbforcing}({\mathcal C}^1)$.

Since ${\mathcal C}^1$ is $\bf q$-directed, there
is some $\bar{\gc} \in {\mathcal C}^1$ such that
$\bar{\gc}^p,\bar{\gc}^q \leq \bar{\gc}$.
By $\bar{\gc}^p \leq_{\bf q} \bar{\gc}$, there $i_p(*)$ and $n_{i_p(*)}$ such that
$\bar{\gc}^p \restriction [n_{i_p(*)}, \omega) \leq_{\rm full } \bar{\gc} 
\restriction [i_p(*), \omega)$,
and $i_p(*), n_{i_p(*)} \geq n^p$.
By $\bar{\gc}^q \leq_{\bf q} \bar{\gc}$, there $i_q(*)$ and $n_{i_q(*)}$ such that
$\bar{\gc}^p \restriction [n_{i_q(*)}, \omega) \leq_{\rm full } \bar{\gc} 
\restriction [i_q(*), \omega)$,
and $n_{i_p(*)} \geq n^p$, w.l.o.g.\ $n_{i_q(*)} \geq n^q$.

By the premise there in some $(n,\bar{\gd}) \in {\bbforcing}({\mathcal C}^2)$
such that $(n,\bar{\gd}) \geq (n^p, \bar{\gc}^p),
(n^q,\bar{\gc}^q)$.
So $\bar{\gd} \restriction [n^p,\infty) \geq_{\rm full} \bar{\gc}^p \restriction [n^p,\infty)$
and  $\bar{\gd} \restriction [n^q,\infty) 
\geq_{\rm full} \bar{\gc}^q \restriction [n^q,\infty)$.
So in particular there is some $m<\omega$ such that
$\bar{\gd} \restriction [n^p,m) \geq_{\rm full} \bar{\gc}^p \restriction [n^p,n_{i_q(*)})$
and  $\bar{\gd} \restriction [n^q,m) \geq_{\rm full} \bar{\gc}^q \restriction [n^q,n_{i_q(*)})$.
Then $((\bar{\gd} \restriction m) \concat \bar{\gc}) \in {\mathcal C}^1$ since
${\mathcal C}^1$ is closed under $\equiv_{\bf q}$ in ${\bf V}$ and all finite sections
of $\bar{\gd}$ are in ${\bf V}$.
  $(m, (\bar{\gd} \restriction m) \concat \bar{\gc}  
\restriction [i_q(*),\omega))$ is a stronger in 
${\bbforcing}({\mathcal C}^1)$ than $(n^p,\bar{\gc}^p)$ and $(n^q,\bar{\gc}^q)$.
\medskip 

6. For ${\mathcal C}^1 \subseteq {\mathcal C}^2$ being $\bf q$-pre-directed, we can
choose $D_1 \subseteq D_2$. In the limit we take the union of the proper filters,
which is again a proper filter, i.e., all intersections of
finitely many of its members are infinite and members of the filter. 
Any finitely many elements from 
$\bigcup_{\alpha < \lambda} {\mathcal C}^\alpha$ have an upper bound whose support
is in the filter.

For the ${\mathcal C}^\alpha$'s being $\bf q$-directed and increasing, it is 
clear that their union is ${\bf q}$-directed.
\proofend

\nothing{\begin{remark}\label{1.14a}
We shall use numerous fact we proved here for the freezing initial segments forcings
${\bbforcing}_{\bf q}({\mathcal C})$ also for the
true creature forcings ${\bf Q}_{\bf q}^{tr}({\mathcal C})$.
The translation is as follows: $(n,\bar{\gc}) $ will be 
replaced by some $(w,\bar{\gc})$ for some $w \in
\pos(\bar{\gc} \restriction n)$. There are finitely many possibilities.
Claim~\ref{1.14} parts 1 4 5, the analogous version of Definition~\ref{1.15} and
the analogous facts to \ref{1.16} to \ref{1.20} hold 
for ${\bf Q}_{\bf q}^{tr}({\mathcal C})$.
\end{remark}}

\begin{definition}\label{1.15}
Let $\bf q$ be a forgetful pcrf.
\begin{myrules}

\item[1.] Let ${\mathcal C}$ be $\bf q$-directed, $\bar{\gd} \in {\mathcal C}$ and 
${\mathcal J} \subseteq {\bbforcing}_{\mathcal C}(\bar{\gd})$ be a maximal antichain of 
${\bbforcing}({\mathcal C})$. By induction on the ordinal $\alpha$ we define when 
$
\rk_{\bf q}(n,m,\bar{\gc},\bar{\gd},{\mathcal J}, {\mathcal C}) =\alpha$
for $\bar{\gc} \in {\bf Q}_{\mathcal C}(\bar{\gd})$ and $n\leq m <\omega$.

Case: $\alpha =0$. $\rk_{\bf q}(n,m,\bar{\gc},\bar{\gd},{\mathcal J}, {\mathcal C}) =0$
iff $(n,\bar{\gc}) \in {\bbforcing}(\mathcal C)$ is above some member of ${\mathcal J}$.

Case:  $\alpha >0$. $\rk_{\bf q}(n,m,\bar{\gc},\bar{\gd},{\mathcal J}, {\mathcal C}) =\alpha$
iff
\begin{myrules}

\item[(a)] for no $\beta < \alpha$ we have
$\rk_{\bf q}(n,m,\bar{\gc},\bar{\gd},{\mathcal J}, {\mathcal C}) =\beta$,

\item[(b)] for all $p'= (m',\bar{\gd}') \in {\mathcal J}$, $\beta <\alpha$,
exists $\bar{\gc}'$ such that
\begin{myrules}

\item[($\alpha$)] $\bar{\gc}' \in {\bf Q}_{\mathcal C}(\bar{\gd})$  (hence $\bar{\gc}'
\leq_{\bf q} \bar{\gd}'$ and vice versa and same with
$\bar{\gc}$, $\bar{\gd}$, because all of them are $\equiv_{\bf q} \bar{\gd}$),

\item[($\beta$)] $\bar{\gc}\restriction m =
 \bar{\gc}' \restriction m$,

\item[($\gamma$)] $(n,\bar{\gc}')$ and $(m',\bar{\gd}')$ are incompatible in
${\bbforcing}_{\mathcal C}(\bar{\gd})$, equivalently (---by (b)($\alpha$) and 
Claim~\ref{1.14}(5) ($\delta$)---) in ${\bbforcing}({\mathcal C})$,

\item[($\delta$)] 
$\beta \leq \rk_{\bf q}(n,m+1,\bar{\gc}',\bar{\gd},{\mathcal J}, {\mathcal C})$.
\end{myrules}
\end{myrules}

\item[2.] $\rk_{\bf q}(\bar{\gd},{\mathcal J}, {\mathcal C}) = \bigcup \{
\rk_{\bf q}(n,m,\bar{\gc},\bar{\gd},{\mathcal J}, {\mathcal C}) +1 \such 
\bar{\gc} \in {\bf Q}_{\mathcal C}(\bar{\gd}), n\leq m<\omega \}$.

\item[3.] We say ${\mathcal C}$ is ${\bf q}$-nice if for all $\bar{\gd}\in K$, ${\mathcal J}
\subseteq {\bbforcing}_{\mathcal C}(\bar{\gd})$
that are predense in ${\bbforcing}({\mathcal C})$,
$\rk_{\bf q}(\bar{\gd},{\mathcal J},{\mathcal C}) < \omega_1$.
Remark: Since ${\bbforcing}(\mathcal C)$ is $\sigma$-centred and since
countable sets have 
$\leq_{\bf q}$ upper bounds, there are such $\bar{\gd}$'s and ${\mathcal J}$'s.

$\forall {\mathcal J} \rk_{\bf q}(\bar{\gd},
{\mathcal J},{\mathcal C}) < \omega_1$
is absolute: 
$(\forall {\mathcal J} \subseteq {\bbforcing}_{\mathcal C}(\bar{\gd}))
(\exists \alpha < \omega_1) 
(\rk_{\bf q}(\bar{\gd},{\mathcal J},{\mathcal C})= \alpha)$
is $\Pi_2^1$ and hence Shoenfield's absoluteness theorem applies. 
Also quantification over
all $\bar{\gd}$ is possible. 
But if ${\mathcal C}$ is expanded, the rank can become infinite.

\item[4.] We say $\mathcal J \subseteq {\bbforcing}_{\mathcal C}(\bar{\gd})$ 
is an explicitly predense set (or maximal antichain in
${\bbforcing}({\mathcal C})$)
if its members are pairwise incompatible and if 
$\rk_{\bf q}(\bar{\gd},{\mathcal J},{\mathcal C}) < \omega_1$. 
\end{myrules}
\end{definition}

Now, with the help of the rank function, we establish criteria for $\mathcal J$
to stay predense in ${\bbforcing}({\mathcal C})$
in any extension of the universe and under any extension
of $\mathcal C$. As usual, we write ${\bf V}^{\mathbb P}$ for ${\bf V}[G]$
with an arbitrary ${\mathbb P}$-generic $G$ over ${\bf V}$.

\begin{claim}\label{1.16}
Assume
\begin{myrules}
\item[(a)] ${\mathcal C}$ is a $\bf q$-directed family.
\item[(b)] $\bar{\gd} \in {\mathcal C}$, and ${\mathcal J} \subseteq 
{\bbforcing}_{\mathcal C}(\bar{\gd})$ is predense in ${\bbforcing}({\mathcal C})$.
\end{myrules}

Then the following six conditions are equivalent:
\begin{myrules}
\item[($\alpha$)] In ${\bf V}^{\rm Cohen}$ there is
some $\bar{\gd}'$ such that $\bar{\gd} \leq_{\bf q} \bar{\gd}' \in {\bf Q_q}^{{\bf V}^{\rm Cohen}}$,
${\mathcal C} \cup \{ \bar{\gd}'\}$ is a ${\bf q}$-pre-directed family, and for
some $n<\omega$,
$(n,\bar{\gd}')$ is incompatible in ${\mathbb Q}({\mathcal C} \cup 
\{\bar{\gd}'\})$ with every member of ${\mathcal J}$.

\item[($\beta$)] Like ($\alpha$) for some atomless forcing notion ${\mathbb P}$.

\item[($\beta'$)] For some atomless forcing notion
${\mathbb P}$, in ${\bf V}^{\mathbb P}$  there are
some $\bar{\gd}'$ and some $n$ and some
$\bf q$-directed ${\mathcal C}'$ such that
${\bf Q}_{\mathcal C}(\bar{\gd}) \subseteq {\mathcal C}'$,
 such that $\bar{\gd} \leq_{\bf q} \bar{\gd}' \in 
{\mathcal C}'$,
and in ${\bbforcing}({\mathcal C}')$ the condition
$(n,\bar{\gd}')$ is incompatible with every member of ${\mathcal J}$.

\item[($\gamma$)] For some atomless forcing notion
${\mathbb P}$, in ${\bf V}^{\mathbb P}$  there is a strongly directed
${\mathcal C}'\supseteq {\mathcal C}$, such that in
${\bbforcing}({\mathcal C}')$, the set ${\mathcal J}$ is not a maximal antichain.

\item[($\delta$)]
$\rk_{\bf q}(\bar{\gd},{\mathcal J},{\mathcal C}) = \infty$.

\item[($\eps$)]
$\rk_{\bf q}(\bar{\gd},{\mathcal J},{\mathcal C}) \geq \omega_1$.
\end{myrules}
\end{claim}

\proof The proof is via a circle ($\eps$) $\Rightarrow$ ($\alpha$) $\Rightarrow$ ($\beta$)
$\Rightarrow$ ($\beta'$) $\Rightarrow$ ($\gamma$) $\Rightarrow$ ($\delta$) $\Rightarrow$ ($\eps$).

\smallskip

($\eps$) $\Rightarrow$ ($\alpha$). $\rk_{\bf q}(\bar{\gd},{\mathcal J}, {\mathcal C})  
\geq \omega_1$.
We choose $\bar{\gd}'' \restriction m $ by induction on $m$.
Let ${\mathcal J} = \{ p'_k\such k \in \omega \}$ be an enumeration of ${\mathcal J}$.
We take $n,m,\bar{\gc}$ such that
$\rk_{\bf q}(n,m,\bar{\gc},\bar{\gd},{\mathcal J},{\mathcal C}) \geq \omega_1$.
Since $\omega \times\omega\times{\bf Q}_{\mathcal C}(\bar{\gd})$ is countable
and since $\omega_1$ is regular, there is such a triple.
Let ${\mathbb P}$ be the forcing adding a Cohen real.
Suppose that $m_0, \bar{\gc}^0, p_k'$ are given and  that
$\omega_1 = \rk_{\bf q}(n,m_0,\bar{\gc}^0,\bar{\gd},{\mathcal J},{\mathcal C})$.
Let $D_{p'_k,m_0,\bar{\gc}^0}$ be the set of 
$\bar{\gc}^1$ such that
\begin{myrules}
\item[(1)] $\bar{\gc}^1 \in {\bf Q}_{\mathcal C}(\bar{\gd})$,
\item[(2)]$\bar{\gc}^1\restriction m_0 = \bar{\gc}^0 \restriction m_0$,
\item[(3)]$(n,\bar{\gc}^1)$ and $p'_k$ are incompatible in 
${\bbforcing}_{\mathcal C}(\bar{\gd})$ and the incompatibility is witnessed
below $m_0+\ell$ for some $\ell$
and (--- so, despite forgetfulness, there are
also many elements outside $\Sigma(\bar{\gc}^1 \restriction [k,m))$
that are in ${\bbforcing}_{\mathcal C}(\bar{\gd})$ and incompatible with $p_k'$---)
\item[(4)] $\forall \beta < \omega_1$, 
$\beta \leq \rk_{\bf q}(n,m_0+\ell,\bar{\gc}^1,\bar{\gd}, {\mathcal J},{\mathcal C})$.
\end{myrules}
Reading definition~\ref{1.15}(b) $\ell$ times shows that $D_{p'_k,m_0,\bar{\gc}^0}$
is is dense in the Cohen forcing.
So, if $G$ is ${\bbforcing}_{\mathcal C}(\bar{\gd})$-generic over $\bf V$ and
$\bar{\gd}' = \bigcup\{ \bar{\gc}\restriction n \such (n,\bar{\gc}) \in G\}$,
then  $(n,\bar{\gd}')$ is incompatible in
${\bbforcing}({\mathcal C} \cup \{\bar{\gd}'\})$ with every member of ${\mathcal J}$.

\smallskip

($\beta$) $\Rightarrow$  ($\beta'$): by Claim~\ref{1.14} Part 5 ($\gamma$):
Incompatibility of $\bar{\gd}'$ with every member of
${\mathcal J}$ is the same in ${\bbforcing}({\mathcal C}\cup \{\bar{\gd}'\})$ 
and in ${\bbforcing}({\mathcal C}')$, because every member of ${\mathcal J}$ is 
$\leq_{\bf q} \bar{\gd}'$.

\smallskip

($\beta'$) $\Rightarrow$  ($\gamma$):
\nothing{
$2^{{\aleph_0}\cdot ac({\mathbb P})}$ 
($ac({\mathbb P})$ shall be a cardinal such that all  antichains of $\mathbb P$
are in size less than or equal to this cardinal) iterations of the forcing ${\mathbb P}$ 
from ($\beta'$) gives us a strongly directed ${\mathcal C}' \supseteq {\mathcal C}
\cup\{\bar{\gd}'\}$.}
We take a model from ($\beta'$) and extend ${\mathcal C}'$ from that model to a strongly 
directed ${\mathcal C}''$. This is possible because $\leq_{\bf q}$ increasing 
countable chains have an upper bound in ${\bf Q}$. Since
${\mathcal C}'$  is directed and all intermediate extensions can chosen to be
directed, it is sufficient to add limits to chains instead of
upper bounds to  countable subsets.

\smallskip

($\gamma$) $\Rightarrow$ ($\delta$)
$\rk_{\bf q}(n,m,\bar{\gd}' \restriction m ,\bar{\gd},{\mathcal J}, {\mathcal C})$ 
is  witness to a descending chain in the ranks. Since $K \subseteq {\mathcal H}(\aleph_0)$
also for all new $\bar{\gd}'$'s, every $\bar{\gd}' \restriction m \in {\bf V}$.
\proofend

\begin{conclusion}\label{1.17}
1. If ${\mathbb P}$ is the forcing for adding $2^{\aleph_0}$ Cohen reals, then 
for any $\bf q$-directed ${\mathcal C}^1 \in {\bf V}$,
 in ${\bf V}^{\mathbb P}$ we can extend ${\mathcal C}^1 $
to a $\bf q$-nice ${\mathcal C}^2$.

2. Assume ${\mathbb P}_1 \lessdot {\mathbb P}_2$ and
 $\Vdash_{{\mathbb P}_1} $ ``$\name{{\mathcal C}^1} $ is $\bf q$-nice'' and 
$\Vdash_{{\mathbb P}_2} $ ``$\name{{\mathcal C}^2} $ is $\bf q$-directed 
and extends  $\name{{\mathcal C}^1}$''.
Then ${\mathbb P}_1 * {\bbforcing}(\name{{\mathcal C}^1}) \lessdot
{\mathbb P}_2 * {\bbforcing}(\name{{\mathcal C}^2})$.

3. If ${\mathbb P}_2/{\mathbb P}_1$ contains as a complete 
subforcing the forcing for adding $2^{\aleph_0}$ Cohen reals, 
then there is $\name{{\mathcal C}^3}$
such that $\Vdash_{{\mathbb P}_2}$ ``${\mathcal C}^3$ is ${\bf q}$-nice 
and extends $\name{{\mathcal C}^1}$'', and 
${\mathbb P}_1 * {\bbforcing}(\name{{\mathcal C}^1}) \lessdot
{\mathbb P}_2 * {\bbforcing}(\name{{\mathcal C}^3})$.

\end{conclusion}

\proof 
1.  
Let $\la (\bar{\gd}_\alpha, {\mathcal J}_\alpha) \such \alpha < 2^{\aleph_0 } \r$
be an enumeration such that
for all $\beta < 2^{\aleph_0 }$, 
$\la (\bar{\gd}_\alpha, {\mathcal J}_\alpha) \such \alpha \geq \beta \r$
enumerates among others  all of 
$ \bigcup \{ K^{{\bf V}[\bar{\gge}]} \times {\mathcal P}(\bbforcing_{{\mathcal C}) \cup
\rge\bar{\gge}}(\bar{\gd}) \such \bar{\gge} \in V[{\rm Cohen}_\beta] \cap 
(\omega^\omega)^\omega, \bar{\gd} \in {\mathcal C}\cup
 \rge{\bar{\gge}}\}$.
Then for step $\alpha< 2^{\aleph_0}$  of an iterative construction take
$(\bar{\gd}_\alpha, {\mathcal J}_\alpha)$ check whether
($\alpha$) of claim 1.16 is true for $(\bar{\gd}_\alpha,{\mathcal J}_\alpha,
{\mathcal C}_\alpha)$, and if so, then add by Cohen forcing $\bar{\gd}'$ as there, and define
${\mathcal C}_{\alpha+1} = {\mathcal C}_\alpha \cup\{\bar{\gd}'\}$.

\smallskip

2. Similar to \cite[Claim 4.7]{Sh:700} or to \cite[Section 5, Claim 5.6]{Sh:707}.
For completeness, we give the proof:
Since we may think of ${\mathbb P}_2$ as ${\mathbb P}_2 / {\mathbb P}_1$, we
can replace ${\mathbb P}_1$ by the trivial forcing. 
Then ${\mathcal C}^1$ is in the ground model. Let $G$ be a ${\mathbb P}_2$-generic 
filter and let ${\mathcal C}^2 = \name{\mathcal C}^2[G]$.
If $p,q \in {\bbforcing}({\mathcal C}^1)$ and $p \leq q$, then they are 
also in ${\bbforcing}({\mathcal C}^2)$ in this order. 
The same holds for incompatibility by  Claim~\ref{1.14}(5).
Now we need to show that maximal antichains in ${\bbforcing}({\mathcal C}^1)$ stay maximal in 
${\bbforcing}({\mathcal C}^2)$:
For this take a maximal antichain ${\mathcal J}$. It is countable and a subset of
${\bbforcing}_{{\mathcal C}^1}(\bar{\gd})$ for some suitable $\bar{\gd}
\in {\bf V}$. Now we
have by ${\mathcal C}^1$ being nice that
$\rk(\bar{\gd},{\mathcal J},{\mathcal C}^1) < \omega_1$.  
This is absolute and holds in $V[G]$ as well. 
So no member of ${\bbforcing}({\mathcal C}^2)$ can be orthogonal to  a member of 
${\mathcal J}$. 

3. This is just the combination of the first two parts.
\proofend

\begin{definition}\label{1.18}
Assume that ${\mathcal C}$ is $\bf q$-directed.
Let $ \name{G}_{{\bbforcing}(\mathcal C)}$ be a name for
the ${{\bbforcing}(\mathcal C)}$-generic filter $G$.
$\name{\bar{\gc}}_{{\bbforcing}(\mathcal C)}$ is the
${\bbforcing}(\mathcal C)$-name
$$\bigcup \{ \bar{\gc} \restriction n \such (n,\bar{\gc}) \in 
\name{G}_{{\bbforcing}(\mathcal C)}\}.
$$
\end{definition}

\begin{claim}\label{1.19}
Assume that ${\mathcal C}$ is $\bf q$-directed.

1. $\Vdash_{{\bbforcing}(\mathcal C)} \mbox{``}\name{\bar{\gc}} \in {\bf Q} 
\mbox{ is a } \leq_{\bf q} \mbox{-upper bound of } {\mathcal C}$''.

\smallskip

2. If $\mathcal C$ is a $\bf q$-nice family, then $\name{\bar{\gc}}$ is a generic
real for ${\bbforcing}(\mathcal C)$, that
is, for $G \subseteq {\bbforcing}(\mathcal C)$ generic over ${\bf V}$, letting
$\bar{\gc}^* = \name{\bar{\gc}}[G]$, we have
\begin{equation}\label{generic}
 G = \{ p \in {\mathbb Q}(\mathcal C) \such
\bar{\gc}^p \restriction n^p = \bar{\gc}^* \restriction n^p
\mbox{ and } \bar{\gc}^p \restriction [n^p,\infty) \leq_{\rm full } \bar{\gc}^* \restriction
 [n^p,\omega)\}.
\end{equation}
\end{claim}

\proof

1. For all $\bar{\gd} \in {\mathcal C}$, the set
$\{ (n^p,\bar{\gc}^p) \in {\bbforcing}(\mathcal C) \such \bar{\gc}^p \geq_{\bf q} \bar{\gd}\}$ is
open and 
dense in ${\bbforcing}(\mathcal C)$, because ${\mathcal C}$ is $\bf q$-directed,
and because $\bar{\gc}^p \geq_{\bf q} \bar{\gd}$ and $(n^q,\bar{\gc}^q) \geq (n^p,\bar{\gc}^p)$
together implies $\bar{\gc}^q \geq_{\bf q} \bar{\gd}$.

\smallskip

2. First we show that $G$ is a subset of the right hand side.
Let $p \in G$. Then $\bar{\gc}^p \restriction n^p \triangleleft \bar{\gc}^*$.
Also for every $m \geq n^p$ there is some $q \in G$ such that $n^q \geq m$
and $q \not\perp p$ hence $\exists m'\geq m$
$\bar{\gc}^p \restriction [n^p,m') \leq_{\rm full} \bar{\gc}^q \restriction [n^p,n^q)
= \bar{\gc}^* \restriction[n^p,n^q)$.
Hence $\bar{\gc}^p \restriction [n^p,\omega) \leq_{\rm full} \bar{\gc}^* \restriction [n^p,\omega)$.

\smallskip

Now for the other direction. Any two elements of the right hand side have
$\bar{\gc}^*$ as a common stronger element.
Hence the right hand side is directed and  hence $G$ cannot be a proper subset of it.
\proofend

\nothing{
We assume that ${\mathcal C}$ is strongly $\bf q$-directed.
If this is not the case, then we use ${\mathcal C}$'s niceness to
extend it to a strongly directed ${\mathcal C}' \supseteq {\mathcal C}$
such that ${\mathcal J}$ stays a maximal antichain.
Then we carry out the proof for ${\bf Q}({\mathcal C}')$.
Intersection both sides of equation~\eqref{generic} will give
the desired result.
\begin{definition}\label{1.20}
(is 1.3 of Mathias)
\begin{multline*}
{\mathcal B}=
\{ U \subseteq {\bf Q}(\mathcal C) \such \forall (n,\bar{\gge}) \in {\bbforcing}(\mathcal C)
\exists (n,\bar{\gge}') \geq (n,\bar{\gge}) \\
 (n,\bar{\gge}')   \Vdash_{{\bbforcing}(\mathcal C)}
\name{\bar{\gc}} \in U \vee
 (n,\bar{\gge}')   \Vdash_{{\bbforcing}(\mathcal C)}
\name{\bar{\gc}} \not\in U )\}.
\end{multline*}
\end{definition}
\begin{proposition}\label{1.21}
${\mathcal B}$ contains all open sets in the Cantor topology.
\end{proposition}
\proof There forgetfulness is used instead of the ultrafilter.
Now the proof for strongly directed ${\mathcal C}$ is
similar to \cite[Theorem 2.0]{mathias:happy}:
We want to characterize ${\bbforcing}({\mathcal C})$-generic sequences of creatures
by properties in $\bf Q$.
For this we define
\begin{definition}\label{1.22}
This is like Definition 2.2. Mathias.
Let $\Delta$ be a dense open subset of ${\bbforcing}(\mathcal C)$.
Let $\bar{\gd}$ be a finite sequence in $K$.
We say $\bar{\gc} \in {\bf Q}$ captures $(\bar{\gd}, \Delta)$ if
$\bar{\gd} = \bar{\gc} \restriction \lg(\bar{\gd})$ and 
\begin{equation}
\begin{split}
&(\forall \bar{\gc}' \geq_{\bf q} \bar{\gc})\mbox{ witnessed by } (i,n))\\
&( \exists \bar{\gd}' \leq_{\rm full} \bar{\gc}' \restriction
[\lg(\bar{\gd}) ,\infty))\,(\exists k)\,
(\lg(\bar{\gd})+k,
\bar{\gd} \concat \bar{\gd}'\restriction k \concat \bar{\gc} \restriction [k,\infty)) \in \Delta.
\end{split}
\end{equation}
\end{definition}
\begin{lemma}\label{1.21}
Let $\Delta$ Be a dense closed subset of ${\bbforcing}({\mathcal C})$ and $\bar{\gd}\in ^m K$.
Then
 there is a $\bar{\gc}$ that captures $(\bar{\gd}, \Delta)$.
\end{lemma}
\proof (this is like Prop.\ 2.3 in Mathias)
For each $\bar{\gd}' \triangleright \bar{\gd}$, be a finite sequence in $K$ we
choose $\bar{\gc}_i =\bar{\gc}_{\bar{\gd}} \in {\bf Q}$ that captures $(\bar{\gd}, \Delta)$ if
possible such that $\bar{\gc}_{i+1} \geq_{\bf q} \bar{\gc}_i$.
We let $\bar{\gc}$ be $\leq_{\bf q}$ above all the $\bar{\gc}_{\bar{\gd}}$.
$$\{ \bar{\gge} \in {\bf Q}(\mathcal C) \such
(m,\bar{\gc}) \leq_{\bbforcing(\mathcal C)}
 (m,\bar{\gge}) \rightarrow \exists k>m (k, \bar{\gge}) \in \Delta \}.$$
This is open in the Cantor topology, hence by the previous proposition there in an
$\bar{\gh} \in {\mathcal C}$ that $(m,\bar{\gh})$ decides it.
Then it decides in the positive sense then $\bar{\gh} \Vdash \name{\bar{\gc}} \in B$
\proofend
End of the proof of the hard direction:
Let $p=(n^p,\bar{\gc}^p)$ be an element of the right hand side.
Now we take a countable model $M$ of a sufficiently large part of ZFC.
Let $\Delta \in M$ be a dense closed subset of $({\bbforcing}({\mathcal C}))^M$.
Working in $M$, pick for each finite $\bar{\gd}$ some
$\bar{\gc}_{\bar{\gd}} \in  {\bf Q}({\mathcal C})$ that captures
$(\bar{\gd},\Delta)$. 
Do it inductively so that the $\bar{\gc}_{\bar{\gd}}$, $\bar{\gd} \in {}^{\omega>} K$ 
are $\leq_{\bf q}$-increasing.
Since ${\mathcal C}$ is strongly ${\bf q}$-directed, there is
some $\bar{\gc}$ that is $\leq_{\bf q}$-stronger than all the $\bar{\gc}_{\bar{\gd}}$.
$\bar{\gc}\leq_{\bf q} \bar{\gc}^*.$
So in $M$
\begin{equation}\label{innen}
\forall \bar{\gc}' \leq_{\bf q} \bar{\gc} \exists
\bar{\gd}' \geq_{\rm full} \bar{\gc}'(\lg(\bar{\gd}) + \lg(\bar{\gd'}), \bar{\gd}
\concat \bar{\gd}' \concat \bar{\gc}\restriction[n,\infty)) \in \Delta.
\end{equation}
Then this is equivalent to
$$
\forall \bar{\gc}' ((H,\prec)_{\bar{\gc}'} \mbox{  is well founded}),$$
 where
$$H_{\bar{\gc}'}= \{ \bar{\gd}' \such
\bar{\gd} \geq_{\rm full} 
\bar{\gc}'
\wedge  (\lg{\bar{\gd}}+\lg{\bar{\gd}'},
 \bar{\gd} \concat \bar{\gd}'\concat \bar{\gc} \restriction [k,\infty))
\not\in \Delta\},$$
ordered by end-extension of finite sequences.
Then $(H,\prec)_{\bar{\gc}^*}$
is well-founded in $V$.
\begin{equation}\label{aussen}
\exists
\bar{\gd}' \geq_{\rm full} \bar{\gc}^*(\lg(\bar{\gd}) + \lg(\bar{\gd'}), \bar{\gd}
\concat \bar{\gd}' \concat \bar{\gc}\restriction[n,\infty)) \in \Delta.
\end{equation}
\proofend}

\begin{claim}\label{1.20}
Assume 
\begin{myrules}
\item[(a)] $\bf q$ is a pcrf,
\item[(b)] ${\mathbb P}_0 \lessdot {\mathbb P}_2$, 
\item[(c)] for $\ell =0,2$ we have that $\name{\mathcal C}^\ell$ is a ${\mathbb P}_\ell$-name
of a $\bf q$-nice family,
\item[(d)] for $\ell =0,2$ let  ${\mathbb P}_{\ell +1}= {\mathbb P}_\ell
* {\bbforcing}(\name{\mathcal C}^\ell) $, and let $\name{\bar{\gc}}^\ell$
be a ${\mathbb P}_{\ell +1}$-name of the generic for
${\bbforcing}(\name{\mathcal C}^\ell)$,
\item[(e)] $\Vdash_{{\mathbb P}_2} \mbox{``}\name{\mathcal C}^0 \subseteq \name{\mathcal C}^2$''.
\end{myrules}

Then 
\begin{myrules}
\item[($\alpha$)] ${\mathbb P}_2 \Vdash \name{\bar{\gc}}^2$ 
is a generic for ${\bbforcing}(\name{\mathcal C}^0)$.
\item[($\beta$)] ${\mathbb P_1} \lessdot {\mathbb P}_3$.
\end{myrules}
\end{claim}

\proof
See \ref{1.17} (2). ${\mathbb P}_0 * {\bbforcing}(\name{{\mathcal C}^0})
\lessdot {\mathbb P}_2 * {\bbforcing}(\name{{\mathcal C}^0})
\lessdot {\mathbb P}_2 * {\bbforcing}(\name{{\mathcal C}^2})$.
\proofend

\section{Forcing bigness notions}

See ``Vive la diff\'{e}rence III'' \cite{Sh:509} for a connection to model theory.

\begin{definition}\label{2.1}
1. A forcing bigness notion is a class $\Gamma$, partially ordered by $\leq_{\Gamma}$, such that
 for all  $({\mathbb P}_0, {\mathbb P}^+_0, \name{\eta}) \in \Gamma$  we have
\begin{myrules}
\item[(a)] ${\mathbb P}_0 \lessdot {\mathbb P}_0^+$ are forcing notions,
\item[(b)] $\name{\eta}$ is a ${\mathbb P}_0^+$-name of a sequence of ordinals
(usually $\in {}^\omega \omega$),
\item[(c)] if ${\mathbb P}\lessdot {\mathbb P}_0$ and ${\mathbb P}^+_0 \lessdot {\mathbb P}^+$,
then  $({\mathbb P}, {\mathbb P}^+, \name{\eta}) \in \Gamma$.
\end{myrules}
If $({\mathbb P}_0, {\mathbb P}_0^+, \name{\eta}) \in \Gamma$, we may say that 
 $({\mathbb P}_0, {\mathbb P}_0^+, \name{\eta})$ is $\Gamma$-big.

\smallskip

2. $\leq_{\Gamma}$ is a quasi order and $\lessdot_{\Gamma} \supseteq \, \leq_{\Gamma}$
where $\lessdot_{\Gamma}$ is the following (quasi) partial order:
$({\mathbb P}_1, {\mathbb P}^+_1,\name{\eta}_1) \lessdot_{\Gamma} 
({\mathbb P}_2, {\mathbb P}^+_2,\name{\eta}_2)$
means that both triples are in $\Gamma$, and
${\mathbb P}_1 \lessdot {\mathbb P}_2$ and ${\mathbb P}^+_1 \lessdot {\mathbb P}_2^+$
and $\name{\eta}_1 = \name{\eta}_2$ and such that
the following diagramme is exact, which means, that
for any $p_1^+ \in {\mathbb P}_1^+ $, $p_2 \in {\mathbb P}_2$, 
\begin{equation}
\label{raute}
\sup_{{\mathbb P}_1} \{ p_1 \in {\mathbb P}_1 \such 
p_1 \perp_{{\mathbb P}^+_1}  p_1^+ \}
=\sup_{{\mathbb P}_2} \{ p_2 \in {\mathbb P}_2 \such p_2 
\perp_{{\mathbb P}_2^+ }  p_1^+ \}.
\end{equation}
$$
\xymatrix{
{\mathbb P}_1^+ \ar[r]^{\lessdot} & {\mathbb P}_2^+\\
 {\mathbb P}_1 \ar[r]^{\lessdot} \ar[u]^\lessdot & {\mathbb P}_2 \ar[u]^\lessdot
}
$$

\smallskip

3. We say $\Gamma$ is a c.c.c.\ forcing bigness notion, if
in addition to the above mentioned properties
the first and second components of all triples in $\Gamma$ are
c.c.c\ forcing notions and if 
\begin{myrules}
\item[(d)] $\langle ({\mathbb P}_i, {\mathbb P}^+_i, \name{\eta}) \such i < \delta \rangle
 \in \Gamma$ 
is continuously increasing in $\leq_{\Gamma}$, then
$({\bigcup \{\mathbb P}_i \such i < \delta\}, \bigcup\{{\mathbb P}^+_i \such i < \delta\},
\name{\eta}) \in \Gamma$ and is $\leq_{\Gamma}$-above all 
$({\mathbb P}_i, {\mathbb P}^+_i, \name{\eta})$.
\end{myrules}

\smallskip

4. For a c.c.c.\ forcing bigness notion $\Gamma$ we say that $\Gamma$ has amalgamation iff:
If $({\mathbb P}_1, {\mathbb P}_1^+, \name{\eta}) \in \Gamma$ and 
${\mathbb P}_1 \lessdot {\mathbb P}_2$ and ${\mathbb P}_2/{\mathbb P}_1$ 
satisfies the  Knaster property, then we can find ${\mathbb P}_2^+$ [satisfying the c.c.c.]
and $f$ such that
\begin{myrules}
\item[($\alpha$)] ${\mathbb P}_1^+ \lessdot {\mathbb P}_2^+$,
\item[($\beta$)] $f$ is a complete embedding of ${\mathbb P}_2$ into ${\mathbb P}_2^+$
over ${\mathbb P}_1$,
\item[($\gamma$)] $(f({\mathbb P}_2), {\mathbb P}_2^+, \name{\eta}) \in \Gamma$,
\item[($\delta$)] $({\mathbb P}_1, {\mathbb P}_1^+, \name{\eta})
\leq_{\Gamma} (f({\mathbb P}_2), {\mathbb P}_2^+, \name{\eta})$.
\end{myrules}

\nothing{
\smallskip
5. $\Gamma$ is a nice forcing bigness notion, iff it is a c.c.c.\ forcing bigness notion with 
amalgamation.}

\smallskip

5. For a c.c.c.\ forcing bigness notion $\Gamma$ we say that $\Gamma$ has free
amalgamation iff:
If $({\mathbb P}_1, {\mathbb P}_1^+, \name{\eta}) \in \Gamma$ and 
${\mathbb P}_1 \lessdot {\mathbb P}_2$ and ${\mathbb P}_2/{\mathbb P}_1$ 
satisfies the Knaster property, then we can choose in 4.\   ${\mathbb P}_2^+=
{\mathbb P}_1^+ \times_{{\mathbb P}_1} {\mathbb P}_2$ [satisfying the c.c.c.]
such that
\begin{myrules}
\item[($\alpha$)] ${\mathbb P}_1^+ \lessdot {\mathbb P}_2^+$,
\item[($\beta$)] $id$ is a complete embedding of ${\mathbb P}_2$ into ${\mathbb P}_2^+$
over ${\mathbb P}_1$,
\item[($\gamma$)] $({\mathbb P}_2, {\mathbb P}_2^+, \name{\eta}) \in \Gamma$.
\end{myrules}
For a definition of  ${\mathbb P}_1^+ \times_{{\mathbb P}_1} {\mathbb P}_2$ 
see \cite{JudahRoslanowski}.
\smallskip

6. We say  $\Gamma$ has the $(<\lambda)$-LS (L\"owenheim Skolem) property when
the following holds:
If $({\mathbb P}, {\mathbb P}^+, \name{\eta}) \in \Gamma$ and  
$Y \subseteq {\mathbb P}^+$, $|Y| < \lambda$ (or alternatively,
for stationarily many $Y \in [{\mathbb P}^+]^{<\lambda}$) we have:
There is
$({\mathbb P}_1, {\mathbb P}_1^+, \name{\eta})$ such that,
\begin{myrules}
\item[($\alpha$)] $({\mathbb P}_1, {\mathbb P}_1^+, \name{\eta}) \in \Gamma$,
\item[($\beta$)] $|{\mathbb P}_1^+| <\lambda$,
\item[($\gamma$)] ${\mathbb P} \cap Y \subseteq {\mathbb P}_1
\lessdot {\mathbb P}$ and 
\item[($\delta$)] $Y \subseteq {\mathbb P}_1^+ \lessdot{\mathbb P}^+$.
\end{myrules}

\smallskip

7. $\Gamma$ is a Knaster forcing bigness notion iff it 
is a c.c.c.\ forcing bigness notion and
in addition, if  $({\mathbb P},{\mathbb P}^+,\name{\eta}) \in \Gamma$
then ${\mathbb P}$ and ${\mathbb P}^+$  satisfy the Knaster condition.
\end{definition}

\begin{observation}\label{2.2}
$\lessdot_{\Gamma}$ is a partial order.
\end{observation}

\begin{definition}\label{2.3}
\begin{myrules}
\item[(1)]
We define the following $\Gamma$:
\begin{myrules}
\item[(a)] $\Gamma_{\rm new} = \{ ({\mathbb P},{\mathbb P}^+,\name{\eta}) \such
{\mathbb P} \lessdot {\mathbb P}^+$ are c.c.c.\ forcing notions and such that
$\Vdash_{{\mathbb P}^+} \name{\eta} \in ({}^\omega \omega)^{{\bf V}[G_{{\mathbb P}^+}]}
\setminus {\bf V}[G_{\mathbb P}] \}$.

\item[(b)] $\Gamma_{\rm Cohen} = \{ ({\mathbb P},{\mathbb P}^+,\name{\eta}) \in
\Gamma_{\rm new}\such
 \name{\eta} $ is forced by ${\mathbb P}^+$ to be a Cohen real
over ${\bf V}^{\mathbb P}\}$.

\item[(c)] $\Gamma_{\rm ud} = \{ ({\mathbb P},{\mathbb P}^+,\name{\eta}) \in
\Gamma_{\rm new}\such
 \name{\eta} $ is forced by ${\mathbb P}^+$ 
not to be dominated by any $\nu \in {\bf V}^{[G_{\mathbb P}]}\}$.
$\eta$ is dominted by $\nu$, abbreviated $\eta \leq^* \nu$, 
if $(\exists n)( \forall k \geq n)( \eta(k) \leq \nu(k))$.
\item[(d)] $\Gamma_{\rm tow} = \{ ({\mathbb P},{\mathbb P}^+,\name{\eta}) 
\in \Gamma_{\rm new} \such
{\mathbb P} \lessdot {\mathbb P}^+$ and 
$\name{\eta} \in ({}^\omega 2)^{{\bf V}[G_{{\mathbb P}^+}]}$
are such that $\eta^{-1}\{1\}$ is infinite and
is forced by ${\mathbb P}^+$  not to  contain any infinite subset
of $\omega$ which belongs to ${\bf V}^{\mathbb P}\}$

\end{myrules}

\item[(2)]
We define when
a name $\name{\eta}$ is $\Gamma$-big.

 \begin{myrules}
\item[(a)] We say $\name{\eta}$ is $\Gamma_{\rm new}$-big over $\nu\in {\bf V}[\mathbb P]$ if
$\Vdash_{{\mathbb P}^+} \name{\eta} \not\in {\bf V}[\nu]$.

\item[(b)] We say $\name{\eta}$ is $\Gamma_{\rm Cohen}$-big over $\nu\in {\bf V}[\mathbb P]$ if
$\Vdash_{{\mathbb P}^+} \name{\eta}$ Cohen-generic over $ {\bf V}[\nu]$.

\item[(c)] We say $\name{\eta}$ is $\Gamma_{\rm ud}$-big over $\nu\in {\bf V}[\mathbb P]$ if
$\Vdash_{{\mathbb P}^+} \name{\eta} \not\leq^* \nu$.

\item[(d)] We say $\name{\eta}$ is $\Gamma_{\rm tow}$-big over $\nu \in {\bf V}[\mathbb P]$ if
$\Vdash_{{\mathbb P}^+} \name{\eta}^{-1}[\{1\}]
$ does not contain an infinite subset from ${\bf V}[\nu]$.
\end{myrules}
\end{myrules}
\end{definition}

\begin{definition}\label{2.4}
\begin{myrules}
\item[(1)]
 For any pseudo creature forcing $\bf q$ we write
\begin{equation}
\begin{split}
\Gamma_{\bf q} = \{({\mathbb P},{\mathbb P}^+,\name{\eta}) \in \Gamma_{\rm new}
\such &
\name{\eta} \in {\bf Q}^{{\bf V}[{\mathbb P}^+]} \wedge\\ 
& \Vdash_{{\mathbb P}^+} \mbox{`` for no } \bar{\gc} \in {\bf Q}^{{\bf V}[\mathbb P]}
\mbox{ we have } \name{\eta} \leq_{\bf q} \bar{\gc}\mbox{''}\}.
\end{split}
\end{equation}
We can allow ${\bf q}$ to be a ${\mathbb P}$-name.

\item[(2)]
We say $\name{\eta}$ is $\Gamma_{\bf q}$-big over $\nu \in {\bf V}[\mathbb P]$ if
for no $\bar{\gc} \in {\bf V}[\nu]$ $\Vdash_{{\mathbb P}^+} \name{\eta} \leq_{\bf q}
\bar{\gc}$.

\end{myrules}
\end{definition}

\begin{claim}\label{2.5}
The $\Gamma$'s from Definitions~\ref{2.3} and \ref{2.4} together with
the partial ordering $\lessdot_\Gamma$ are c.c.c.\ forcing bigness notions with amalgamation.
Moreover, if 
$(\forall \mu < \lambda) (\mu^{\aleph_0} < \lambda) $ 
then they have the  $(<\lambda)$-LS.
\end{claim}

\proof
(a) We need to check 
all the items in Definition~\ref{2.1} parts 1.\ and 3.\ and 4.
1.(a),(b),(c) follow directly from the definition of
$\Gamma_{\rm new}$.

Now for part 3: we first need to show that the union of 
c.c.c.\ forcings ${\mathbb P}_i$, $i<\delta$, is again c.c.c.
Since the chain is $\lessdot$ increasing, we can fix complete
projections for $i<j <\delta$ ${\rm pr}_{j,i}: {\mathbb P}_j \to {\mathbb P}_i$
(see \cite{abraham:handbook} for the characterization of $\lessdot$
via complete projections)
that form a commutative system.
Let $\{p_\alpha \such \alpha <\omega_1\} \subseteq \bigcup \{ {\mathbb P}_i \such i<\delta \}$.
We show that this is not an antichain. Suppose that $p_\alpha \in {\mathbb P}_{i(\alpha)}$
and such that for $\alpha <\beta$, $i(\alpha) \leq i(\beta)$.
Since ${\mathbb P}_0 $ has the c.c.c, there are
$\alpha < \beta$ such that
${\rm pr}_{i(\alpha),0}(p_\alpha) ||_0 \,  {\rm pr}_{i(\beta),0}(p_\beta)$.
But then $p_\alpha ||_\beta \, p_\beta$  and hence $p_\alpha ||_\delta \,  p_\beta$.

Next we have to show that ${\mathbb P_i} \lessdot 
\bigcup \{ {\mathbb P}_i \such i<\delta \}$.
This is again easy chasing arrows.

Now we need to show that $\bigcup \{ {\mathbb P}_i \such i<\delta \}
\lessdot \bigcup \{ {\mathbb P}_i^+ \such i<\delta \}$ and that
\begin{equation}\label{bigmatrix}
\xymatrix{
{\mathbb P}_1^+ \ar[r]^{\lessdot} & {\mathbb P}_2^+ \ar[r]^{\lessdot} & {\dots} 
\ar[r]^(0.3){\lessdot} &
\bigcup \{ {\mathbb P}_i^+ \such i<\delta \} \\
 {\mathbb P}_1 \ar[r]^{\lessdot} \ar[u]^{\lessdot} & {\mathbb P}_2 \ar[u]^{\lessdot}
\ar[r]^{\lessdot}
& {\dots} \ar[r]^(0.3)\lessdot & 
\bigcup \{ {\mathbb P}_i \such i<\delta \} \ar[u]^{\lessdot}
}
\end{equation}

together with its complete embeddings and their inverse
 projections ${\rm pr}_{j,i}\colon  {\mathbb P}_j \to
{\mathbb P}_i$, ${\rm pr}_{j^+,i^+}\colon {\mathbb P}^+_j \to {\mathbb P}^+_i$,
$ {\rm pr}_{i^+,i} \colon {\mathbb P}^+_i \to {\mathbb P}_i$ i exact
in all rectangles.

Let $\{p_\alpha \such \alpha <\omega_1\} \subseteq \bigcup \{ {\mathbb P}_i \such i<\delta \}$
be a maximal antichain. We show that it is also a maximal antichain in 
$\bigcup \{ {\mathbb P}_i^+ \such i<\delta \}$. Suppose $p \in {\mathbb P}_i^+$
is incompatible with all the $p_\alpha$'s. But then
$f_{i^+,i}(p) $ is incompatible in $\bigcup \{ {\mathbb P}_i \such i<\delta \}$
with the alledged maximal antichain.

Now we have to show that also the rectangles with limits in the right hand 
corners are exact: Let ${\mathbb P} = \bigcup_{\alpha<\delta} 
{\mathbb P}_\alpha $  and ${\mathbb P}^+ = \bigcup_{\alpha<\delta} 
{\mathbb P}^+_\alpha $  and let $i<\delta$, $p^+_i \in {\mathbb P}_i^+$.
Then  
\begin{equation*}
\begin{split}
&
\sup_{{\mathbb P}_i} \{ p_i \in {\mathbb P}_1 \such 
p_i \perp_{{\mathbb P}^+_i}  p_i^+ \}
\\
& =\sup_{{\mathbb P}_{\alpha}} \{ p_{\alpha} 
\in {\mathbb P}_{\alpha} \such p_{\alpha} 
\perp_{{\mathbb P}_{\alpha}^+ }  p_1^+ \}
\\
& =\sup_{{\mathbb P}_\delta} \{ p_\delta \in {\mathbb P}_\delta \such p_\delta 
\perp_{{\mathbb P}_\delta^+ }  p_1^+ \},
\end{split}
\end{equation*}

because the rectangle from $i$ to $\alpha$ is exact and because 
${\mathbb P}_\delta^+$ is the direct limit of the $P_\alpha$'s.
But are the {\em limit} triples $(\bigcup_{\alpha<\delta} 
{\mathbb P}_\alpha, \bigcup_{\alpha<\delta} 
{\mathbb P}^+_\alpha,\name{\eta})$ in $\Gamma$? 
Say we formulate it for  $\Gamma_{\rm new}$
and the proof of
this item of Definition~\ref{2.1} can be used also for
$\Gamma_{\rm Cohen}$, $\Gamma_{\rm ud}$, $\Gamma_{\rm tow}$, $\Gamma_{\bf q}$ (from Claim~\ref{2.5}).

The essential part is to show (for (a) for $\Gamma_{\rm new}$) that
$$
\Vdash_{\bigcup\{{\mathbb P}_i^+ \such i<\delta\}} 
\name{\eta} \not\in {\bf V}^{\bigcup\{{\mathbb P}_i \such i<\delta\}}.$$

For this we use ``Tools for your forcing construction'' \cite{Goldstern93} p.\ 358,
the part for the c.c.c.\ forcings and the finite support iteration
(taking unions of increasing notions of forcings is the same as finite
support iteration).
Again we can use common properties of $\Gamma_{\rm new}$, $\Gamma_{\rm Cohen}$,
$\Gamma_{\rm ud}$,
$\Gamma_{\rm tow}$, $\Gamma_{\bf q}$.

We write  
${\bigcup\{{\mathbb P}_i^+ \such i<\delta\}} = {\mathbb P}^+_\delta$.
${\bigcup\{{\mathbb P}_i \such i<\delta\}} = {\mathbb P}_\delta$.

Note that all the time, for  all our five $\Gamma$'s, for the limit 
property, we use the same fact:
$$\Vdash_{{\mathbb P}^+} \mbox{``}
 \name{\eta} \mbox{ is in $\Gamma$-relation with all the reals of } 
V^{\mathbb P}\mbox{''}$$ 
is of the form
$$\Vdash_{{\mathbb P}^+} (\forall \nu \in \omega^\omega \cap V^{\mathbb P}) 
\neg \name{\eta} R \nu,$$
where $R= R _\Gamma$ is an $F_\sigma$ relation.
There is a preservation theorem for
this situation, a modification of Theorem~8.4 of \cite{Goldstern93}.

\begin{hypothesis}\label{2.6}
We assume that $R = \sqsubset = \bigcup\{ \sqsubset_n \such n <\omega \}$ and
for all $n$ and $\eta \in {\mathbb R}$ the set
$\{ g \in {\bf V}^{\mathbb P} \such \eta \sqsubset_n g\}$ is closed.
\end{hypothesis}

This is true for
$R=R_\Gamma$ being the equality for item (a), being $\leq^*= \bigcup\{\{
(f,g) \such f\restriction [n,\infty)
\leq g \restriction[n,\infty)\} \such n<\omega\}$ for item (b), 
being $\sqsubset_{\rm Cohen}$ (see \cite[Chapter 7]{BJ}) for item (c) and
$\leq_{\rm tow} = \bigcup \leq_{{\rm tow},n}$
with $\eta \leq_{{\rm tow},n} g$ iff $\eta^{-1}\{1\} \cap g^{-1}\{1\} \cap [n,\infty) 
= \emptyset$.

\begin{theorem}\label{2.7} 
(Modification of Theorem~8.4 of \cite{Goldstern93}).
Let $\la {\mathbb P}^+_\alpha \such \alpha < \delta \r$ 
be a $\lessdot$ increasing continuous sequence of c.c.c.\ forcings and
let ${\mathbb P}_\delta$ be their direct limit. 
Let $\la {\mathbb P}_\alpha, \such \alpha < \delta \r$ 
be a sequence such that
all the rectangles from the diagramme~\eqref{bigmatrix}
are exact.
Let $\name{\eta}$ be a ${\mathbb P}_0^+$-name
such that 
$$\forall \alpha < \delta 
\Vdash_{{\mathbb P}_\alpha^+} \forall g \in {\bf V}^{{\mathbb P}_\alpha}
(\name{\eta} \not\sqsubset g).
$$
Then $$
\Vdash_{{\mathbb P}_\delta^+} \forall g \in {\bf V}^{{\mathbb P}_\delta}
(\name{\eta} \not\sqsubset g).
$$
\end{theorem}

\proof
Let $\la \delta_i \such i < \cf(\delta) \r$ be sequence of ordinals converging 
to $\delta$. We assume that $\cf(\delta) = \omega$, because otherwise
there is nothing to show because by the c.c.c.\ each real will 
appear for the first time  at a stage with cofinality less of equal $\omega$.
Assume the conclusion is false, so there is a condition $p_0 \in {\mathbb P}^+_\delta$ and a 
${\mathbb P}_\delta$-name $\name{g}$ such that
$p_0 \Vdash_{{\mathbb P}_\delta^+} \name{\eta} \sqsubset \name{g}$.
We find a condition $p_1\leq p_0$ and an integer $n$ such that
$p_1 \Vdash \name{\eta} \sqsubset_n \name{g}$.
Since $p_1$ has finite support there is some $i<\cf(\delta)$ such that
$p_1 \in {\mathbb P}_{\delta_i}^+$.

Let $G$ be ${\mathbb P}^+_\delta$ -generic over ${\bf V}$, and
let $G_{\delta_i}$ be $G \cap {\mathbb P}_{\delta_i}$.
Clearly in ${\bf V}[G_{\delta_i}]$ with $p_1 \in G_{\delta_i}$ we
have that
\begin{equation}\label{positivfuern}
\Vdash_{\delta_i,\delta} \name{\eta} \sqsubset_n \name{g}.
\end{equation}

On the other hand $\name{\eta}[G_{\delta_i}] \not\sqsubset_n x$ for
all $x \in V[G_{\delta_i}]$.
Let  ${\rm tree}(\name{g}[G_{\delta_i}])$ be the
tree of possible  evaluations of $\name{g}[G_{\delta_i}]$, all possible 
evaluations, 
not only generic ones. 
Take $x \in ({\rm tree}(\name{g}[G_{\delta_i}]))^{{\bf V}[G_{\delta_i}]}$.

We consider the set $\{ h \such \name{\eta} \not\sqsubset_n h\}\subseteq {\bf V}[G_\delta]$. 
This is an open set containing $x$. 
So there is some $k<\omega$ such that
$$\Vdash_{{\mathbb P}_{\delta_i,\delta}}
\forall h \, (h\restriction k = \name{g}[G_{\delta_i}]
\restriction k \rightarrow
\name{\eta}[G_{\delta_i}] \not\sqsubset_n h).$$
Let $q \in P_\delta/G_{\delta_i}$ be a condition forcing 
$\name{g} \restriction k =
x\restriction k$.
Then also $q$
 forces in ${\mathbb P}_{\delta_i,\delta}$ that  $\name{\eta}
\not\sqsubset_n \name{g}$ in contradiction to \eqref{positivfuern}.
\proofendof{limit\, part\, of\, \ref{2.5}}

\smallskip

Now for the {\em amalgamation}:
Let $({\mathbb P}_1,{\mathbb P}_1^+,\name{\eta}) \in \Gamma$,
let ${\mathbb P}_1 \lessdot {\mathbb P}_2$, and let
${\mathbb P}_2 / {\mathbb P}_1$ have the Knaster property.
Then choose $f \colon {\mathbb P}_2 \to {\mathbb P}_2'$ such
that $f \restriction {\mathbb P}_1 = id$ and ${\mathbb P}_2' \cap {\mathbb P}_1^+ ={\mathbb P}_1$.
Since ${\mathbb P}_2/ {\mathbb P}_1$ has the Knaster property 
${\mathbb P}_2^+ =
{\mathbb P}_1^+
\times_{{\mathbb P}_1} {\mathbb P}_2=
{\mathbb P}_1 * ({\mathbb P}_1^+/{{\mathbb P}_1} \times
{\mathbb P}_2/{{\mathbb P}_1})$ has the c.c.c.

Then ${\mathbb P}_1^+ \lessdot {\mathbb P}_2^+$.
and $id \colon {\mathbb P}_2 \to {\mathbb P}_2^+$ is a complete embedding over ${\mathbb P}_1$.

\begin{theorem}\label{2.8} 
Let $({\mathbb P}_1,{\mathbb P}_1^+,\name{\eta}) \in \Gamma$, $R_\Gamma = \sqsubset$,
let ${\mathbb P}_1 \lessdot {\mathbb P}_2$, and let
${\mathbb P}_2 / {\mathbb P}_1$ have the Knaster property.
Let $\name{\eta}$ be a ${\mathbb P}_1^+$-name
such that 
$$ \Vdash_{{\mathbb P}_1^+} \forall g \in {\bf V}^{{\mathbb P}_1}
(\name{\eta} \not\sqsubset g).
$$
Then for ${\mathbb P}_2^+ = {\mathbb P}_1^+ \times_{{\mathbb P}_1} {\mathbb P}_2$
we have that 
$$\Vdash_{{\mathbb P}_2^+} \forall g \in {\bf V}^{{\mathbb P}_2}
(\name{\eta} \not\sqsubset g).
$$
\end{theorem}

\proof
Let $({\mathbb P}_1,{\mathbb P}_1^+,\name{\eta}) \in \Gamma$,
let ${\mathbb P}_1 \lessdot {\mathbb P}_2$, and let
${\mathbb P}_2 / {\mathbb P}_1$ have the Knaster property.
Then choose an isomorphism $f \colon {\mathbb P}_2 \to {\mathbb P}_2'$ such
that $f \restriction {\mathbb P}_1 = id$ and ${\mathbb P}_2'
 \cap {\mathbb P}_1^+ ={\mathbb P}_1$.
Since ${\mathbb P}_2/ {\mathbb P}_1$ has the Knaster property ${\mathbb P}_2^+
=
{\mathbb P}_1^+
\times_{{\mathbb P}_1} {\mathbb P}'_2$ has the c.c.c.

Then ${\mathbb P}_1^+ \lessdot {\mathbb P}_2^+$
and $f \colon {\mathbb P}_2 \to {\mathbb P}_2^+$ is a complete embedding over ${\mathbb P}_1$.
Then $\Vdash_{{\mathbb P}_2^+} (\forall g
\in {\bf V}^{{\mathbb P}_2} \bigwedge_{n<\omega} \eta \not\sqsubset_n g)$
since otherwise there would be some $n <\omega$ 
and  $g\in {\bf V}^{{\mathbb P}_2}$ and $p_2^+$ such that
$p_2^+ \Vdash_{{\mathbb P}_2^+} \eta \sqsubset_n g$.
Let $G_2^+$ be ${\mathbb P}_2^+$-generic over $V$ and let $p_2^+ \in G_2^+$
and $G_2 = G_2^+\cap {\mathbb P}_2$.
Then $T= \{ s \in \omega^{<\omega} \such (\exists g \in V^{{\mathbb P}_2})
(s \subseteq g \, \wedge \, 
\eta[G_2^+] \sqsubset_n g) \}$ is an ill-founded tree in $V[G_2^+]$.
But since $\sqsubset_n$ is absolute and since $\eta[G_2^+] =
\eta[G_2^+\cap {\mathbb P}_1^+]$ we have that the
$T$ is ill-founded is $V[G_2^+\cap {\mathbb P}_1^+]$ as well
and $G_2^+ \cap {\mathbb P}_1^+ $ is ${\mathbb P}_1^+$-generic over ${\bf V}$, 
and this is a
contradiction.
\proofendof{amalgamation \, part \, of \, \ref{2.5}}

\medskip

Now for the LS property: Since it has the c.c.c., 
$\name{\eta}$ is determined by
countably many countable antichains.

$|Y| = \mu < \lambda$, and hence by closing ${\mathbb P} \cap Y$ under the 
suprema in  ${\mathbb P}$
of all subsets of  ${\mathbb P}$, we get ${\mathbb P}_1 \lessdot {\mathbb P}$
such that ${\mathbb P} \cap Y  \subseteq {\mathbb P}_1$.
We can require that $|{\mathbb P}_1| < \lambda$, because 
every supremum is determined by countably many elements and hence
there are less or equal than $\mu^{\aleph_0}< \lambda$ many tasks.
We do the same for $Y$ in ${\mathbb P}^+$ and for ${\mathbb P}_1$ 
and thus find ${\mathbb P}_1^+$.
Then we automatically have ${\mathbb P}_1 \lessdot {\mathbb P}_1^+$ and
${\mathbb P} \times_{{\mathbb P}_1}
{\mathbb P}_1^+\lessdot {\mathbb P}^+$, and thus get
the diagramme

\vskip 5mm
\hskip 2cm
\xymatrix{
{\mathbb P} \ar[r]^{\lessdot} & {\mathbb P}^+ \\
 {\mathbb P}_1 \ar[u]^{\lessdot} 
\ar[r]^{\lessdot} & {\mathbb P}_1^+ \ar[u]^{\lessdot}\\
 {\mathbb P}\cap Y \ar[r]^{\subseteq} \ar[u]^{\subseteq} & Y \ar[u]^{\subseteq}
}
\vskip 5mm

whose upper square is exact.
This was the proof for $\Gamma_{\rm new}$ and all other $\Gamma$'s, because it
did not use any specific property but Hypothesis~\ref{2.6} \proofendof{\ref{2.5}}

\smallskip

\begin{definition}\label{2.9}
Let $\bf q$ be a {\bf forgetful} pcrf, and let $\Gamma$ be a forcing bigness notion.
\begin{myrules}
\item[1.]
We say ``$\bf q$ is orthogonal to $\Gamma$'' and write ${\bf q} \perp \Gamma$, if
the following holds: If 
\begin{myrules}

\item[(a)] $({\mathbb P}_1,{\mathbb P}_1^+, \name{\eta}) \in \Gamma$ and

\item[(b)] ${\mathbb P}_1 \lessdot {\mathbb P}_2$ and ${\mathbb P}_2 /{\mathbb P}_1$ is Cohen$_A$
for some set $A$ and ${\mathbb P}_2 \cap {\mathbb P}^+_1 = {\mathbb P}_1$,

\item[(c)] $\name{{\mathcal C}_2}$ is a ${\mathbb P}_2$-name and ${\mathbb P}_2$
forces that it is $\bf q$-nice,
\end{myrules}
then we can find $({\mathbb P}_2^+,f, \name{{\mathcal C}^+_2} )$ such that

\begin{myrules}

\item[($\alpha$)] ${\mathbb P}_1^+ \lessdot {\mathbb P}_2^+$,

\item[($\beta$)] $f$ is a complete embedding of ${\mathbb P}_2$ into ${\mathbb P}_2^+$ over
${\mathbb P}_1$,

\item[($\gamma$)] $\Vdash_{{\mathbb P}_2^+} \mbox{``} f({\name{\mathcal C}_2}) \subseteq
\name{{\mathcal C}_2}^+ \wedge \name{{\mathcal C}_2}^+ \mbox{ is } {\bf q}\mbox{-nice''}$,

\item[($\delta$)] $({\mathbb P}_1,{\mathbb P}_1^+,\name{\eta}) \leq_{\Gamma}
 (f({\mathbb P}_2),{\mathbb P}_2^+,\name{\eta}) \leq_{\Gamma}
 (f({\mathbb P}_2) *{\bbforcing}(f(\name{{\mathcal C}_2})),{\mathbb P}_2^+
* {\bbforcing}(\name{{\mathcal C}_2}^+),\name{\eta})$.
\end{myrules}

We say that ${\bf q} \perp \Gamma$ with a simple witness, if 
in the conclusion we can choose 
${\mathbb P}_2^+ = ({\mathbb P}_1^+ \times_{{\mathbb P}_1} {\mathbb P}_2) 
\times {\rm Cohen}_\lambda$,
$\lambda = (|{\mathbb P}_1^+| + |{\mathbb P}_2|)^{\aleph_0}$, $f = id$.

\medskip

\item[2.] We say that $({\mathbb P}, {\mathbb P}^+,\name{\eta}, \name{\mathcal C},
\name{{\mathcal C}^+})$ is a $({\bf q}, \Gamma)$-tuple if
\begin{myrules}

\item[(a)] $({\mathbb P},{\mathbb P}^+,\name{\eta}) \in \Gamma$,

\item[(b)] $\name{{\mathcal C}}$ is a ${\mathbb P}$-name and ${\mathbb P}$
forces that it is $\bf q$-nice,

\item[(c)] $\name{{\mathcal C}^+}$ is a ${\mathbb P}^+$-name and ${\mathbb P}^+$
forces that it is $\bf q$-nice,

\item[(d)] $\name{{\mathcal C}} = \name{{\mathcal C}^+} \cap {\bf Q}^{{\bf V}[{\mathbb P}]}$

\item[(e)] $({\mathbb P},{\mathbb P}^+,\name{\eta}) \leq_{\Gamma}
 ({\mathbb P} * {\bbforcing}(\name{{\mathcal C}}),{\mathbb P}^+
* {\bbforcing}(\name{{\mathcal C}}^+),\name{\eta})$.
\end{myrules}
\medskip

\item[3.] We say ``$\bf q$ is strongly orthogonal
to $\Gamma$'' and write ${\bf q} \perp_{st} \Gamma$, if: If
\begin{myrules}

\item[(a)] $\Gamma$ is a forcing bigness notion,

\item[(b)] $\bf q$ is a forgetful pcrf,

\item[(c)] $({\mathbb P}_1, {\mathbb P}_1^+,\name{\eta}, \name{{\mathcal C}_1},
\name{{\mathcal C}_1^+})$ is a $({\bf q}, \Gamma)$-tuple and if
${\mathbb P}_1 \lessdot {\mathbb P}_2$ and 
${\mathbb P}_2 /{\mathbb P}_1$ satisfies the Knaster condition and 
$\Vdash_{{\mathbb P}_2} \name{{\mathcal C}_2}$ is a $\bf q$-nice family extending 
$\name{{\mathcal C}_1}$,
\end{myrules}
then we can find $({\mathbb P}_2^+,f, \name{{\mathcal C}^+_2} )$ such that
\begin{myrules}

\item[($\alpha$)] ${\mathbb P}_1^+ \lessdot {\mathbb P}_2^+$,

\item[($\beta$)] $f$ is a complete embedding of ${\mathbb P}_2$ into ${\mathbb P}_2^+$ over
${\mathbb P}_1$,

\item[($\gamma$)] $\Vdash_{{\mathbb P}_2^+} \mbox{``} f({\name{\mathcal C}_2}) \subseteq
\name{{\mathcal C}^+_2}\mbox{''}$,

\item[($\delta$)] $({\mathbb P}_2,{\mathbb P}_2^+,\name{\eta}, 
f(\name{{\mathcal C}_2}),\name{{\mathcal C}_2}^+)$ is a $({\bf q},\Gamma)$-tuple,

\item[($\eps$)] $({\mathbb P}_1 * {\bbforcing}(\name{{\mathcal C}}_1),{\mathbb P}_1^+
* {\bbforcing}(\name{{\mathcal C}}_1^+),\name{\eta})
\leq_\Gamma 
 ({\mathbb P}_2 * {\bbforcing}(f(\name{{\mathcal C}}_2)),{\mathbb P}_2^+
* {\bbforcing}(\name{{\mathcal C}}_2^+),\name{\eta})$.
\end{myrules}

\medskip

\item[4.] The order $\leq$ on the family of $({\bf q},\Gamma)$-tuples is defined by
$$
({\mathbb P}_1,{\mathbb P}_1^+,\name{\eta}, 
\name{{\mathcal C}_1},\name{{\mathcal C}_1^+}) \leq 
({\mathbb P}_2,{\mathbb P}_2^+,\name{\eta}, 
\name{{\mathcal C}_2},\name{{\mathcal C}_2^+})
$$
if both are $({\bf q},\Gamma)$-tuples and
\begin{myrules}

\item[($\alpha$)] ${\mathbb P}_1  \lessdot {\mathbb P}_2 $,

\item[($\beta$)] ${\mathbb P}_1^+ \lessdot {\mathbb P}_2^+$,

\item[($\gamma$)] $\name{\eta_1} = \name{\eta_2}$,

\item[($\delta$)] $\Vdash_{{\mathbb P}_2} \mbox{``} {\name{\mathcal C}_1} \subseteq
\name{{\mathcal C}_2}\mbox{''}$,

\item[($\eps$)] $\Vdash_{{\mathbb P}_2^+} \mbox{``} \name{{\mathcal C}_1^+} \subseteq
\name{{\mathcal C}_2^+}\mbox{''}$.

\item[($\zeta$)] $({\mathbb P}_1 * {\bbforcing}(\name{{\mathcal C}}_1),{\mathbb P}_1^+
* {\bbforcing}(\name{{\mathcal C}}_1^+),\name{\eta})
\leq_\Gamma 
 ({\mathbb P}_2 * {\bbforcing}(f(\name{{\mathcal C}}_2)),{\mathbb P}_2^+
* {\bbforcing}(\name{{\mathcal C}}_2^+),\name{\eta})$.

\end{myrules}
\end{myrules}

\end{definition}

\begin{claim}\label{2.10}
Assume that $\Gamma$ is a c.c.c.\ forcing bigness notion with amalgamation.
A sufficient condition for ${\bf q} \perp \Gamma$ is the following:

$\boxtimes$ If
\begin{myrules}
\item[(a)] $({\mathbb P},{\mathbb P}^+,\name{\eta}) \in \Gamma$,
\item[(b)] $\Vdash_{\mathbb P} \mbox{``} \name{\bar{\gc}} \in {\bf Q}$''
or just for some $({\mathbb P}_0,\name{\mathcal C})$ we have that
${\mathbb P} \lessdot {\mathbb P}_0$ and
$\Vdash_{{\mathbb P}_0} \mbox{``}\name{\mathcal C} $ is $\bf q$-nice'' and
$\Vdash_{{\mathbb P}_0} \mbox{``}\name{\bar{\gc}} $ is generic for 
$({\bf Q}(\name{\mathcal C}))^{{\bf V}[G_{\mathbb P}]}$'',
\end{myrules}
then we can find 
$\name{{\mathcal C}^+}$ and a 
${\mathbb P}^+$-name of a $\bf q$-nice family to which $\bar{\gc}$
belongs such that
if ${\mathcal J}$ is a ${\mathbb P}$-name of a subset of ${\bbforcing}_{\mathcal C}(\bar{\gc})$ 
which is explicitly predense then
$\Vdash_{{\mathbb P}^+}\mbox{``} {\mathcal J}$ 
is a predense subset of ${\bbforcing}({\mathcal C}^+)$''
and $({\mathbb P} \ast \bbforcing(\name{\mathcal C}), {\mathbb P}^+ \ast
\bbforcing(\name{\mathcal C}^+),\name{\eta}) \in \Gamma$.

\end{claim}

\proof
Suppose ${\mathbb P}_1$, 
${\mathbb P}_1^+$, $\name{\eta}$, ${\mathbb P}_2$ and 
$\name{{\mathcal C}_2}$ are given as in the premises 
(a), (b) and (c) of Definition~\ref{2.9} Part 1.
We have to find $({\mathbb P}_2^+, f, \name{{\mathcal C}_2^+})$ as there.

Since $\Gamma$ has the amalgamation property, there is 
$({\mathbb P}_1,{\mathbb P}_1^+,\name{\eta}) \leq_\Gamma (f({\mathbb P}_2),
{\mathbb P}_2^+, \name{\eta})$. 
Now we read the condition given in 
Claim~\ref{2.10} for
$(f({\mathbb P}_2), {\mathbb P}_2^+, \name{\eta})$ and 
$\name{{\mathcal C}_2}$ and find $\name{{\mathcal C}_2^+}$ 
such that 
 $({\mathbb P}_1,{\mathbb P}_1^+,\name{\eta}) \leq_{\Gamma}
 (f({\mathbb P}_2),{\mathbb P}_2^+,\name{\eta}) \leq_{\Gamma}
 (f({\mathbb P}_2) *{\bbforcing}(f(\name{{\mathcal C}_2})),{\mathbb P}_2^+
* {\bbforcing}(\name{{\mathcal C}_2^+}),\name{\eta})$.
\proofend

\section{Long low c.c.c.\ iteration and ${\bf q} \perp \Gamma$}\label{S3}

\begin{hypothesis}\label{3.1}
$\lambda = \cf(\lambda)$ satifies
$$(\forall \alpha < \lambda)(|\alpha|^{\aleph_0} < \lambda),$$
and $\zeta(*)$ is a limit of uncountable cofinality.
\end{hypothesis}

\begin{definition}
\label{3.2}
Let $\bar{\Gamma} = \la \Gamma_\zeta \such \zeta < \zeta(*)\r$ be sequence of c.c.c.\ forcing 
bigness notions.
\begin{myrules}
\item[(0)] $\bar{\Gamma}$ is nice if each $\Gamma_\zeta$ is, i.e.
$\Gamma_\zeta$ is a c.c.c.\ forcing bigness notion, has amalgamation and the $(<\lambda)$-L.S.\/
$\bar{\Gamma}$ has free amalgamation if each $\Gamma_\zeta$ has.

\item[(1)] $K=K_{\lambda,\zeta(*)}(\bar{\Gamma})$ is the family of $\itername=
\la {\mathbb P}_\zeta,\name{\eta}_\zeta\such \zeta < \zeta(*)\r =
\la {\mathbb P}^\itername_\zeta,\name{\eta}^\itername_\zeta\such \zeta 
< \zeta(*)\r$
such that 

\begin{myrules}

\item[(a)] ${\mathbb P}_\zeta$ is a c.c.c.\ forcing notion of cardinality $<\lambda$ for $\zeta
\leq \zeta(*)$.
\item[(b)] $\la {\mathbb P}_\zeta \such \zeta < \zeta(*)\r$ is $\lessdot$-increasing,
and if $\xi <\zeta(*)$ is a limit ordinal, then ${\mathbb P}_{<\xi} =
\bigcup\{ {\mathbb P}_\zeta \such \zeta < \xi \} \lessdot {\mathbb P}_\xi$, if
$\aleph_0 < \cf(\xi)$.
\item[(c)] 
$\name{\eta}_\zeta$ is a ${\mathbb P}_{\zeta+1}$ name such that
$({\mathbb P}_\zeta,{\mathbb P}_{\zeta +1},\name{\eta}_\zeta)$ is $\Gamma_\zeta$-big.
\end{myrules}

We stipulate ${\mathbb P}^\itername_{\zeta(*)}=
\bigcup \{{\mathbb P}^\itername_{\zeta} \such \zeta 
< \zeta(*)\}$
and we may write ${\mathbb P}_{\zeta}[\itername]$ instead of
${\mathbb P}^\itername_{\zeta}$.

\item[(1a)] We say $\itername$ is $\lambda$-standard if ${\itername} 
\in {\mathcal H}(\lambda)$ and we write
$K^* =K^*_{\lambda,\zeta(*)}(\bar{\Gamma}):= K_{\lambda, \zeta(*)}(\bar{\Gamma}) \cap
{\mathcal H}(\lambda)$ when $\bar{\Gamma}$, $\lambda$, $\zeta(*)$
are clear from the context.

\item[(2)]
We define three binary relations $\leq$, $\leq_{\pr}$ and $\leq_{ap}$ on $K^*$:
\begin{myrules}
\item[(a)]
$\iterp \leq_{pr} \itername$ iff $\zeta < \zeta(*)$ implies $\name{\eta}^\iterp_\zeta
=\name{\eta}^\itername_\zeta$ and 
${\mathbb P}^\iterp_\zeta \lessdot {\mathbb P}^\itername_\zeta$,

\item[(b)]
$\iterp \leq \itername$ iff there is some $\xi < \zeta(*)$ such
that  $\zeta \in [\xi, \zeta(*))$ the relations $\name{\eta}^\iterp_\zeta
=\name{\eta}^\itername_\zeta$ and 
${\mathbb P}^\iterp_\zeta \lessdot {\mathbb P}^\itername_\zeta$ hold,

\item[(c)]
$\iterp \leq_{ap} \itername$  iff 
$\iterp \leq \itername$ and $(\forall^* \zeta < \zeta(*))(\name{\eta}^\iterp_\zeta
=\name{\eta}^\itername_\zeta \wedge
{\mathbb P}^\iterp_\zeta = {\mathbb P}^\itername_\zeta)$,

\item[(d)] $\iterp \equiv \itername$ iff $\iterp \leq_{ap} \itername \leq_{ap} \iterp$.
\end{myrules}

\item[(3)] If $\la \itername_\alpha \such \alpha < \delta \r$ is 
$\leq_{pr}$-increasing, let $\itername_\delta=
\bigcup_{\alpha< \delta} \itername_\alpha$ be
defined as
$\la \bigcup_{\alpha<\delta} {\mathbb P}^{\itername_\alpha}_\zeta, \name{\eta}_\zeta
\such \zeta<\zeta(*)\r$.

\item[(4)] ${\mathfrak K}^*$ as a forcing notion is $(K^*,\leq)$ and the generic is $\name{\mathbb P} =
\bigcup\{{\mathbb P}^\itername_{\zeta(*)}\such \itername \in \name{G}_{K^*}\}$.

\item[(5)] If for all $\zeta < \zeta(*)$, $\Gamma_\zeta = \Gamma$, then we may write
$\Gamma$ instead of $\bar{\Gamma}$.

\item[(6)] We say $\la {\itername}_\alpha \such \alpha < \alpha^*\r$ is $\leq_{pr}$-continuous if
\begin{myrules}
\item[(a)] $\alpha <\beta < \alpha^* \Rightarrow \itername_{\alpha} \leq_{pr} 
\itername_\beta$.
\item[(b)] If  $\beta < \alpha^*$ is a limit ordinal
then $\itername_\beta =
\bigcup_{\alpha < \beta} \itername_\alpha$, i.e.,
${\mathbb P}_\zeta^{\itername_\beta}=
\bigcup\{
{\mathbb P}_\zeta^{\itername_\alpha} \such \alpha<\beta\}$
for every $\zeta<\zeta(*)$.
\end{myrules}
\item[(7)]
We replace continuous by ``weakly continuous'' if in clause (b)
above only
$\bigcup\{
{\mathbb P}_\zeta^{\itername_\alpha} \such \alpha<\beta\}
\lessdot {\mathbb P}_\zeta^{\itername_\beta}$
is demanded.
\end{myrules}
\end{definition}

\begin{center}
\setlength{\unitlength}{0.6cm}
\begin{picture}(21,10)
\put(1,1.5){\vector(1,0){19}}

\put(19.5,0.8){$\lambda$}
\put(2.8,0.8){$\iterp$}
\put(3.8,0.8){$\iterq$}
\put(1,2.5){\vector(1,0){19}}
\put(1,3.5){\vector(1,0){19}}
\put(1,4.5){\vector(1,0){19}}
\put(1,5.5){\vector(1,0){19}}
\put(1,6.5){\vector(1,0){19}}

\put(-0.2,9){$\zeta(*)$}

\put(1,1.5){\vector(0,1){8}}
\put(2,1.5){\vector(0,1){8}}
\put(3,1.5){\vector(0,1){8}}
\put(4,1.5){\vector(0,1){8}}
\put(4,9){$\la {\mathbb P}^\iterq_\zeta, \name{\eta}^\iterq_\zeta 
\such \zeta < \zeta(*)\r$}

\put(20,6.3){$\Gamma_\zeta$}

\put(4,8){$({\mathbb P}_\zeta, {\mathbb P}_{\zeta+1}, \name{\eta}) 
\in \Gamma_\zeta$}

\put(3,0){Figure 1: A sketch of $K_{\lambda,\zeta(*)}(\bar{\Gamma})$}
\end{picture}
\end{center}

\medskip

\begin{claim}\label{3.3}
\begin{myrules}
\item[(1)]
$\leq$ and $\leq_{pr}$ are partial orders on $K^*$.

\item[(2)] $\equiv$ is an equivalence relation on $K^*$ and
the forcing $(K^*,\leq)$ is equivalent to  $(K^* / \equiv, \leq)$.

\item[(3)] If $\Gamma$ is nice and if
$\delta<\lambda$ is a limit ordinal 
and $\la \itername_\alpha \such \alpha < \delta \r$ is $\leq_{pr}$-increasing and continuous
as a sequence of members in $K^*$ and $\itername_\delta
=\bigcup_{\alpha<\delta} \itername_\alpha$ then $\itername_\delta \in K^*$ and 
$\la \itername_\alpha \such \alpha \leq \delta \r$ is $\leq_{pr}$-increasing and continuous.
Instead of continuous we can demand weakly continuous.

\item[(4)]
$\iterp \equiv \itername \Rightarrow \iterp \leq \itername$ and $\iterp \leq_{pr} 
\itername \Rightarrow\iterp \leq \itername$.

\item[(5)] As each $\Gamma_\zeta$ has the 
$(<\lambda)$-L.S., clearly $K^*_{\lambda, \zeta(*)} \neq \emptyset$.

\item[(6)] If each $\Gamma_\zeta$, $\zeta < \zeta(*)$, has free amalgamation and $\itername
\in K^*_{\lambda,\zeta(*)} $ and $\xi < \zeta(*)$ and ${\mathbb P}^\itername_\xi
\lessdot {\mathbb P}
\in {\mathcal H}(\lambda)$ and
${\mathbb P}$ is c.c.c.\ and $\Vdash_{{\mathbb P}^\itername_{\zeta(*)}}
\mbox{``}{\mathbb P}/(\name{G} \cap {\mathbb P}^\itername_\xi)$ satisfies the c.c.c.''.
\em{then} we can find some $(\itername^+,f)$ such that
\begin{myrules}
\item[(a)] $\itername \leq_{pr} \itername^+ \in K^*$,

\item[(b)] $\zeta < \xi \Rightarrow {\mathbb P}^{\itername^+}_\zeta 
= {\mathbb P}_\zeta^{\itername}$

\item[(c)] $f$ is a $\lessdot$-complete embedding of ${\mathbb P}$ into 
${\mathbb P}^{\itername^+}_\xi$
over ${\mathbb P}^\itername_\xi$:
$$
\xymatrix{
{\mathbb P}_{\xi+1}^{\iterq } \ar[rr]^{\lessdot} && {\mathbb P}_{\xi+1}^{\iterq^+}\\
 {\mathbb P}_\xi^\iterq \ar[u]^\lessdot \ar[r]^\lessdot & {\mathbb P}  \ar[r]^{\lessdot, f} 
 & {\mathbb P}_\xi^{\iterq^+} \ar[u]^\lessdot
}
$$
\end{myrules}
\end{myrules}
\end{claim}
\proof
This follows almost immediately from the definitions. \proofend

\begin{claim}\label{3.4}
\begin{myrules}
\item[(1)] If $\iterp \leq \itername$ then for some $\itername' \equiv \itername$ we have
$\iterp \leq_{pr} \itername'$.

\item[(2)]
If $\name{G}_{K^*} \subseteq K^*_{\lambda,\zeta(*)} $ is generic over ${\bf V}$, then
$\name{\mathbb P} =
\bigcup\{
{\mathbb P}^\iterp_{\zeta(*)} \such \iterp \in \name{G}_{K^*}\}$ is
a c.c.c.\ forcing notion of cardinality $\leq \lambda$ and
$\itername \Vdash {\mathbb P}^\itername_{\zeta(*)} \lessdot 
\name{\mathbb P}$, but it locally looks like
${\mathbb P}^\itername_{\zeta(*)}$.
\item[(3)] If $\bar{\Gamma}$ is nice and if $\iterp \in K^*$ and if 
$\xi < \zeta(*)$ and if we define $\itername$ by 
$${\mathbb P }_\zeta^\itername = \left\{
\begin{array}{ll}
{\mathbb P}_\zeta^\iterp & \mbox{ if } \zeta < \xi,\\
{\mathbb P}^\iterp_\zeta * (^{\omega >} \omega, \triangleleft) & \mbox{ if }
\zeta \in [\xi, \zeta(*)),
\end{array}
\right.
$$ then $\iterp \leq_{pr} \itername \in K^*$. 

\end{myrules}
\end{claim}

\begin{claim}\label{3.5}
Let $G_{{\mathfrak K}^*} \subseteq {\mathfrak K}^*$ be generic over ${\bf V}$ and let ${\mathbb P}
= \name{\mathbb P}[G_{{\mathfrak K}^*}]$. Then in ${\bf V}[G_{\mathfrak K}]$:
\begin{myrules}
\item[(a)] ${\mathbb P}$ is a c.c.c.\ forcing notion of cardinality 
$\lambda$ adding 
$\lambda$ reals.
\item[(b)] If $\name{\nu}$ is a ${\mathbb P}$-name for a real
and $\iterq^* = \langle {\mathbb P}^*_\zeta, \name{\eta}_\zeta \such \zeta < \zeta(*) \r
\in G_{{\mathfrak K}^*}$, then for
every sufficiently large $\zeta < \zeta(*)$ 
for some ${\mathbb P}_\zeta \lessdot {\mathbb P}_\zeta^+ \lessdot {\mathbb P}$ we have
\begin{myrules} 
\item[($\alpha$)] $\name{\nu}$ is a ${\mathbb P}_\zeta$-name,
\item[($\beta$)] $({\mathbb P}_\zeta,{\mathbb P}_\zeta^+,\name{\eta}_\zeta) \in \Gamma_\zeta$.
\end{myrules}
\item[(b')] So, if $\zeta < \zeta(*) \Rightarrow 
\Gamma_\zeta = \Gamma_{ud}$ and $\iterq^*$ is as above, then
$\iterq^* \Vdash_{\mathbb P} \mbox{``} \{ \name{\eta}_\zeta \such \zeta < \zeta(*) \} $ is unbounded''
(hence $\gb \leq \zeta(*)$).
\item[(c)] Similar for other relations on ${}^\omega 2$.
\end{myrules}
\end{claim}
\proof
Item (b): 
Since  ${\name{\mathbb P}}$ has the c.c.c., for $m <\omega$, let 
$\{ p_{m,n} \such n <\omega \} \subseteq {\mathbb P}$ be a maximal
antichain forcing the value of $\name{\nu} \restriction m$.
${\mathfrak K}^*$ is $\lambda$-complete and $\lambda > \aleph_1$,
hence there is some $\iterq \in G_{K}$ such that 
$\{ p_{m,n} \such n,m \in \omega \} \subseteq {\mathbb P}^\iterq$
and $\iterq \geq \iterq^*$.

Since ${\mathbb P}^\iterq $ is the union of the increasing chain
$\la {\mathbb P}_\zeta^\iterq
\such  \zeta < \zeta(*) \r$ and since  $\cf(\zeta(*)) > \aleph_0$,
 for some $\zeta_0 < \zeta(*)$
we have that 
$\{ p_{m,n} \such n,m\in \omega \} \subseteq \name{\mathbb P}_{\zeta_0}^\iterq$.

As $\iterq^* \in G$, $\iterq^* \leq \iterq$ we
have for some $\zeta_1 < \zeta(*)$ that
$\bigwedge_{\xi \in [\zeta_1,\zeta(*))}({\mathbb P}_\xi^{\iterq^*} \lessdot
{\mathbb P}_\xi^{\iterq}
\wedge \name{\eta}_\xi^{\iterq^*}= \name{\eta}_\xi^\iterq).$

So let $\zeta \in [\max(\zeta_0,\zeta_1),\zeta(*))$.
Since $\iterq \in G_{K^*}$, we have
 ${\mathbb P}_{\zeta}^\iterq \lessdot {\mathbb P}_{\zeta+1}^\iterq
\lessdot \name{{\mathbb P}}[G_{K^*}]$.
Now since 
$\name{\nu}$ is a ${\mathbb P}_\zeta^\iterq$-name, by the demands on $\iterq$ all follows.
\proofend

\begin{definition}\label{3.6} 
Let $\iterq_* \in {\mathfrak K}^* $, $ \zeta_* <\zeta(*)$
and $\zeta_* + \zeta(*) = \zeta(*)$, for notational simplicity,
and $p_* \in {\mathbb P}_\zeta^{\iterq_*}$.

Define $K^*_{\geq \iterq_*}= \{ \iterq \in K \such \iterq_* \leq \iterq\}$ and define
the function $F$ with domain $K^*_{\geq \iterq_*}
$ as follows: if $\iterq \in K^*_{\iterq_*} $ then $F(\iterq) = 
\la {\mathbb P}_\zeta' \such \zeta < \zeta(*)\r$ where 
${\mathbb P}_\zeta'
= 
{\mathbb P}^*[\iterq_{\zeta_* + \zeta}]\restriction
\{ p \in {\mathbb P}^*[\iterq_{\zeta_* + \zeta}] \such
{\mathbb P}^*[\iterq_{\zeta_* + \zeta}] \models p_* \leq p\}$ and
lastly $K^*_{\iterq_*,\zeta_*,p_*} = \rge(F)$.
\end{definition}

\begin{claim}\label{3.7}
1. $K^*_{\geq \iterq_*}, K^*_{\iterq_*,\zeta_*,p_*} \subseteq K^*$.

2. $F$ preserves the order $\leq, \leq_{pr}, \leq_{ap}$ and their negations and the forcing 
$K^*_0$ above $(\iterq_*,p_*)$ is the same as the forcing
$K^*_{\iterq_*,\zeta_*,p_*}$.
\end{claim}

\begin{claim}\label{3.8} 
Let ${\bf q} \perp \bar{\Gamma}$.
Let each $\Gamma_\zeta$ have free amalgamation.
Let ${\mathfrak K}^* =(K^*_{\lambda,\zeta(*)}(\Gamma), \leq)$. 
Then ${\mathfrak K}^* * \name{{\mathbb P}}$ 
is a $(<\lambda)$-strategically complete forcing
of cardinality $\leq \lambda^{<\lambda}$ and in 
${\bf V}^{{\mathfrak K}^* * \name{{\mathbb P}}}$,
in the game $\Game_{\kappa}({\bf q})$ the empty player does not have 
a winning strategy.
\end{claim}

\proof Assume towards a contradiction that we have
for some $\kappa =\cf(\kappa) \in (\aleph_0,\lambda)$ and some
 $\name{St}$ that 
$\iterq \Vdash_{{\mathfrak K}^*} p \Vdash_{\name{\mathbb P}}$
``$\name{St}$ is a name for a winning strategy for 
the empty player in the game $\Game_\kappa({\bf q})$''.

Since ${\mathbb P}$ may be replaced by ${\mathbb P}_{\geq p}$, we may 
assume $p=1$.
We now choose by induction on $\alpha < \kappa$ a triple $(\iterq_\alpha,
\name{\gc}_\alpha, \name{\gd}_\alpha)$
satisfying

\begin{myrules}
\item[($t_1$)] $\iterq_\alpha \in K^*$ is $\leq_{pr}$-increasing weakly continuous (in $K^*$),
\item[($t_2$)] $\iterq_0 = \iterq$,
\item[($t_3$)] $\name{\gc}_\alpha$ is a ${\mathbb P}_0^{\iterq_\alpha}$-name and 
$\name{\gd}_\alpha$ is a ${\mathbb P}_0^{\iterq_\alpha}$-name
(not just ${\mathbb P}_\eps^{\iterq_\alpha}$-names for some $\eps$),
\item[($t_4$)] $\iterq_\alpha \Vdash_{{\mathfrak K}^*} 1 \Vdash_{\name{\mathbb P}}
\mbox{``}\la \name{\gc}_\beta,\name{\gd}_\beta \such \beta < \alpha \r$ is
a play of $\Game_\kappa (\bf q)$ in which the empty player
uses the  strategy $\name{St}$''.
\end{myrules}

Assume that we arrive at stage $\alpha < \kappa$. First we let 
$\iterp_{\alpha,0}$ be $\iterq$ 
if $\zeta= 0$, and we let 
$\iterp_{\alpha,0}$ be  $\iterq_{\alpha-1}$ if $\alpha$ is a successor ordinal,
if $\alpha$ is a limit ordinal we take the direct limit that is the union of 
$\la \iterq_\beta \such \beta < \alpha\r$.
By Claim~\ref{3.3}, part (3), it exists and is
 $\leq_{pr}$-above every $\iterq_\beta$ for $\beta < \alpha$. We now choose 
$\iterp_{\alpha,1}$,
$\bar{{\mathcal C}}= 
\la \name{\mathcal C}_{\alpha,\zeta} \such \zeta <\zeta(*)\r$.

\nothing{
\smallskip
{\bf Version 1}
We first work with ``each $\Gamma_\zeta$ 
has free amalgamation'' (see \ref{2.1}(5)).
}
Let us define $\iterp_{\alpha,1}$ by choosing $\lambda_{\alpha,\zeta} =
|{\mathbb P}_\zeta[\iterp_{\alpha,0}]|^{\aleph_0}$ and ${\mathbb P}_{\alpha,\zeta} = 
\{ f \such f $ is a finite function from $\lambda_{\alpha,\zeta}$ to 
$\{0,1\}\}$ and
${\mathbb P}_{\zeta}[\iterp_{\alpha,1}] =
{\mathbb P}_\zeta[\iterp_{\alpha,0}] \times {\mathbb P}_{\alpha,\zeta}$
and $\name{\eta}^{\iterp_{\alpha,0}} = \name{\eta}^{\iterp_{\alpha,1}}$.
Since $\Gamma_\zeta$ has free amalgamation $\iterp_{\alpha,0} \leq_{pr}
\iterp_{\alpha,1} \in K^*$.
Clearly $\Vdash_{{\mathbb P}_{0}[\iterp_{\alpha,1}]}
\la \name{\gc}_\beta,\name{\gd}_\beta \such \beta < \alpha \r$ is
an initial segment of a  play of $\Game_\kappa (\bf q)$.

So, by absoluteness, 
 $\Vdash_{{\mathbb P}_{0}[\iterp_{\alpha,1}]}
\la \name{\gd}_\beta \such \beta < \alpha \r$ is
$\leq_{\bf q}$-linearly ordered.
Hence by \ref{1.14} there is a ${\mathbb P}_0[\iterp_{\alpha,1}]$-name
$\name{\mathcal C}_{\alpha}$ of 
a $\bf q$-nice family which (is forced to) include
 $\{\name{\gd}_\beta \such \beta < \alpha\}$.

We now can choose $\name{\mathcal C}_{\alpha, \zeta}$ by induction on $\zeta<\zeta(*)$
such that 
\begin{myrules}
\item[(a)]$\Vdash_{{\mathbb P}_{\zeta}
[\iterp_{\alpha, 1}]} \name{\mathcal C}_{\alpha,\zeta} $ 
is a $\bf q$-nice family,
\item[(b)] $\name{\mathcal C}_{\alpha,0} = \name{\mathcal C}_\alpha$ and $\eps < \zeta 
\rightarrow
\Vdash_{{\mathbb P}_{\zeta}[\iterp_{\alpha, 1}]}\name{\mathcal C}_{\alpha,\eps} 
\subseteq \name{\mathcal C}_{\alpha,\zeta} $, 
\item[(c)] if $\zeta = \eps +1$ then 
$({\mathbb P}_{\zeta(*)}^{\iterp_{\alpha,1}} * {\bbforcing}_{\bf q}(\name{\mathcal C}_{\alpha,\eps}),
{\mathbb P}_{\zeta(*)}^{\iterp_{\alpha,1}} * {\bbforcing}_{\bf q}(\name{\mathcal C}_{\alpha,\zeta}),
\name{\eta}_\zeta^{\iterp_{\alpha,0}}) \in \Gamma_\eps$ is
above $\iterp_{\alpha,1}$.
\end{myrules}
We carry  out the induction by using 
the Definition~\ref{2.9} Part 1 of 
$\perp$ 
and thereafter  Claim~\ref{1.17} to get from
directedness to niceness. Or we use
Claim~\ref{2.10}.

\medskip

Now define $\iterp_{\alpha,2} \in K^*$ by choosing ${\mathbb P}_\zeta[\iterp_{\alpha,2}] =
({\mathbb P}_\zeta[\iterp_{\alpha,1}]) * {\bbforcing}(\name{{\mathcal C}}_{\alpha,\zeta})$.
Clearly $\iterp_{\alpha,1} \leq_{pr} \iterp_{\alpha,2}\in K^*$ so 
${\mathbb P}_{\zeta(*)}[\iterp_{\alpha,2}] =
({\mathbb P}_{\zeta(*)}[\iterp_{\alpha,1}]) * {\bbforcing}(\name{{\mathcal C}}_{\alpha,\zeta(*)})$,
where $\name{\mathcal C}_{\alpha,\zeta(*)} = 
\bigcup \{ \name{\mathcal C}_{\alpha,\zeta} \such \zeta < \zeta(*)\}$.

It satisfies the c.c.c by Claim~\ref{1.17}(1).
Lastly let us choose a move 
$\name{\bar{\gc}}_\alpha$  for the
non-empty player the generic of ${\bbforcing}(\name{{\mathcal C}}_{\alpha,\zeta(*)})$,
which is an  ${\bbforcing}(\name{{\mathcal C}}_{\alpha,\zeta})$-name for every 
$\zeta < \zeta(*)$ that is, it is $\bigcup\{\bar{\gc}^p \restriction n^p \such
p=(n^p,\bar{\gc}^p)
\in \name{G}({\bbforcing}_{\bf q}({\mathcal C}_{\alpha,\zeta}))\}$. So
half of clause ($t_3$) holds.

So there is $\iterp_{\alpha,3} \in K$ which is by $\leq_{K^*}$ 
above $\iterp_{\alpha,2}$ and forces 
$\name{St}(\la \name{\gc}_\beta \such \beta \leq \alpha \r)$
to be equal to some ${\mathbb P}_\zeta^{\iterp_{\alpha,\zeta}}$-name 
$\name{\bar{\gd}}^\alpha$.

However we cannot in general choose 
$\iterq_\alpha= \iterp_{\alpha,3}$, $\name{\bar{\gd}}_\alpha = \name{\bar{\gd}}^\alpha$ 
as possible $\name{\bar{\gd}}^\alpha$ is not a ${\mathbb P}_0[\iterp_{\alpha,3}]$-name.
Since $\zeta(*)$ is uncountable,
for some $\zeta_\alpha <\zeta(*)$ the name $\name{\bar{\gd}}^\alpha$ is a 
${\mathbb P}_{\zeta_{\alpha}}[\iterp_{\alpha,3}]$-name.
So we can find $\iterp_{\alpha,4}$ such that $\iterp_{\alpha,3} \leq_{\pr}
\iterp_{\alpha,4}$ and  $\la{\mathcal C}^*_{\alpha,\zeta} \such \zeta \in [\zeta_\alpha,
\zeta(*))\r$ such  that for every $\zeta \in [\zeta_\alpha,
\zeta(*))$ we have ${\mathcal C}^*_{\alpha,\zeta}$ is a ${\mathbb P}_\zeta[\iterp_{\alpha,4}]$-name 
of a $\bf q$-nice family to which $\bar{\gd}^\alpha$ is forced to belong increasing with 
$\zeta$ such that
$$({\mathbb P}_\eps[\iterp_{\alpha,4}],
{\mathbb P}_{\eps+1}[\iterp_{\alpha,4}], \name{\eta}_\eps)
\leq_\Gamma
({\mathbb P}_\eps[\iterp_{\alpha,4}] * {\bbforcing}_{\bf q}(\name{\mathcal C}^*_{\alpha,\eps}),
{\mathbb P}_{\eps+1}[\iterp_{\alpha,4}] * {\bbforcing}_{\bf q}
(\name{\mathcal C}^*_{\alpha,\eps+1}),
\name{\eta}_\eps).
$$

Let $\iterp_{\alpha,5} \in K^*$ be defined as follows:
${\mathbb P}_\eps[\iterp_{\alpha,5}]$ is
${\mathbb P}_\eps[\iterp_{\alpha,4}]$ if $\eps < \zeta_\alpha$ and is
${\mathbb P}_\eps[\iterp_{\alpha,4}] * {\bbforcing}_{\bf q}({\name{\mathcal C}}^*_{\alpha,\eps})$
 if $\eps \in [\zeta_\alpha,\zeta(*))$.
${\mathcal C}^*_{\alpha,\zeta(*)}  =
 \bigcup\{ {\mathcal C}^*_{\alpha,\zeta} \such \zeta < \zeta(*)\}$.

\medskip

Now clearly $\iterp_{\alpha,4} \leq_{\pr} \iterp_{\alpha,5}$ and let $\name{\bar{\gd}}_\alpha$ 
be the generic of ${\bbforcing}_{\bf q}(\name{\mathcal C}_{\alpha,\zeta(*)}^*)$.
Now by Claim~\ref{1.20} we have that
$\Vdash_{{\mathbb P}_0[\iterp_{\alpha,5}]/{\mathbb P}_\eps[\iterp_{\alpha,1}]}$
 ``$\name{\bar{\gd}}_\alpha$ is generic for 
${\bbforcing}_{\bf q}(\name{\mathcal C}^*_{\alpha,\zeta(*)})$''.
Lastly we define $\iterq_\alpha$ by
\begin{equation}\label{qalpha}
{\mathbb P}_\eps[\iterq_\alpha] = \left\{
\begin{array}{ll}
({\mathbb P}_\eps[\iterp_{\alpha,5}] * 
\mbox{  (the subforcing of${\bbforcing}_{\bf q}(\name{\mathcal C}^*_{\alpha,\zeta(*)})$}&\\
\mbox{  generated by }
 \name{\bar{\gd}}_\alpha), &  \mbox{ if } \eps \leq  \zeta_\alpha,\\
{\mathbb P}_\eps[\iterp_{\alpha,5}], & \mbox{ if } \eps >  \zeta_\alpha.
\end{array}
\right.
\end{equation}

Now in step $\kappa$, also 
$\bar{\gc}_\kappa$ exists and 
is $\leq_{\mathbb P}$, below all   the $\bar{\gd}_\alpha$, $\alpha <\kappa$, which is
a contradiction to $St$ being a winning strategy for the empty player.
\proofend

\nothing{
$$ * * * $$
{\bf Version 2} {\bf OPEN WORK}\\
No extra assumption on $\Gamma_\zeta$. So we assume only the $\Gamma_\zeta  \perp {\bf q}$ for
each $\zeta< \zeta(*)$ and that $\Gamma_\zeta$ has amalgamation.

This needs some extra assumption to ensure the c.c.c.\ in limit $\zeta$ of cofinality 
$\aleph_1$.
Probably we should add ``${\mathbb P}_{\alpha,\zeta}/{\mathbb P}_\zeta[\iterq_{\alpha,1}]$ 
 has the Knaster property''. THINK.
It is not essential for the rest.
We choose 
by induction on $\zeta<\zeta(*)$ a pair $({\mathbb P}_{\alpha,\zeta}, 
\name{\mathcal C}_{\alpha,\zeta})$ such that
\begin{myrules}
\item[(a)] ${\mathbb P}_{\alpha,\zeta}$ is a c.c.c.\ forcing notion  and 
element of ${\mathcal H}(\lambda)$,
\item[(b)] if $\xi < \zeta$ then ${\mathbb P}_{\alpha,\zeta} \lessdot {\mathbb P}_{\alpha,\zeta}$,
\item[(c)] ${\mathbb P}_\zeta[\iterp_{\alpha,0} ]\lessdot {\mathbb P}_{\alpha,\zeta}$,
 moreover the quotient in Knaster,
\item[(d)] ${\mathbb P}_{\alpha,\zeta} \cap {\mathbb P}_{\zeta(*)}[\iterp_{\alpha,\zeta}] 
= {\mathbb P}_\zeta[\iterp_{\alpha,\zeta}]$,
\item[(e)] for each $\zeta< \zeta<(*)$ there is a ${\mathbb P}_{\alpha,\zeta}'$ 
such that ${\mathbb P}_{\alpha,\zeta}' \lessdot {\mathbb P}_\zeta$,
and $\{{\mathbb P}_{\alpha,\xi} \such \xi < \zeta\} \subseteq {\mathbb P}_{\alpha,\zeta}'$,
$|{\mathbb P}_{\alpha,\zeta}'| = |{\mathbb P}_{\alpha,\zeta}'|^{\aleph_0}$
and $\lambda_{\alpha,\zeta} =  |{\mathbb P}_{\alpha,\zeta}'|^{\aleph_0}$
and ${\mathbb P}_{\alpha,\zeta} = {\mathbb P}_{\alpha,\zeta}' 
\ast {\rm Cohen}_{\lambda_{\alpha,\zeta}}$.
Here the Cohen forcing is with finite partial functions from $\lambda_{\alpha\zeta}$ to 2,
ordered by inclusion.
\item[(f)] $\name{\mathcal C}_{\alpha,\zeta}$ is a ${\mathbb P}_{\alpha,\zeta}$-name 
of a $\bf q$-nice family including
$\{ \name{\bar{\gc}}_\beta \such \beta < \alpha\}$,
\item[(g)] if $\xi < \zeta$ then $\Vdash_{{\mathbb P}_{\alpha,\zeta}}
\name{\mathcal C}_{\alpha,\xi} \subseteq \name{\mathcal C}_{\alpha,\zeta}$,
\item[(h)] if $\zeta= \xi +1$ then
$$({\mathbb P}_\xi[\iterp_{\alpha,0}],
{\mathbb P}_\zeta[\iterp_{\alpha,0}],\name{\eta}_\xi)
\leq_\Gamma ({\mathbb P}_{\alpha,\xi} \ast {\bbforcing}_{\bf q}(\name{{\mathcal C}}_{\alpha,\xi}),
{\mathbb P}_{\alpha,\zeta} \ast {\bbforcing}_{\bf q}(\name{{\mathcal C}}_{\alpha,\zeta}),
\name{\eta}_\xi)$$
usually
$$({\mathbb P}_\xi[\iterp_{\alpha,0}],
{\mathbb P}_\zeta[\iterp_{\alpha,0}],\name{\eta}_\xi)
\leq_\Gamma 
({\mathbb P}_{\alpha,\xi},
{\mathbb P}_{\alpha,\zeta},\name{\eta}_\xi)
\leq_\Gamma
({\mathbb P}_{\alpha,\xi} \ast {\bbforcing}_{\bf q}(\name{{\mathcal C}}_{\alpha,\xi}),
{\mathbb P}_{\alpha,\zeta} \ast {\bbforcing}_{\bf q}(\name{{\mathcal C}}_{\alpha,\zeta}),
\name{\eta}_\xi)$$
\end{myrules}
Carrying out the induction we ingore case (d), as it is just a matter of renaming and 
we describe it by cases.
Case 1: $\zeta=0$. 
Let 
${\mathbb P}_{\alpha,\zeta} = ({\mathbb P}_0[\iterp_{\alpha,0}])  
\ast {\rm Cohen}_{\lambda_{\alpha,0}}$, where
 $\lambda_{\alpha,0} =  |{\mathbb P}_0[\iterp_{\alpha,0}] |^{\aleph_0}$.
Now clearly 
\begin{equation*}
\begin{split}
(\iterp_{\alpha,0}, \emptyset) \Vdash_{{\mathfrak K}^* \ast \name{\mathbb P}}
 & \la (\name{\bar{\gc}}_\beta,\name{\bar{\gd}}_\beta\such \beta < \alpha \r
\mbox{ is an initial segment of a play of $\Game_\kappa(\bf q)$},
\\
&
\mbox{hence }
\{ \name{\bar{\gc}}_\beta,\name{\bar{\gd}}_\beta \such \beta < \alpha \} \subseteq {\bf Q}_{\bf q}
\mbox{ is a $\leq_{\bf q}$-linearly ordered set.''}
\end{split}
\end{equation*}
As $\name{\bar{\gc}}_\beta$, $\name{\bar{\gd}}_\beta$ for $\beta < \alpha$ are 
$({\mathbb P}_0[\iterp_{\alpha,0}]$-names clearly by absoluteness
\begin{equation*}
\begin{split}
\Vdash_{{\mathfrak K}^* \ast \name{\mathbb P}}
 & \{ \name{\bar{\gc}}_\beta,\name{\bar{\gd}}_\beta \such \beta < \alpha \} 
\subseteq {\bf Q}_{\bf q}
\mbox{ is a $\leq_{\bf q}$-linearly ordered set.''}
\end{split}
\end{equation*}
Hence this holds also for $\Vdash_{{\mathbb P}_{\alpha,0}}$, 
hence by Claim~\ref{1.17} (1) there is a ${\mathbb P}_{\alpha,\zeta}$-name 
$\name{\mathcal C}_{\alpha,\zeta}$ of a $\bf q$-nice family containing 
$\{ \name{\bar{\gc}}_\beta,\name{\bar{\gd}}_\beta \such \beta < \alpha \}$.
Case 2: $\zeta = \xi+1$ FILL
Case 3: $\zeta$ is a limit of countable cofinality. FILL
Case 4: $\zeta$ is a limit of uncountable cofinality. FILL
{\bf not yet worked through.}}

\begin{theorem}\label{3.9}
1. Assume that ${\bf q} \perp \bar{\Gamma}$, as in \ref{1.1},
$\bar{\Gamma} $ be a sequence of Knaster forcing bigness notions
$\Gamma_\zeta$ with free amalgamation and the $(<\lambda)$-LS. as in 2.1,
$K^*= K^*_{\lambda,\zeta(*)}(\bar{\Gamma})$. Then  recalling $\invgm(\bf q)$
from Definition~\ref{0.1}(4)
we have
$\Vdash_{{\mathfrak K}^* * \name{\mathbb P}} \invgm(\bf q) = \lambda$.

2. If $\bar{\Gamma}$ has free amalgamation and $\zeta(*) = \sup\{ \zeta < \zeta(*) \such
\Gamma_\zeta = \Gamma\}$, then 
$\Vdash_{{\mathfrak K}^* * \name{\mathbb P}} $ ``for some $Y \subseteq {}^\omega \omega$
$|Y| = \cf(\zeta(*))$, for every $\nu \in {}^\omega \omega$
(of $\nu \in {}^\gamma \omega$ for some $\gamma <\zeta(*)$) some
$\eta \in Y$ is $\Gamma$-big over $\{\nu\}$''. For $\Gamma =\Gamma_{\rm ud}$ 
this means $\eta \not\leq^* \nu$. 
\end{theorem}

\proof 1. Assume towards a contradiction that we have
$\kappa = \cf(\kappa) \in (\aleph_0,\lambda)$
and $\iterq \in {\mathfrak K}^*$ such that

$$\iterq \Vdash_{{\mathfrak K}^*} p \Vdash_{\name{\mathbb P}} \invgm(\bf q) = \kappa.$$

So, for some $\name{St}$ we have that 
$\iterq \Vdash_{{\mathfrak K}^*} p \Vdash_{\name{\mathbb P}}$
``$\name{St}$ is a name for a winning strategy for 
the empty player in the game $\Game_\kappa(\bf q)$''.
But this contradicts the previous claim.

2. is Claim~\ref{3.5}.
\proofend

Equation~\eqref{qalpha} is probably easier to 
state in terms of complete boolean algebras, because we speak about 
the generated subforcing. The following remark gives the
translation.
 
\begin{remark}\label{3.10}
1. If $\iterq \in K$ and if $\hat{\iterq}$ is defined by $\name{\eta}_\zeta^{\hat{\iterq}} =
\name{\eta}_\zeta^\iterq$ and
${\mathbb P}_\zeta^{\hat{\iterq}} = $ the completion of
${\mathbb P}_\zeta^\iterq$  , then 
$\iterq \leq_{\pr} \hat{\iterq}$.

2. For a forcing notion ${\mathbb P}$ let its completion $\hat{\mathbb P}$
consist of all 
$\{{\mathcal J} \such {\mathcal J} \subseteq {\mathbb P} \mbox{ is an antichain}\}$, 
ordered by ${\mathcal I} \leq {\mathcal J}$ iff $(\forall p \in {\mathcal I}) 
(\exists q \in {\mathcal J})(p \leq_{\mathbb P} q)$.

Identifying $\{p\}$ with $p$ we have that  ${\mathbb P} \lessdot \hat{\mathbb P}$,
${\mathbb P}$ is dense in $\hat{\mathbb P}$ and $\hat{\mathbb P}$ is a quasi order, 
there are usually many equivalent elements.
\end{remark}

\section{Examples}\label{S4}

\begin{definition}\label{4.1}
The Hechler forcing in the creature framework.
${\bf q}_H = {\bf q}=  (K,\Sigma,\norm,\val)$ is given by
$K=\{\gc =(\mdn,\mup,h)\such (\exists n<\omega) \mdn= n, \mup= n+1,
h\colon \{n\} \to \omega \}$.
$\dom(\Sigma)= K$ and 
$\Sigma(\mdn,\mup,h) = \{(\mdn,\mup,g) \such g(\mdn) \geq  h(\mdn)\}$,
and $\norm(\mdn,\mup,h)= \mdn$, $\val(\gc) = \{(n,h(n))\}$.
\end{definition}

\begin{observation}\label{4.2}
1. ${\bf q}_H$ is a forgetful pcrf as in Definitions~\ref{1.1} and \ref{1.2}.  
\\
2. $\invcm({{\bf q}}_H) = \invgm({{\bf q}}_H) = \gb$.
\\
3. $\leq_{\bf q}$ is the usual $\leq^*$ on ${}^\omega \omega$ if we identify
$\bar{\gc} = \la \gc_n \such n<\omega \r
\in {\bf Q}$  with $f_{\bar{\gc}} \in {}^\omega \omega$, $f_{\bar{\gc}}(\mdn)=
h(\mdn)$ if $(\mdn,\mup,h) \in \rge(\bar{\gc})$.
\\
4. ${\bf q}_H$ is not so interesting, as in  our framework for
${\mathcal C}$ being ${\bf q}$-nice,
${\bbforcing}_{\bf q}(\mathcal C)$ will add
a dominating real and hence increase $\gb$.
\end{observation}

\begin{definition}
${\bf q}_{\rm gr} =(K,\Sigma,\norm,\val)$ where
\begin{myrules}
\item[(a)] $K=\{\gc \such \gc=(\mdn^\gc,\mup^\gc,u^\gc), \mdn^\gc < \mup^\gc < \omega,
\emptyset \neq u^\gc \subseteq [\mdn^\gc,\mup^\gc)\}$
\item[(b)] $\gc\in \Sigma(\gc_0, \dots,\gc_{k-1})$ iff $\gc_\ell
\in K $ and for $\ell < k-1$, $\mup^{\gc_\ell} = \mdn^{\gc_{\ell+1}}$,
and $\mdn^{\gc} = \mdn^{\gc_{0}}$ and
$\mup^{\gc} = \mup^{\gc_{k-1}}$
and for some non-empty $v \subseteq\{0,\dots,k-1\}$ we have
$u^\gc = \bigcup_{\ell \in v} u^{\gc_\ell}$.
\item[(c)]
$\norm(\gc) = \mdn^\gc$ for $\gc \in K$, $\val(\gc)= u_\gc$.
\end{myrules}
\end{definition}

Remark: ${\bf Q}^{tr}_{\bf q_{\rm gr}}$
 is called Matet forcing, see  \cite{Blasstoronto}.
Close to Hindman \cite{hindman:sums}.

\begin{observation}\label{4.4}
1. ${\bf q}_{\rm gr}$ is a forgetful pcrf.\\
2. $\invgm({\bf q}_{\rm gr}) \leq \gro$.
\end{observation}

\proof 1. Straightforward.\\
2. Let $\kappa= \gro$. By the definition of $\gro$ we can find 
$\la {\mathcal A}_\alpha \such \alpha < \kappa \r$ such that
\begin{myrules}
\item[(a)] ${\mathcal A}_\alpha \subseteq [\omega]^{\aleph_0}$,
\item[(b)] if $A \subseteq^* B\in {\mathcal A}_\alpha$ and $A$ is infinite, 
then $A\in {\mathcal A}_\alpha$,
\item[(c)] if $n_0 <n_1< \dots$ then for some infinite set $w \subseteq\omega$ we have
$\bigcup\{[n_i,n_{i+1}) \such i \in w \} \in {\mathcal A}_\alpha$,
\item[(d)] for no $B \in [\omega]^{\aleph_0}$ do we have $(\forall \alpha < \kappa)
(\exists A \in {\mathcal A}_\alpha) (B \subseteq A)$.
\end{myrules}

Now we describe a strategy ${\bf ST}$ for the empty player in the game
 $\Game_\kappa({\bf q}_{\rm gr})$. If the play up to the $\beta$-th move 
has been $\la (\bar{\gc}_\beta,\bar{\gd}_\beta)
\such \beta <\alpha\r$ and the non-empty player
chooses $\bar{\gc}_\alpha$, then the empty player chooses
$\bar{\gd}_\alpha$ such that
\begin{equation}
\label{stern}
\bar{\gc}_\alpha \leq_{{\bf q}_{\rm gr}}
\bar{\gd}_\alpha \, \wedge \, \bigcup\{u[\bar{\gd}_{\alpha,n}] \such n<\omega \} 
\in {\mathcal A}_\alpha.
\end{equation}

This is possibe as if we let $m_n = \mdn[\bar{\gc}_{\alpha,n}]$ 
then $m_n <m_{n+1} < \dots$ and $m_0 = 0$ hence for 
some infinite $w \subseteq \omega$ we have $\bigcup\{[m_i,m_{i+1})
\such i \in w\} \in {\mathcal A}_\alpha$ and 
letting $i_0 < i_1 < i_2 < \dots$ list $w$, we choose $\gd_{\alpha,\ell}
=
(m_{i_{\ell -1} +1}, m_{i_\ell +1}),
u[\bar{\gd}_{\alpha, \ell }]= u[\bar{\gc}_{\alpha.i_\ell}]$,
stipulating $m_{i_{-1} +1} = 0$ so
$\la \gd_{\alpha,\ell} \such \ell < \omega \r$ is as required.

\smallskip
Suppose that   $\la (\bar{\gc}_\beta,\bar{\gd}_\beta)
\such \beta <\kappa\r$ have been chosen and there is $\bar{\gc}_\kappa
= \la \gc_{\kappa,n} \such n < \omega \r
\in {\bf Q}^{tr}_{{\bf q}_{\rm gr}}$ such that $\alpha < \kappa \rightarrow 
\bar{\gc}_\alpha \leq_{{\bf q}_{\rm gr}}
\bar{\gc}_\kappa$. We shall derive a contradiction.

Let the set $B = \bigcup\{u[\gc_{\kappa,n}] \such n<\omega \}$.
This an infinite subset of $\omega$ by the definition of ${\bf q}_{\rm gr}$.
Also for each $\alpha < \kappa$ clearly $\bar{\gc}_\kappa$ is above $\bar{\gc}_{\alpha+1}$,
hence above $\gd_\alpha$, hence $B \subseteq^* \bigcup \{u[\gd_{\alpha,n}] \such n < \omega \}
\in {\mathcal A}_\alpha$, so $B \in {\mathcal A}_\alpha$.
Hence $B \in \bigcup \{{\mathcal A}_\alpha \such \alpha < \kappa \}$, contradicting 
the choice of $\la {\mathcal A}_\alpha \such \alpha < \kappa \r$
\proofendof{\ref{4.4}}

\nothing{
The rest of this section is not needed for $\gb^+ < \gro$.
\begin{claim}\label{4.5}
Assume
\begin{myrules}
\item[(a)] $({\mathbb P}_1, {\mathbb P}_1^+,\name{\eta}) \in \Gamma_{\rm tow}$,
\item[(b)] ${\mathbb P}_1 \lessdot {\mathbb P}_2$ and ${\mathbb P}_2 /{\mathbb P}_1$ 
satisfies the Knaster condition,
\item[(c)] $\Vdash_{{\mathbb P}_2} $''$ \name{\mathcal C}$ is a nice ${\bf q}_{\rm gr}$-family'',
\item[(d)] $\lambda = (|{\mathbb P}_1^+| + |{\mathbb P}_2|)^{\aleph_0}$,
\item[(e)] ${\mathbb P}' = {\mathbb P}_1^+ \times_{{\mathbb P}_1} {\mathbb P}_2$,
\item[(f)] ${\mathbb P}_2^+ = {\mathbb P}' \times {\rm Cohen}_\lambda$.
\end{myrules}
for notational simplification  ${\mathbb P}_2 \cap {\mathbb P}_1^+ = {\mathbb P}_1$
and ${\mathbb P}_2, {\mathbb P}_1^+ \lessdot {\mathbb P}'$ and we identify $p \in {\mathbb P}_1^+$
with $(p,\emptyset)$.
Then
\begin{myrules}
\item[($\alpha$)]
${\mathbb P}'$ satisfies the c.c.c\ and has cardinality $\leq \lambda^{\aleph_0}$,
\item[($\beta$)]
there is a ${\mathbb P}^+_2$-name $\name{{\mathcal C}}^+$ such that
$\Vdash_{{\mathbb P}_2^+} \name{{\mathcal C}}^+$ is a ${\bf q}_{\rm gr}$-nice family extending 
$\name{{\mathcal C}}$,
and $({\mathbb P}_1, {\mathbb P}_1^+, \name{\eta}) \leq_{\Gamma_{\rm tow}}
 ({\mathbb P}_2 * {\bbforcing}_{{\bf q}_{\rm gr}}(\name{\mathcal C}), {\mathbb P}_2^+
 * {\bbforcing}_{{\bf q}_{\rm gr}}(\name{\mathcal C}^+), \name{\eta})$.
\end{myrules}
\end{claim}
\proof We can force first with ${\mathbb P}_1$ so without loss of generality
\begin{equation}
{\mathbb P}_1 \mbox{ is trivial, so } {\mathbb P}_2^+ = {\mathbb P}_1^+ \times {\mathbb P}_2
\times {\rm Cohen}_\lambda
\end{equation}
so we can write it in this form.
We can also present ${\rm Cohen}_\lambda$ as 
$\{f \such f$  a finite function from $\lambda $  to ${}^{\omega >} \omega\}$
ordered naturally and let $\eta_\alpha$  for $\alpha <\lambda$  be
$$
\bigcup\{ f(\alpha) \such p \in \name{G}_{{\rm Cohen}_\alpha}, \alpha \in \dom(p)\}.
$$
Let $G_1^+ \times G_2 \times G_{\rm Cohen} \subseteq {\mathbb P}_1^+ \times
{\mathbb P}_2
\times {\rm Cohen}_\lambda$ be generic over $V$.
In $V[G_2]$ let 
$\la (\bar{\gc}_\alpha, {\mathcal J}_\alpha, \name{h}_\alpha ) \such 
\alpha < \lambda\r$ list the set of tuples $(\bar{\gc}, {\mathcal J}, h)$ such that
\begin{equation}
\begin{split}
&\bar{\gc}_\alpha \in {\mathcal C} = \name{\mathcal C}[G_2],\\
&\bar{\mathcal J}_\alpha = \la {\mathcal J}_n \such n<\omega \r \\
&{\mathcal J}_n \subseteq {\bbforcing}_{{\bf q}_{\rm gr},{\mathcal C}}(\bar{\gc}_\alpha)
\mbox{ is predense in }({\bbforcing}_{{\bf q}_{\rm gr}}(\mathcal C))^{{\bf V}[G_2]},\\
&\bar{h} = \la h_n \such n <\omega \r, h_n\colon {\mathcal J}_n \to {}^{\omega >} \omega,\\
&\mbox{  be such that for some }
{\bbforcing}_{{\bf q}_{\rm gr}}({\mathcal C})\mbox{-name } \name{f}
\mbox{ of a member of } {}^\omega \omega\\
&\mbox{ we have } (n<\omega \wedge q \in {\mathcal J}_n)
\rightarrow q \Vdash_{ {\bbforcing}_{{\bf q}_{\rm gr}}(\mathcal C)} \name{f} \restriction n
= h_n(q).
\end{split}
\end{equation}
{\bf FILL  as in 1.14 (3) 1.17 (1).\\
But we need to preserve also $\Gamma_{\rm tow}$!
This is a new complaint, because up to
now I did not take care of the fact, that now we have $\Gamma_{\rm tow}$ and not $\Gamma_{\rm ud}$.}
\begin{definition}\label{4.6}
1. For $n\geq 1$ let ${\bf q}_n =(K_n, \Sigma_n, \norm_n)$ be
\begin{myrules}
\item[(a)] $\gc \in K_n$ iff $\gc = (\mdn,\mup, k_0, \dots k_{n-1})$ with
$\mdn <\mup < \omega$, $k_\ell \in [\mdn,\mup )$ and $k_\ell < k_{\ell+1}$ for $\ell <n-2$.
\item[(b)] $\bar{c} \in \Sigma_n(\gc_0,\dots , \gc_{n-1})$ iff 
$n \geq 1$, $\mup(\gc_m) =  \mdn(\gc_{m+1})$ for $m < n-2$ and for $\ell < n$,
$k_\ell(\gc_n) \in \{k_\ell(\gc_m) \such m < n \}$.
\item[(c)] $\norm_n(\gc) = \mdn(\gc)$.
\end{myrules}
2. Define $F_\ell \colon {\bf Q}_n \to [\omega]^{\aleph_0}$ by $F_\ell (\bar{\gc}) =
\{ k_\ell (\gc_n) \such n < \omega \}$.
\end{definition}
Remark: $\bf q_n$ is a way to formalize ${\gh}({}^n([\omega]^{\aleph_0}))$.
\begin{claim}\label{4.7}
1. ${\bf q}_n$ fits Definition~\ref{1.1}
2. $\invgm(\bf q) \leq \kappa$ iff in the following game 
the empty player has a winning strategy.
A play lasts $\kappa+1$ moves. In the $\alpha$-th move the non-empty player chooses 
$A_{\alpha,0}, A_{\alpha, 1}, \dots, 
A_{\alpha, n-1}
\in [\omega]^{\alpha_0}$ such that $(\beta < \alpha \wedge
\ell < n)
\rightarrow A_{\alpha,\ell} \subseteq^* B_{\alpha,\ell}$ and then the empty player chooses
 $B_{\alpha,0}, B_{\alpha, 1}, \dots, 
B_{\alpha, n-1}
\in [\omega]^{\alpha_0}<$  such that 
$B_{\alpha,\ell} \subseteq^* A_{\alpha,\ell}$.
3. $n<m \rightarrow {\bf q}_n \not\leq {\bf q}_m$.
4. $\bar{\gc} \leq_{{\bf q}_n} \bar{\gd}$ implies that $F_\ell (\bar{\gd}) 
\subseteq F_\ell(\bar{\gc})$.
5. If ${\mathcal C}$ is pre-${\bf q}_n$-directed and $\ell <n$, 
then the intersection of any finitely many members of $F_\ell(\mathcal C) 
:= \{F_\ell(\bar{\gc}) \such  \bar{\gc} \in {\mathcal C} \}$ is infinite,
and $\subseteq \omega$. 
6. If ${\mathcal C}$ is ${\bf q}_n$-directed and $\ell <n$, 
then $F_\ell(\mathcal C)$ is directed by $\supseteq^*$.
7. If ${\mathcal C}$ is a ${\bf q}_n$-nice family  and $\ell <n$, 
then $F_\ell(\mathcal C)$ generates a Ramsey ultrafilter on $\omega$,
if $\mathcal C$ is hereditary (i.e. downward closed under $<^*$), then $F_\ell(\mathcal C)$ is
equal to a Ramsey ultrafilter on $\omega$ and the 
$F_\ell(\mathcal C)$, $\ell <n$, are pairwise $\leq_{RK}$-incomparable.
8. The inverses of the implications in 6. and in 7.
\end{claim}
\proof{\bf EASY}
\begin{claim}\label{4.8}
${\bf q}_{n_1} \perp {\bf q}_{n_2}$
or ${\bf q}_{n_1} \not\leq {\bf q}_{n_2}$ .
\end{claim}
{\bf Saharon, $\not\leq$ and $\leq$ are not yet defined.}
}

\section{On the consistency of $\gb \ll \gro$}\label{s5}

We now turn to a specific problem.

\begin{theorem}\label{5.1}
Let $\aleph_0 < \kappa = \cf(\kappa) < \lambda = \lambda^{<\lambda}$.
Then for some notion of forcing ${\mathfrak K} \ast \name{\mathbb P}$
of cardinality $\lambda$ in ${\bf V}^{{\mathfrak K} \ast \name{\mathbb P}}$
we have $\gb = \kappa$ and $\gro = \lambda$.
\end{theorem}

\begin{definition}\label{5.2}
For ${\mathcal C} \subseteq {\bf Q}^{tr}_{{\bf q}_{\rm gr}}$ 
directed and a a filter $D({\mathcal C})
=        \{ \{\mdn(\gc_n) \such n \in \omega\} \such 
\bar{\gc} = \la \gc_n\such n \in \omega \r \in {\mathcal C}\}$
 we define
\begin{equation}
\begin{split}
{\bf Q}^{tr}_D =&\{(w,\bar{\gc}) \such
w \in [\omega]^{<\aleph_0},\bar{\gc} =
\la \gc_n \such n < \omega \r \in {\mathcal C},
\mup(\gc_n) \leq \mdn(\gc_{n+1}),\\
&\sup(w) <\mdn(\gc_0), \{ \mdn(\gc)_n \such n<\omega\}
\in D\}
\end{split}
\end{equation}
and order like in $\leq^{tr}_{\bf q}$ in the true ordering.
\end{definition}

We use the Definitions from Section~1 for this ${\bf q}_{\rm gr}$ and
for ${\bf Q}^{tr}_{{\bf q}_{\rm gr}}$.

\begin{lemma}\label{5.3}
${\bf q}_{\rm gr} \perp \Gamma_{\rm ud}$.
\end{lemma}

\proof This is a relative of the Main Lemma on page 263 in 
\cite{BsSh:257}. We show how to find our scenario in
Blass' and Shelah's work:
If 
\begin{myrules}

\item[(a)] $({\mathbb P}_1,{\mathbb P}_1^+, \name{\eta}) \in 
\Gamma_{\rm ud}$ and

\item[(b)] ${\mathbb P}_1 \lessdot {\mathbb P}_2$ and ${\mathbb P}_2 /{\mathbb P}_1$ is Cohen$_A$
for some set $A$ and ${\mathbb P}_2 \cap {\mathbb P}^+_1 = {\mathbb P}_1$,

\item[(c)] $\name{{\mathcal C}_2}$ is a ${\mathbb P}_2$-name and ${\mathbb P}_2$
forces that it is ${\bf q}_{\rm gr}$-nice,
\end{myrules}
then we can find $({\mathbb P}_2^+,f, \name{{\mathcal C}^+_2} )$ such that

\begin{myrules}

\item[($\alpha$)] ${\mathbb P}_1^+ \lessdot {\mathbb P}_2^+$,

\item[($\beta$)] $f$ is a complete embedding of ${\mathbb P}_2$ into ${\mathbb P}_2^+$ over
${\mathbb P}_1$,

\item[($\gamma$)] $\Vdash_{{\mathbb P}_2^+} \mbox{``} f({\name{\mathcal C}_2}) \subseteq
\name{{\mathcal C}_2}^+ \wedge \name{{\mathcal C}_2}^+ \mbox{ is } 
{\bf q}_{\rm gr}\mbox{-nice''}$,

\item[($\delta$)] $({\mathbb P}_1,{\mathbb P}_1^+,\name{\eta}) \leq_{\Gamma_{\rm ud}}
 (f({\mathbb P}_2),{\mathbb P}_2^+,\name{\eta}) \leq_{\Gamma_{\rm ud}}
 (f({\mathbb P}_2) *{\bbforcing}(f(\name{{\mathcal C}_2})),{\mathbb P}_2^+
* {\bbforcing}(\name{{\mathcal C}_2}^+),\name{\eta})$.
\end{myrules}

Indeed we can choose with a simple witness: 
${\mathbb P}_2^+ = ({\mathbb P}_1^+ \times_{{\mathbb P}_1} {\mathbb P}_2) 
\times {\rm Cohen}_\lambda$,
$\lambda = (|{\mathbb P}_1^+| + |{\mathbb P}_2|)^{\aleph_0}$, $f = id$.

Since Cohen reals do not add dominating functions,
$({\mathbb P}_1,{\mathbb P}_1^+,\name{\eta}) \leq_{\Gamma_{\rm ud}}
 (f({\mathbb P}_2),{\mathbb P}_2^+,\name{\eta})$.
Now for the second part we use the Main Lemma with
$M= V^{{\mathbb P}_2}$, $M'=V^{{\mathbb P}_2^+}$.
$Q({\mathcal U})$ has the r\^{o}le of 
$\bbforcing_{{\bf q}_{\rm gr}}({\mathcal C}_2)$ and
$Q({\mathcal U}')$  has the r\^{o}le of 
$ \bbforcing_{{\bf q}_{\rm gr}}({\mathcal C}^+_2)$.
The ${\bf q}_{\rm gr}$-niceness of ${\mathcal C}$ 
gives the desired maximality of ${\mathcal C}$ like being an ultrafilter.
We first look for ${\mathcal C}^+$ being ${\bf q}$-directed. 
If this were not possible, then finitely many elements would 
witness this. So this translates into Blass Shelah proof.
In the end we can extend the found ${\bf q}$-directed ${\mathcal C}^+$ 
to a ${\bf q}$-nice ${\mathcal C}^+$ using
Cohen reals as in Conclusion~\ref{1.17} Part 1.
\proofend

\proof \ref{5.1}:
We choose $\bar{\Gamma} = \la \Gamma_i \such i < \kappa \r$ so $\zeta(*) = \kappa$ and $i< \kappa 
\rightarrow \Gamma_i = \Gamma_{\rm ub}$. Using our present $\lambda$ we
take the forcing notion
${\mathfrak K}^* * \name{\mathbb P}$. It does not collapse any cardinals
 by \ref{3.6} and no cofinality 
is changed and clearly for some $(\iterq,\emptyset) \in
{\mathfrak K}^* * {\mathbb P}$ we have
$(\iterq,\emptyset) \Vdash_{{\mathfrak K}^* * {\mathbb P}} 2^{\aleph_0} = \lambda \, \wedge \,  MA_{<\kappa}$,
so $\gb \geq \kappa$ in the extension. 
Now by Claim~\ref{3.7}, for every $(\iterq,\emptyset)$,
$(\iterq,\emptyset) \Vdash_{{\mathfrak K}^* * {\mathbb P}} 
\{\name{\eta}_i \such i<\kappa \} $ is 
unbounded, hence $\gb = \kappa$.
By Claim~\ref{3.9} $\invgm({\bf q}_{\rm gr}) = \lambda$ and by Fact~\ref{4.4}
Part 2,
 $\Vdash_{{\mathfrak K}^* * {\mathbb P}}\gro \geq \lambda$, hence
since $2^\omega =\lambda$, $\gro = \lambda$.
\proofend

\nothing{
\section{Towards many cardinal invariants}\label{S6}
We would like to have many cardinal invariants with distinct values simultaneously.
We present frames for this, but the exact orthogonality conditions between 
the various ${\bf Q}$'s are not dealt with.  We 
 still have 
the problem of $\bf Q$'s interacting with itself if we like not just 
to control $\invgm{\bf Q}$, but even ${\rm Spec}(\bf Q)$.
We think to control them using forcing notion as in \S 3, but we use partial 
squares to control the various games on the various members, so for ${\rm Spec}(\bf Q)$ 
so the non-losing partial strategy is along a square sequence implicite in  \ref{3.8} 2.'s proof.
\begin{hypothesis}\label{6.1}
\begin{myrules}
\item[(a)] $\lambda = \lambda^{<\lambda} > \aleph_1$,
\item[(b)] $\leq^*$ is a well-ordering of ${\mathcal H}(\lambda)$,
\item[(c)] $A^*= \la A_\alpha^* \such \alpha < \lambda\r$, 
$A_\alpha^* \subseteq {\rm Reg} \cap \lambda \setminus \{\aleph_0\}$
\item[(d)] ${\bf F}$ is a function forcing
\begin{equation}
\begin{split}
\dom(\bf F) \subseteq & \{(\bar{\mathbb P},\name{\bf q}) \such 
\bar{\mathbb P} = \la {\mathbb P}_\beta \such \beta \leq\alpha \r \in {\mathcal H}(\lambda)\\
&
\mbox{ is a $\lessdot$-increasing sequence of c.c.c\ forcings},\\
& \name{\bf Q} \mbox{ is a canonical } {\mathbb P}_\alpha\mbox{-name as in 
Definition~\ref{1.1}. for } \beta <\alpha\}
\end{split}
\end{equation}
this includes the trivial $\name{\bf Q}$ it will be helpful to be idle sometimes.
${\bf F}(\bar{\mathbb P},\name{Q})$ is a $\bar{\mathbb P}_{\rm last}$-name of an ordinal 
$<\lambda$ such that if $(\bar{\mathbb P}^1,\name{\bf Q}^1),
(\bar{\mathbb P}^2,\name{\bf Q}^2) \in \dom(\bf F)$,
and if $\bar{\mathbb P}^1 \trianglelefteq \bar{\mathbb P}_\alpha^2$, 
$\alpha_\ell = \lg(\bar{{\mathbb P}}^\ell)-1$ then 
$$
\Vdash_{{\mathbb P}^2_{\alpha_2}} \mbox{``if } {\bf Q}^1 = {\bf Q}^2 \mbox{ or more then }
{\bf F}({\bar{\mathbb P}}^1,\name{\bf Q})
= {\bf F}(\bar{{\mathbb P}}^2, \name{\bf Q})\mbox{''}.$$
\end{myrules}
\end{hypothesis}
\begin{definition}\label{6.2}
We define a forcing nothin ${\mathbb  P}= {\mathbb P}_{\lambda, {\bf F}, 
\zeta(*)}$. We have $\iterp \in {\mathbb P}$ if 
$$
\iterp = \left\langle 
\alpha, \la {\mathbb P}_{\zeta,\beta} \such \zeta < \zeta(*), \beta \leq \alpha \r,
\la C_\beta \such \beta \leq  \alpha \r, F_1, F_2 \right\rangle,$$
such that
\begin{myrules}
\item[(a)] $\iterp \in {\mathcal H}(\lambda)$,
\item[(b)] ${\mathbb P}_{\zeta,\beta}$ is a forcing notion,
\item[(c)] for each $\beta \leq \alpha$, ${\mathbb P}_{\zeta,\beta}$ is increasing
with $\zeta$,
\item[(d)] $\la {\mathbb P}_{\zeta(*),\beta}\such \beta < \alpha \r$ is 
$\lessdot$-increasing continuous,
\item[(e)]
if $\gamma \leq \beta \leq \alpha$ then for every sufficiently large $\zeta < \zeta(*)$
${\mathbb P}_{\zeta,\gamma}\lessdot {\mathbb P}_{\zeta,\beta}$,
\item[(f)] if $\beta \leq \alpha$ is a limit ordinal then for some 
$\zeta < \zeta(*)$ and club $C$ of $\beta$'s for every $\zeta \in [\zeta_*,\zeta(*))$ 
the sequence $\la {\mathbb P}_{\zeta,\gamma} \gamma \in C \cup \{\beta\} \r$ is $\lessdot$-increasing continuous,
\item[(g)] $c_\beta \in \beta \setminus\{0\}$ $C_\beta$ is a closed subset of 
$\beta$, but possibly has a last element, and $\otp(C_\beta) \leq\kappa_\beta$,
\item[(h)] if $\gamma \in C_\beta$ then $C_\gamma = C_\beta \cap \gamma$,
\item[(i)] if $\gamma \in C_\beta$ and $\zeta \leq \zeta(*)$, then 
${\mathbb P}_{\zeta,\gamma} \lessdot {\mathbb P}_{\zeta,\beta}$,
\item[(j)] if $C_\beta = \sup(C_\beta)$ and $\zeta \leq \zeta(*)$, then 
${\mathbb P}_{\zeta,\gamma} \lessdot {\mathbb P}>_{\zeta,\beta}$,
\item[(k)] $\dom(F_\ell) \subseteq \alpha$, 
$F(\beta) \in \name{{bf Q}}$, $\name{\bf Q}$ is a 
${\mathbb P}_{\zeta(*),\alpha }$-name as in \ref{1.1} and we
 write 
$\name{\bf Q}_\beta$ for $F_1(\beta)$ moreover $\name{\bf Q}_\beta$ is a 
${\mathbb P}_{0,\beta}$-name,
\item[(l)] $F_2(\beta)$ is a regular cardinal $\in [\aleph_1,\lambda)$, 
and we write $\kappa_\beta$, $ \kappa_\beta(\iterq)$ for $F_2(\beta)$,
$\Vdash_{{\mathbb P}_{\zeta(*),\alpha}} \kappa_\beta \in {\bf F}(F_1(\beta))$, 
if $F_1(\beta)$ is well-defined,
\item[(m)] if $\gamma  \in C_\beta$ then $(F_1(\beta), F_2(\beta))=
(F_1(\gamma),F_2(\gamma))$, so if the right-hand side is well defined, then the left-hand
 side is defined,
\item[(n)] if $\beta \in C^-:= \{\beta \such \beta > \sup(C_\beta)\}$, then
$\name{\bar{\gc}}_\beta$, $\name{\bar{\gd}}_\beta$ are ${\mathbb P}_{0,\beta}$-names,
\item[(o)] if $\beta \in \dom(F_1)$, $\beta > \sup(C_\beta)$,  then $\la
\name{\bar{\gc}}_\gamma \name{\bar{\gb}}_\gamma \such \gamma \in {\rm nacc}(C_\beta) \cup
\{\beta\}) \r$
is a play in $\Game_{\kappa_{\beta}}(F_1(\beta))$,
\item[(p)] if $\beta \in C^-$ and $\gamma = \max(C_\beta)$ then 
$\name{\bar{\gc}}_\beta$ induces a ${\bf Q}_\gamma$-nice family 
${\mathcal C}_beta$ of ${\bf V}^{{\mathbb P}_\gamma}$
and ${\mathcal C}_\beta$ is a ${\mathbb P}_{\zeta(*),\gamma}$-name and 
${\mathcal C}_{\zeta,\beta}= {\mathcal C}_\beta =\name{\mathcal C}_\beta
\cap {\bf V}^{{\mathbb P}_{\zeta,\gamma}}$ is
a ${\mathbb P}_{\zeta,\gamma}$-name of a $\name{\bf Q}_\beta$-nice family,
\item[(q)] if $\beta <\alpha$, $\cf(\beta )= \theta^+$, 
then for stationarily many $\gamma < \beta$ we have $\cf(\gamma) = \theta \wedge \gamma = \sup(C_\gamma)$.
\end{myrules}
\end{definition}
Concluding Remark: The scheme should be clear, but we still have no good case where it is finished.
}

\bibliographystyle{plain}
\bibliography{../sh/lit,../sh/shelah}

\def\germ{\frak} \def\scr{\cal} \ifx\documentclass\undefinedcs
  \def\bf{\fam\bffam\tenbf}\def\rm{\fam0\tenrm}\fi 
  \def\defaultdefine#1#2{\expandafter\ifx\csname#1\endcsname\relax
  \expandafter\def\csname#1\endcsname{#2}\fi} \defaultdefine{Bbb}{\bf}
  \defaultdefine{frak}{\bf} \defaultdefine{mathfrak}{\frak}
  \defaultdefine{mathbb}{\bf} \defaultdefine{mathcal}{\cal}
  \defaultdefine{beth}{BETH}\defaultdefine{cal}{\bf} \def\bbfI{{\Bbb I}}
  \def\mbox{\hbox} \def\text{\hbox} \def\om{\omega} \def\Cal#1{{\bf #1}}
  \def\pcf{pcf} \defaultdefine{cf}{cf} \defaultdefine{reals}{{\Bbb R}}
  \defaultdefine{real}{{\Bbb R}} \def\restriction{{|}} \def\club{CLUB}
  \def\w{\omega} \def\exist{\exists} \def\se{{\germ se}} \def\bb{{\bf b}}
  \def\equivalence{\equiv} \let\lt< \let\gt> \def\cite#1{[#1]}
\begin{thebibliography}{10}

\bibitem{abraham:handbook}
Uri Abraham.
\newblock Proper forcing.
\newblock In Matthew Foreman, Akihiro Kanamori, and Menachem Magidor, editors,
  {\em Handbook of Set Theory}. Kluwer, To appear.

\bibitem{BJ}
Tomek Bartoszy\'{n}ski and Haim Judah.
\newblock {\em {{Set Theory, On the Structure of the Real Line}}}.
\newblock A K Peters, Wellesley, Massachusetts, 1995.

\bibitem{Blasstoronto}
Andreas Blass.
\newblock {Applications of superperfect forcing and its relatives}.
\newblock In Juris Stepr\={a}ns and Steve Watson, editors, {\em Set Theory and
  its Applications}, volume 1401 of {\em Lecture Notes in Mathematics}, pages
  18--40, 1989.

\bibitem{BlassShelah}
Andreas Blass and Saharon Shelah.
\newblock {There may be P$_{\aleph_1}$- and P$_{\aleph_2}$-points and the
  Rudin-Keisler ordering may be downward directed}.
\newblock {\em Ann. Pure Appl. Logic}, 33:213--243, [BsSh:242], 1987.

\bibitem{BsSh:257}
Andreas Blass and Saharon Shelah.
\newblock Ultrafilters with small generating sets.
\newblock {\em Israel Journal of Mathematics}, 65:259--271, 1989.

\bibitem{Goldstern93}
Martin Goldstern.
\newblock Tools for your forcing construction.
\newblock In Haim Judah, editor, {\em {Set Theory of the Reals}}, volume~6 of
  {\em Israel Mathematical Conferences Proceedings}, pages 305--360, 1993.

\bibitem{hindman:sums}
Neil Hindman.
\newblock Finite sums from sequences within cells of a partition of ${N}$.
\newblock {\em J. Combin. Theory Ser. A}, 17:1--11, 1974.

\bibitem{JudahRoslanowski}
Haim Judah and Andrzej Ros{\l}anowski.
\newblock On {S}helah's amalgamation.
\newblock In Haim Judah, editor, {\em {Set Theory of the Reals}}, volume~6 of
  {\em Israel Mathematical Conferences Proceedings}, pages 385--415, 1993.

\bibitem{RoSh:470}
Andrzej Ros{\l}anowski and Saharon Shelah.
\newblock {\em {Norms on Possibilities I: Forcing with Trees and Creatures}},
  volume 141 (no. 671) of {\em Memoirs of the American Mathematical Society}.
\newblock 1999, [RoSh:470].

\bibitem{Sh:700}
Saharon Shelah.
\newblock Are $\mathfrak a$ and $\mathfrak d$ your cup of tea? [{S}h:700].
\newblock {\em Acta Mathematica}, accepted.

\bibitem{Sh:707}
Saharon Shelah.
\newblock Long iterations for the continuum, [{S}h:707].
\newblock {\em Archive for Mathematical Logic}, submitted.

\bibitem{Sh:509}
Saharon Shelah.
\newblock Vive la diff\'erence {III}.
\newblock {\em Israel Journal of Mathematics}, submitted, [Sh:509].

\bibitem{Sh:207}
Saharon Shelah.
\newblock On cardinal invariants of the continuum.
\newblock In {\em Axiomatic set theory (Boulder, Colo., 1983)}, volume~31 of
  {\em Contemp. Mathematics}, pages 183--207, [Sh:207]. Amer. Math. Soc.,
  Providence, RI, 1984.
\newblock Proceedings of the Conference in Set Theory, Boulder, June 1983; ed.
  Baumgartner J., Martin, D. and Shelah, S.

\bibitem{Sh:f}
Saharon Shelah.
\newblock {\em {Proper and Improper Forcing, 2nd Edition}}.
\newblock Springer, 1998.

\end{thebibliography}
\end{document}